\documentclass[11pt,reqno]{amsart}
\usepackage{amsbsy,amsfonts,amsmath,amssymb,amscd,amsthm}
\usepackage{mathrsfs,float}
\usepackage[mathcal]{euscript}
\usepackage{subfigure,enumitem,graphicx,epstopdf}
\usepackage[a4paper,text={6.1in,8in},centering]{geometry}
\usepackage{pxfonts,epstopdf}
\usepackage[colorlinks=true,linkcolor=blue,citecolor=green]{hyperref}
\usepackage{srcltx}
\usepackage[margin=20pt,font=small,labelfont=bf,textfont=it]{caption}

\hypersetup{linktocpage}

\small\normalsize
\theoremstyle{plain}
\newtheorem{mainthm}{Theorem}

\newtheorem{thm}{Theorem}[section]
\newtheorem{cor}[thm]{Corollary}
\newtheorem{lem}[thm]{Lemma}
\newtheorem{claim}[thm]{Claim}
\newtheorem{prop}[thm]{Proposition}
\newtheorem{addendum}[thm]{Addendum}
\newtheorem{defi}[thm]{Definition}
\newtheorem{scho}[thm]{Scholium}
\theoremstyle{definition}
\newtheorem{exap}[thm]{Example}
\newtheorem{rem}[thm]{Remark}

\newcommand{\eqdef}{\stackrel{\scriptscriptstyle\rm def}{=}}
\DeclareMathOperator{\IFS}{IFS}

\makeatletter
\def\l@part{\@tocline{0}{-2pt}{1pc}{}{}}
\def\l@section{\@tocline{1}{-2pt}{1pc}{4.6em}{}}
\renewcommand{\tocpart}[3]{%
  \indentlabel{\@ifnotempty{#2}{\makebox[2.3em][l]{%
    \ignorespaces#1 #2.\hfill}}}\bf{#3}}
\renewcommand{\tocsection}[3]{%
  \indentlabel{\@ifnotempty{#2}{\hspace*{2.3em}\makebox[2.3em][l]{%
    \ignorespaces#1 #2.\hfill}}}#3}

\makeatother

\begin{document}
\title[Robust  tangencies]{Robust tangencies of large codimension}
\author[Barrientos]{Pablo G. Barrientos}
\address{\centerline{Instituto de Matem\'atica e Estat\'istica, UFF}
    \centerline{Rua M\'ario Santos Braga s/n - Campus Valonguinhos, Niter\'oi,  Brazil}}
\email{pgbarrientos@id.uff.br}
\author[Raibekas]{Artem Raibekas}
\address{\centerline{Instituto de Matem\'atica e Estat\'istica, UFF}
   \centerline{Rua M\'ario Santos Braga s/n - Campus Valonguinhos, Niter\'oi,  Brazil}}
\email{artem@mat.uff.br}

\subjclass[2010]{58F15, 58F17, 53C35.}
\keywords{homoclinic tangencies, heterodimensional tangencies, bundle tangencies, blenders}
\thanks{The first author was supported by MTM2014-56953-P and the second  by FAPERJ-INST 111.595/2014.}

\begin{abstract}
We construct $C^2$-robust homoclinic and heterodimensional tangencies of large codimension
inside transitive partially hyperbolic sets.
\end{abstract}

 \maketitle \thispagestyle{empty}


\section{Robust cycles, tangencies and transitivity in dynamical systems}
\label{s:dynamics} The study of uniformly hyperbolic dynamics
began with the works of Anosov and Smale in the 60s, and nowadays
these systems are well understood from the topological and
statistical perspectives (see~\cite{Kat95}). Although uniform
hyperbolicity was first believed to involve a dense subset of
$C^r$-diffeo\-mor\-phisms of a compact manifold, it soon emerged
that this was not true~\cite{AS70,New70,Si72}. There are two
important mechanisms that yield
\emph{non-hyperbolic dynamics}: heterodimensional cycles and
homoclinic tangencies. A
diffeomorphism $f$ has a \emph{heterodimensional cycle} associated
with two transitive hyperbolic sets if these sets have different
indices (dimension of the stable bundle) and their invariant
manifolds meet cyclically. On the other hand, $f$ has a
\emph{homoclinic tangency} if there is a pair of points $P$ and
$Q$ belonging to the same transitive hyperbolic set so that the unstable
invariant manifold of $P$ and the stable invariant manifold of $Q$
have a non-transverse intersection $Y$. The positive integer $
c_T= \dim T_Y W^u(P)\cap T_Y W^s(Q)$ is called the
\emph{codimension of the tangency} (at $Y$) and we say that the
tangency is \emph{large} (or degenerated) if $c_T\geq 2$.

A well-known and important conjecture by Palis \cite{P00} claims
that these
cycles 
  are $C^r$-dense in the complement of hyperbolic systems. In
contrast with hyperbolic systems,  heterodimensional cycles and
homoclinic tangencies are in general not robust. However,
important examples of open subsets of non-hyperbolic
diffeomorphisms arise from \emph{robust cycles}, i.e., from
heterodimensional cycles or homoclinic tangencies that persist
under perturbations.
Namely,
a $C^r$-diffeomorphism $f$ has a

\begin{enumerate}[leftmargin=0.6cm, itemsep=0.2cm]
\item[-] \emph{$C^r$-robust heterodimensional cycle} if there are
transitive hyperbolic sets $\Lambda_1$ and $\Lambda_2$ of $f$ with
different indices and a $C^r$-neighborhood $\mathscr{U}$ of $f$ such
that any $g\in \mathscr{U}$ has a heterodimensional cycle
associated with the continuations of $\Lambda_1$ and $\Lambda_2$
for $g$ respectively;
\item[-] \emph{$C^r$-robust homoclinic tangency} (of codimension $c_T>0$) if
there are a transitive hyperbolic set $\Lambda$  and a
$C^r$-neighbor\-hood $\mathscr{U}$ of $f$ such that any $g\in
\mathscr{U}$ has a homoclinic tangency (of codimension $c_T>0$)
associated with the continuation of $\Lambda$ for $g$.
\end{enumerate}


In~\cite{New70,New79}, Newhouse constructed surface
diffeomorphisms having $C^2$-robust homoclinic tan\-gen\-cies
and these results were extended to higher dimensions
in~\cite{GTS93,PV94,Ro95}.
\mbox{Also in~\cite{Ao08,BD12}}
robust homoclinic tangencies are built in the $C^1$-topology and
dimension higher than $2$. In dimension 2, results by Moreira in
\cite{Moreira2011} imply that there are no $C^1$-robust homoclinic
tangencies associated with hyperbolic basic sets. However, in all
these constructions the homoclinic tangencies have codimension
$c_T=1$.

The study of bifurcations involving heterodimensional cycles
leading to robustly non-hyperbolic dynamics was started in the
pioneering papers by D\'iaz~\cite{D95a,D95b}. New examples of
robustly non-hyperbolic transitive diffeomorphisms were
constructed by Bonatti and D\'iaz in~\cite{BD96}.  The
non-hyperbolicity of these examples
comes from the existence of a robust heterodimensional cycle and, in particular, from
the coexistence of periodic points with different indices inside the same transitive set.
The construction of Bonatti and D\'iaz can be thought as a
generalization of previous examples due to Shub and
Ma\~ne~\cite{Sh71,Ma78} using a better understanding of the
mechanism that provides the robustness of  
cycles and transitivity. They named this local plug, which is a
special transitive hyperbolic set, as \emph{blender}.

A strong version of the Palis conjecture~\cite{BD12} states that
every diffeomorphism can be approximated either by a hyperbolic
one or by one with a $C^1$-robust cycle. In view of this, it is
important to understand and describe mechanisms generating robust
cycles. This is one of the contributions of this paper concerned
with two main issues: robust heterodimensional cycles of any
co-index (absolute value of the difference of indices of the
hyperbolic sets cyclically connected) and robust tangencies of any
codimension.
Namely, we deal with the \emph{construction} of
$C^2$-robust tangencies (of large codimension) inside transitive
partially hyperbolic sets  in dimension $d\geq 4$. By a partially
hyperbolic set we understand a closed invariant set whose tangent
bundle admits a dominated splitting into three non-trivial
continuous vector subbundles $E^s \oplus E^c \oplus E^u$ so that
$E^s$ is uniformly contracted and $E^u$ is uniformly expanded.

The derivative of the map naturally induces a cocycle on the
tangent bundle. Robust heterodimensional cycles constructed for
this induced map on the tangent bundle will ''project'' to robust
tangencies for the original diffeomorphism (see Section 2 for a
geometrical outline of this idea). To obtain the robust cycle for
the tangent bundle dynamics (or namely in \emph{Grassmannian
bundle manifolds}), we will use methods similar to those presented
in~\cite{BKR14}. We construct blenders (with a higher dimensional
center direction) for the induces dynamics which will require the
induced map to be $C^1$, and thus the original diffeomorphism to
be $C^2$. In contrast, the known  examples in higher dimensions
with robust tangencies are based on either a dimension reduction
using normally-hyperbolic manifolds or on the construction of
blenders (with one-dimensional center direction) in the
\emph{ambient manifold}. In addition, in the recent works of
Berger~\cite{Ber16,Ber16b} robust tangencies were constructed for
open sets of parametric families such that the order of  the
unfolding of the tangency with respect to the parameter is
controlled. For this purpose blenders were created in the
\emph{space of jets}. We observe that these are still codimension
one tangencies. Our objective in this article is the construction
of higher codimensional tangencies and this involves different
ideas such as the dynamics in the \emph{tangent bundle}.

\subsection{Homoclinic tangencies}
In the following theorem we obtain robust homoclinic tangencies of
large codimension. We observe that the manifold is not necessarily
compact, and in this case the set of diffeomorphisms is endowed
with the compact open $C^r$-topology.

\begin{mainthm}
\label{thm:cor1} Every manifold of dimension $d\geq 4$ admits a
diffeomorphism with a transitive partially hyperbolic set having a
$C^2$-robust homoclinic tangency of codimension $c_T$ which can be
chosen to be any integer $0<c_T\leq \lfloor (d-2)/2 \rfloor$.
\end{mainthm}

Tangencies in the above result can have large codimension for $d
\geq 6$. Notice that in general large homoclinic tangencies
require $d\geq 4$. However, if we restrict the dynamics to
invariant partially hyperbolic sets, then the above result, on the
existence of robust tangencies of large codimension, is actually
optimal with respect to the dimension of the manifold. This is
because partially hyperbolic sets would require at least two extra
dimensions for the stable ($E^s$) and the unstable ($E^u$)
directions. Still, the problem of the existence of large
tangencies for lower dimensional manifolds (and far from partially
hyperbolic sets) remains open.

\subsection{Heterodimensional tangencies} The study of non-transverse intersections between
stable and unstable manifolds of two hyperbolic sets with
different indices led to the notion of a heterodimensional
tangency, formally defined below. The dynamical consequences of
heterodimensional tangencies were first studied in~\cite{DNP06}.
But the first examples of robust heterodimensional tangencies
associated with a heterodimensional cycle were given later on by
Kiriki and Soma in~\cite{KS12}.

A $C^r$-diffeomorphism $f$ of a manifold $\mathcal{M}$  has a
\emph{heterodimensional tangency} if there are tran\-sitive
hyperbolic sets $\Lambda_1$ and $\Lambda_2$, points $P\in
\Lambda_1$, $Q\in \Lambda_2$ and $Y\in W^u(P)\cap
W^s(Q)$~such~that
$$
 \dim T_Y W^u(P) + \dim T_Y W^s(Q)  > \dim \mathcal{M} \quad \text{and}
 \quad
 T_Y \mathcal{M} \not = T_Y W^u(P) + T_Y W^s(Q).
$$
Observe that the above condition implies $i_1<i_2$ where $i_1$ and
$i_2$ are the indices of $\Lambda_1$ and $\Lambda_2$ respectively.
We also call the \emph{codimension of the tangency} (at $Y$) the positive integer
\begin{align*}
         c_T &= \dim \mathcal{M} - \big[\dim T_Y W^u(P) + \dim T_Y W^s(Q) -
         \dim T_Y W^u(P)\cap T_YW^s(Q)\big]. 
\end{align*}
The codimension measures how far the tangencial intersection is
from a transverse intersection. The above heterodimensional
tangency is said to be \emph{$C^r$-robust} (of codimension
$c_T>0$) associated with $\Lambda_1$ and $\Lambda_2$ if every
small enough $C^r$-perturbation $g$ of $f$ has an
heterodimensional tangency (of codimension $c_T$) associated with
the continuations of $\Lambda_1$ and~$\Lambda_2$~for~$g$.

The examples in~\cite{KS12} are $C^2$-robust heterodimensional
tangencies of codimension $c_T=1$. The following result provides
the first example of robust heterodimensional tangencies of large
codimension. We need first some notation. Let us consider a
manifold $M$ of dimension $c\geq 3$ and a $C^r$-diffeomorphism $F$
of a manifold $N$  having a hyperbolic set $\Lambda \subset N$ of
index $i_F$  conjugated to a full shift with a sufficiently large
number of symbols that will depend only on the dimension of $M$.
Both of the manifolds, $M$ and $N$, are not necessarily compact.

\begin{mainthm}
\label{thmC} For any pair of integers $0< k < c-1$ and
$0<i_{cs}<c-k$, there is an arc
$\{f_{\varepsilon}\}_{\varepsilon\geq 0}$ of $C^r$-diffeomorphisms
of $N \times M$ such that $f_{0} = F \times \mathrm{id}$ and for
every $\varepsilon
>0$, any small enough $C^2$-perturbation $g$ of $f_{\varepsilon}$ has
\begin{enumerate}[leftmargin=0.6cm]
\item[-] a transitive partially hyperbolic set $\Delta_g \subset N \times M$ homeomorphic to $\Lambda \times
M$ and
\item[-] a 
heterodimensional tangency (in $\Delta_g$) between basic
sets of indices $i_F+i_{cs}$ and $i_F+i_{cs}+k$.
\end{enumerate}
The codimension of the tangency can be chosen to be any integer $
0<c_T\leq \min\{i_{cs},c-i_{cs}-k\}$. Thus we can get
tangencies of codimension $c_T=\lfloor(c -1)/2\rfloor$ which are
large for $c \geq 5$.
\end{mainthm}

The statement of the above theorem emphasizes the skew-product
construction of the maps and the fact that the tangencies can be
embedded in a partially hyperbolic invariant set. But tangencies
are by themselves local objects and can be created in any
manifold (it is enough to take $N$ and $M$ of the previous theorem as local charts).
These considerations imply as a consequence a similar
result as Theorem~\ref{thm:cor1} for heterodimensional tangencies.
In particular, we get that any manifold of dimension $d\geq 5$
admits a diffeomorphism with a transitive partially hyperbolic set
having inside a $C^2$-robust heterodimensional tangency of
codimension $c_T=\lfloor (d-3)/2 \rfloor$. Similarly as homoclinic
tangencies, heterodimensional tangencies of large codimension need
the dimension to be at least $5$. This can be checked by
considerations on the co-index $k$ of the transitive hyperbolic
sets involved since $c_T = \dim T_YW^u(P)\cap T_Y W^s(Q) - k$.
However, once again notice that to get a large heterodimensional
tangency inside a partially hyperbolic set requires  $d\geq
7$.

\subsection{Bundle tangencies}
\label{s:tangencies} The non-hyperbolicity of a transitive set
with homoclinic tangencies or heterodimensional cycles comes, in
particular, from the existence of a direction in the tangent space
that is exponentially contracted for forward and backward
iterations. Let $f$ be a diffeomorphism of a manifold
$\mathcal{M}$ and consider a point $x \in \mathcal{M}$. We call a
unitary vector $v\in T_x \mathcal{M}$ a \emph{tangent direction}
if there exist constants $C>0$ and $0<\lambda<1$ so that
$$
   \|Df^{n}(x)v\|\leq C\lambda^{|n|} \quad \text{for all } n\in
   \mathbb{Z}.
$$

Recently, in an announcement~\cite{G14}, Gourmelon introduced the
notion of a bundle tangency, which allows to unify the different
types of tangencies: homoclinic and heterodimensional. Namely, a
$C^r$-diffeomorphism $f$ has a \emph{bundle tangency} between the
unstable manifold $W^u(\Lambda_1)$ and the stable manifold
$W^s(\Lambda_2)$ of transitive hyperbolic sets $\Lambda_1$ and
$\Lambda_2$ respectively if there are points $P\in \Lambda_1$,
$Q\in \Lambda_2$ and $Y \in W^u(P)\cap W^s(Q)$ such that
\begin{equation}
\label{eq:bandle1} d_T \eqdef \dim T_Y W^u(P)\cap T_Y W^s(Q) >0
\quad \text{and} \quad T_Y \mathcal{M} \not = T_Y W^u(P)+T_Y
W^s(Q).
\end{equation}
The integer $d_T$ is called the \emph{dimension} of the bundle
tangency (at $Y$). Notice that $d_T$ is the maximum number of
independent tangent directions in $T_Y \mathcal{M}$. Note
that~\eqref{eq:bandle1} is equivalent to
$$
d_T>0 \quad \text{and} \quad
\mathrm{ind}^u(\Lambda_1)+\mathrm{ind}^s(\Lambda_2)-d_T < \dim
\mathcal{M}.
$$
This forces that
$d_T > \max\{0, i_\alpha - i_\omega\}$ 
where $i_\alpha=\mathrm{ind}^s(\Lambda_2)$
   and  $i_\omega=\mathrm{ind}^s(\Lambda_1)$.

The \emph{codimension} of the bundle tangency is defined by
$$
   c_T \eqdef \dim \mathcal{M} -
   [\mathrm{ind}^u(\Lambda_1)+\mathrm{ind}^s(\Lambda_2)-d_T] = d_T
   - (i_\alpha-i_\omega).
$$

Homoclinic and heterodimensional tangencies are particular cases
of bundle tangencies, where for homoclinic tangencies one just
needs to take $\Lambda_1=\Lambda_2$. In fact, all these cases
require that
$\mathrm{ind}^{u}(\Lambda_1)+\mathrm{ind}^s(\Lambda_2)\geq \dim
\mathcal{M}$. Then it must hold that $i_\omega\leq i_\alpha$ and
hence the codimension of tangency is given by $c_T=d_T -k$ where
$k=i_\alpha-i_\omega\geq 0$ is the co-index between $\Lambda_1$
and $\Lambda_2$. Having in mind this relation between the
dimension of the tangency and the co-index of the involved
hyperbolic sets, the robust tangencies of Theorems~\ref{thm:cor1}
and~\ref{thmC} follow from the next result.

\begin{mainthm}
\label{thmF} For any integer $0<i_1,i_2 < c$  there is an arc
$\{f_{\varepsilon}\}_{\varepsilon \geq 0}$ of
$C^r$-diffeomorphisms of $N\times M$ such that $f_{0} = F \times
\mathrm{id}$ and for every $\varepsilon
>0$, $f_{\varepsilon}$ has a
\begin{enumerate}[leftmargin=0.6cm]
\item[-] $C^2$-robust bundle tangency  between the unstable and the
stable manifold of basic sets $\Gamma^1_\varepsilon$ and
$\Gamma^2_\varepsilon$ contained in $\Lambda\times M$ of indices
$i_\alpha=i_F+i_1$ and $i_\omega=i_F+i_2$ respectively.
\end{enumerate}
The dimension $d_T$ of the tangency can be chosen to be any
integer so that
$$
\max\{0,i_2-i_1\}<d_T\leq \min\{c-i_1,i_2\}.
$$
Moreover, these arcs can be taken in such a way that
$f_\varepsilon$ also has a $C^2$-robust bundle tangency between
the stable manifold of $\Gamma^1_\varepsilon$ and the unstable of
$\Gamma^2_\varepsilon$.
\end{mainthm}

\subsection{Symbolic skew-products}
Besides the construction of robust cycles for partially hyperbolic
diffeomorphisms, we introduce blenders and (robust) cycles for
bi-Lipschitz skew-product homeomorphisms. Namely we work with
skew-products over a shift of finite symbols and whose fiber maps
depend H\"older continuously on the base points, that we call
\emph{symbolic skew-products}. The above theorems are, in fact, a
consequence of the realization, via a straightforward application
of known results \cite{GI99,Go06,IN10,BKR14}, of a more general
construction in this setting (see
Theorem~\ref{thm:mainA-symbolic-setting} and
Theorem~\ref{thm:mainC-symbolic-setting}). Moreover, we give
\emph{criteria} to construct robust heterodimensional cycles
(Theorem~\ref{thm-cycles1}) and robust tangencies
(Theorem~\ref{thm:tangencias}).


\subsection{Structure of the paper}
In the next section, we explain the geometric ideas behind the
construction of blenders and our new criterion for robust
tangencies of large codimension. These constructions will actually
be made in the abstract set of the so-called \textit{symbolic skew
products}. Thus, section $3$ gives the necessary definitions and
properties of symbolic skew-products and defines partial
hyperbolicity and blenders in this context. Section $4$ provides
the criteria for constructing robust cycles and tangencies in
symbolic skew-products. In section $5$, it is explained how to
build an actual example satisfying the previous criteria, but
still in the symbolic setting. Finally, in section $6$, it is
shown how to transfer the examples from the symbolic setting to a
realization in smooth manifolds.

\section{Outline of the method to yield robust tangencies of large codimension}
\label{s:pictures}
In this section we explain the 
ideas behind the
construction of robust tangencies of large codimension. To
accomplish this task, we need first to introduce the class of
blenders that we will consider and explain how these local
mechanisms are used to construct robust cycles.

 \begin{figure}
 \footnotesize
  \begin{center}
   \begin{picture}(400,145)
\put(10,0){\subfigure[projection onto $E^{ss}\oplus E^{cs}\oplus
E^{uu}$]{\label{fig1a}
\includegraphics[scale=0.5]{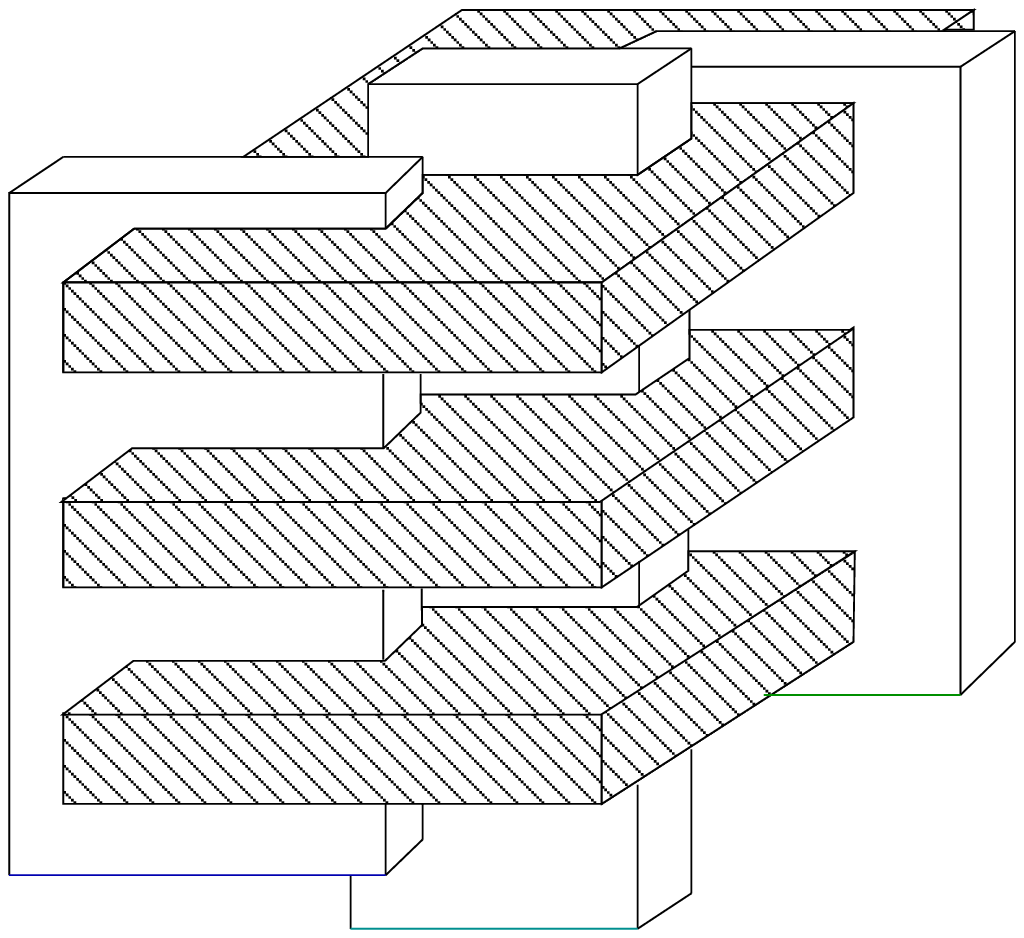}}}
\put(26,83){$\mathcal{H}_1$} \put(26,52){$\mathcal{H}_2$}
\put(26,21){$\mathcal{H}_3$} \put(5,117){$f^{n_1}(\mathcal{H}_1)$}
\put(115,10){$f^{n_2}(\mathcal{H}_2)$}
\put(160,125){$f^{n_3}(\mathcal{H}_3)$}
\put(230,0){\subfigure[projection onto $E^{ss}\oplus E^{cu}\oplus
E^{uu}$]{\label{fig1b}\includegraphics[scale=0.5]{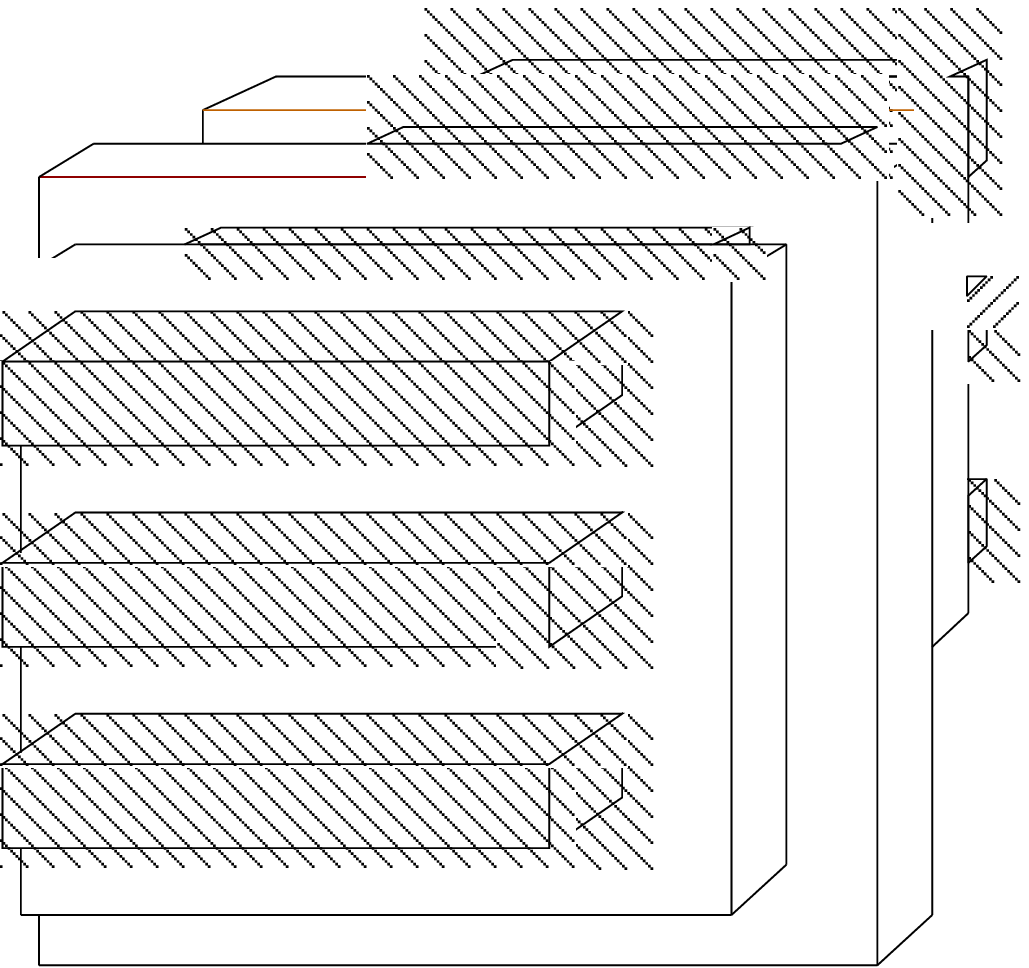}}}
\put(236,78){$\mathcal{H}_1$} \put(236,49){$\mathcal{H}_2$}
\put(236,20){$\mathcal{H}_3$}
\put(195,98){$f^{n_1}(\mathcal{H}_1)$}
\put(369,10){$f^{n_2}(\mathcal{H}_2)$}
\put(224,130){$f^{n_2}(\mathcal{H}_3)$}
 \end{picture}
~\\[1.2cm]
\begin{picture}(250,120)
\put(-20,90){\vector(1,0){30}} \put(-37,110){$E^{uu}$}
\put(-20,90){\vector(0,1){30}} \put(13,87){$E^{cs}$}
\put(-20,90){\vector(-1,-1){18}} \put(-50,60){$E^{ss}$}
\put(50,20){\subfigure[projection onto $E^{cs}\oplus
E^{cu}$]{\label{fig1c}\includegraphics[scale=0.7]{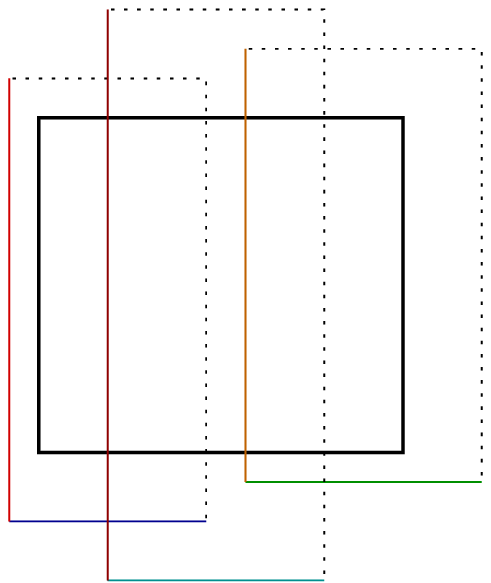}}}
\put(35,40){$f^{n_1}(\mathcal{H}_1)$}
\put(135,25){$f^{n_2}(\mathcal{H}_2)$}
\put(165,90){$f^{n_3}(\mathcal{H}_3)$} \put(76,102){$\mathcal{B}$}
\put(240,90){\vector(1,0){30}}
\put(223,110){$E^{ss}$}
\put(240,90){\vector(0,1){30}}
\put(275,87){$E^{cu}$}
\put(240,90){\vector(-1,-1){18}}
\put(210,60){$E^{uu}$}
\end{picture}
 \caption{This figure shows the case of a $cs$-blender in dimension $4$ with $\dim E^{cs}>0$, $\dim E^{cu}>0$.
 Figure~\ref{fig1a}
 represents the superposition domain $\mathcal{B}=\mathcal{H}_1\cup\mathcal{H}_2\cup \mathcal{H}_3$ and $f^{n_1}(\mathcal{H}_1)\cup f^{n_2}(\mathcal{H}_2) \cup f^{n_3}(\mathcal{H}_3)$
 in the section $E^{ss}\oplus E^{cs}\oplus
 E^{uu}$ and Figure~\ref{fig1b} in the section $E^{ss}\oplus E^{cu}\oplus
 E^{uu}$. Figure~\ref{fig1c} is the projection onto the central direction, $E^c=E^{cs}\oplus E^{cu}$, showing the corresponding
 covering property.
 }\label{fig1}
 \end{center}
 \end{figure}

\subsection{Blenders}
In our constructions, we will consider a particular class of
blenders obtained by the covering property criterion:

 A compact invariant set $\Gamma$ of a diffeomorphism $f$ of
a manifold $\mathcal{M}$ is a \emph{$cs$-blender} if
\begin{enumerate}[itemsep=0.2cm]
\item there is an open neighborhood $\mathcal{U}$ of $\Gamma$ such that $\Gamma$ is the maximal invariant set in $\overline{\mathcal{U}}$,
\item the set $\Gamma$ is transitive, (uniformly) hyperbolic with a dominated
splitting
$$T_\Gamma \mathcal{M} = E^{ss}\oplus E^c \oplus E^{uu}, \quad \text{where} \ E^c=E^{cs}\oplus E^{cu} \ \ \text{with} \ \ \dim E^{cs}>0 \ \ \text{and}$$
$E^{ss}\oplus E^{cs}$ and $E^{cu}\oplus E^{uu}$ being the
contracting and the expanding bundle respectively;
\item there are an open set $\mathcal{B} \subset \mathcal{U}$,
called the \emph{superposition domain},
 and integers
$n_1, \dots, n_k \in \mathbb{N}$ such that
$f^{n_1}(\mathcal{B})\cup \dots \cup f^{n_k}(\mathcal{B})$
intersects $\mathcal{B}$ as is shown in Figure~\ref{fig1} and explained below.
\end{enumerate}




The superposition domain in Figure~\ref{fig1} is
the set $\mathcal{B}=\mathcal{H}_1\cup \mathcal{H}_2 \cup
\mathcal{H}_3$ where we write $\mathcal{H}_i= \mathsf{H}_i\times
B$ so that $\mathsf{H}_i$ is an open ''rectangle`` in
$E^{ss}\oplus E^{uu}$ and $B$ is an open set in $E^c$.
Let $B_i$ be the projections of $f^{n_i}(\mathcal{H}_i)$
onto~$E^c$. The \emph{covering property} is defined by the
condition
$$
\overline{B} \subset \bigcup_{i} B_i,
$$
where the Lebesgue number of the cover is large enough in relation
to the transverse variation of the strong stable leafs (see
Theorem~\ref{cs-thm}). We will refer to $B$ as the \emph{blending
region}.

\begin{figure}
\vspace{0.2cm} \footnotesize
\begin{center}
\begin{picture}(500,290)
\put(30,0){\includegraphics{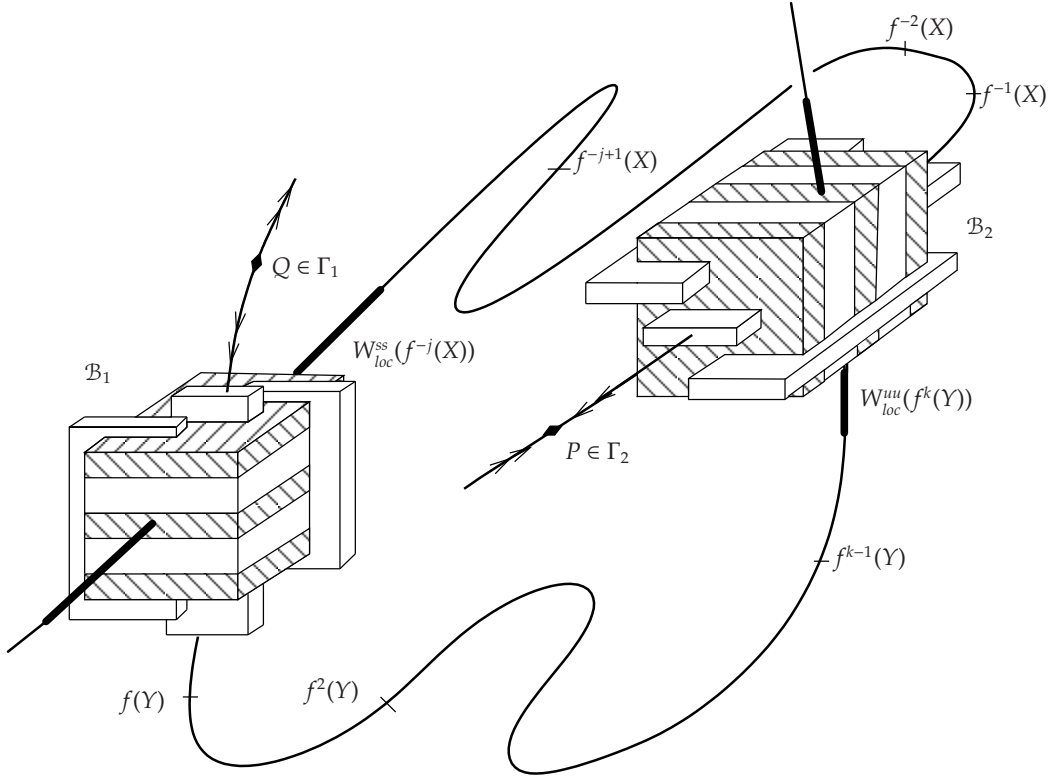}}
\put(240,120){$P\in\Gamma_2$} \put(395,255){$f^{-1}(X)$}
\put(360,280){$f^{-2}(X)$} \put(243,230){$f^{-j+1}(X)$}
\put(160,158){$W^{ss}_{loc}(f^{-j}(X))$}
\put(130,190){$Q\in\Gamma_1$} \put(73,26){$f(Y)$}
\put(142,30){$f^2(Y)$} \put(340,80){$f^{k-1}(Y)$}
\put(350,140){$W^{uu}_{loc}(f^{k}(Y))$}
\put(60,150){$\mathcal{B}_1$} \put(390,205){$\mathcal{B}_2$}
\end{picture}
\caption{Criterion for robust cycles}\label{fig4}
\end{center}
\end{figure}

Let $\Lambda^u=\Lambda^u(\mathcal{B};f)$ be the set of points $P \in \mathcal{B}$ such that
$f^{-n_{i_j}}\circ \dots \circ f^{-n_{i_1}}(P)\in \mathcal{B}$ 
for some sequence of integers $n_{i_j} \in \{n_1,\dots,n_k\}$.
Observe that when all of $n_i=1$ for all $i=1,\dots,k$, then
$\Lambda^u$ is simply the maximal forward invariant set in
$\mathcal{B}$, and furthermore it holds that
$$\Lambda^u \subset W^{u}_{loc}(\Gamma) \eqdef
  \{P\in\mathcal{M}: \ f^{-n}(P)\in \mathcal{U}  \  \text{for $n\geq 0$}\}
  =\bigcap_{n\geq 0}f^n(\mathcal{U}) .
$$ The main property of a superposition domain $\mathcal{B}$ of
a $cs$-blender $\Gamma$ for $f$ is the following:
 \begin{equation}
 \label{(B)}
\Lambda^u \cap W^{ss}_{loc}(P) \not=\emptyset \quad \textit{for all $P \in \mathcal{B}$.}
\end{equation}
Indeed, the covering property implies the intersection property:
given any $P\in \mathcal{B}$, by the covering property (see
figure~\ref{fig1a}), the local strong stable manifold
$\mathcal{D}^s=W^{ss}_{loc}(P)$ of $P$
goes through
some $f^{n_{i_1}}(\mathcal{H}_{i_1})$.  Then the pre-image
$f^{-n_{i_1}}(\mathcal{D}^s)$ contains a new disc
$\mathcal{D}_1^s$, which is the strong stable manifold of another point in $\mathcal{B}$. 
This allows us to repeat the process inductively and get a
sequence of discs
$$\mathcal{D}_j^s
\subset f^{-n_{i_j}}\circ \dots \circ f^{-n_{i_1}}(\mathcal{D}^s) \quad \text{in} \ \  \mathcal{B}.
$$ Hence,
by the hyperbolicity of $f$ on $\mathcal{B}$, the nested sequence
$f^{n_{i_1}}\circ\dots \circ f^{n_{i_j}}(\mathcal{D}^s_{j})
\subset \mathcal{D}^s$ defines a point in $ \Lambda^u\cap
W^{ss}_{loc}(P)$.  Moreover, since the covering property of a
$cs$-blender persists under $C^1$-perturbations of $f$, the above
property~\eqref{(B)} is $C^1$-robust.


\begin{rem}[$cu$-blender and $double$-blender]
\emph{A $cu$-blender (for $f$) is a $cs$-blender for $f^{-1}$.}
Similarly, the $cu$-covering property implies that $C^1$-robustly
\begin{equation}
 \label{(B2)}
\Lambda^s \cap W^{uu}_{loc}(P) \not=\emptyset \quad \textit{for all $P \in \mathcal{B}$
\quad where $\Lambda^s=\Lambda^s(\mathcal{B};f)\eqdef\Lambda^u(\mathcal{B};f^{-1})$.}
\end{equation}
\emph{We say that $\Gamma$ is a $double$-blender if it is simultaneously
both, a $cs$-blender and a $cu$-blender for $f$.}
\end{rem}


\subsection{Criterion for robust cycles}

\begin{figure}
 \footnotesize
\begin{center}
\begin{picture}(300,180)
\put(30,0){\includegraphics[scale=0.7]{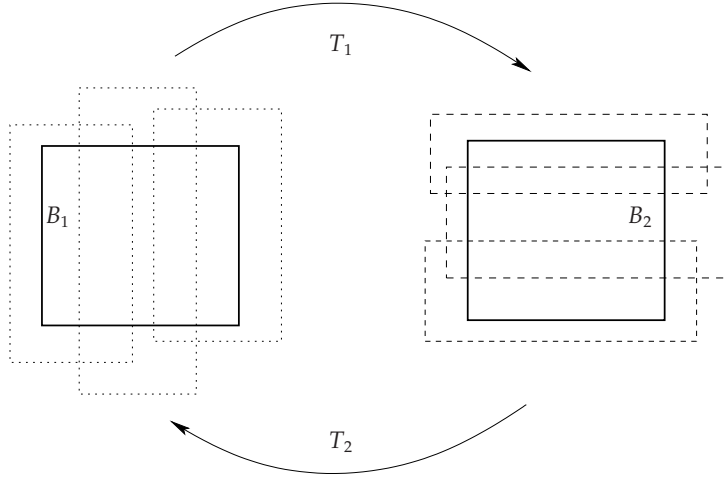}}
\put(150,160){$T_1$} \put(150,10){$T_2$} \put(44,95){$B_1$}
\put(262,95){$B_2$}
\end{picture}
 \caption{Robust cycle between a $cs$ and $cu$-blender
as seen on $E^c$. The sets $B_1\equiv \mathcal{B}_1$ and
$B_2\equiv \mathcal{B}_2$ are, respectively, a $cs$ and
$cu$-blending region. The maps $T_1\equiv f^j$ and $T_2\equiv f^k$
are the two transitions.}
\label{fig5}
\end{center}
\end{figure}

Let $f$ have a
$cs$-blender $\Gamma_1$ and a $cu$-blender $\Gamma_2$ with
superposition domains $\mathcal{B}_1$ and $\mathcal{B}_2$
respectively. Assume there are $P\in \Gamma_2$ and $Q\in
\Gamma_1$ so that
$$
\text{$W^{ss}(P)$ and $W^{uu}(Q)$ meets $\mathcal{B}_1$ and
$\mathcal{B}_2$ as is outlined in Figure~\ref{fig4}.}
$$
This implies there are $X\in \mathcal{B}_2 \cap
W^{s}(\Gamma_2)$, $Y \in \mathcal{B}_1 \cap W^{u}(\Gamma_1)$
and integers $j,k\in \mathbb{N}$ so that
\begin{equation}
 \label{(RC1)}
 \text{$W^{ss}_{loc}(f^{-j}(X))$ crosses
 $\mathcal{B}_1$  and   $W^{uu}_{loc}(f^{k}(Y))$ crosses $\mathcal{B}_2$,
 as shown in Figure~\ref{fig4}}.
\end{equation}
Then~\eqref{(B)} and \eqref{(B2)} imply
that
$$
   W^s(\Gamma_1) \cap W^u(\Gamma_2)\not=\emptyset  \quad \text{and}
   \quad W^u(\Gamma_1) \cap W^s(\Gamma_2)\not=\emptyset \quad
   \text{$C^1$-robustly}.
$$
In Figure~\ref{fig5} the heterodimensional
cycle is re-interpreted via the projection mappings on $E^c$.

\subsection{Criterion for robust tangencies}

\begin{figure}
\footnotesize
\begin{center}
\begin{picture}(400,250)
\put(0,0){\includegraphics[scale=0.8]{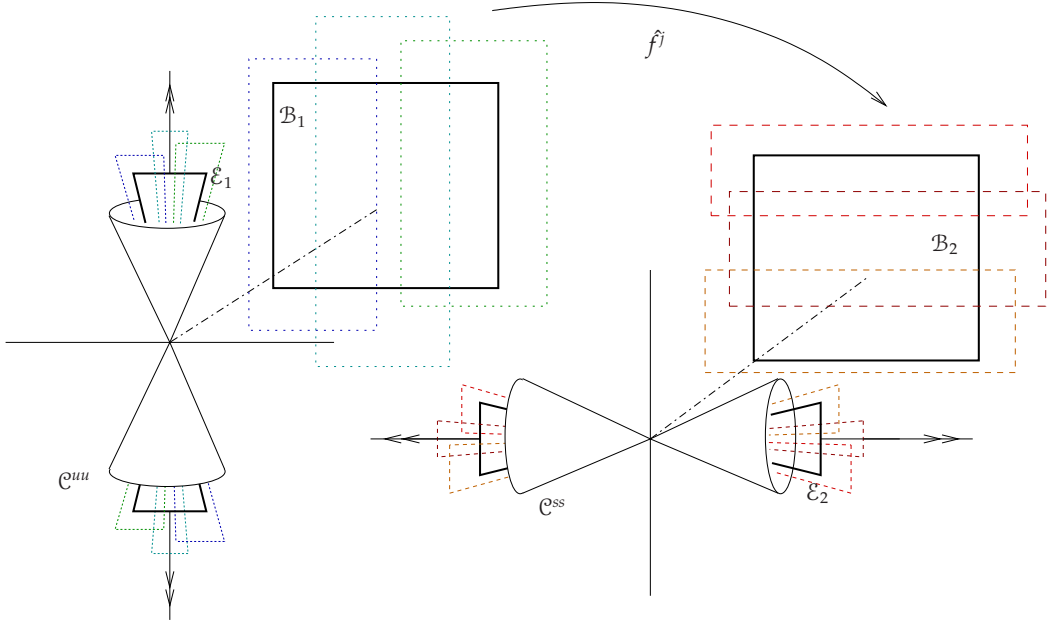}}
\put(20,50){$\mathcal{C}^{uu}$} \put(77,165){$\mathcal{E}_1$}
\put(103,188){$\mathcal{B}_1$} \put(347,140){$\mathcal{B}_2$}
\put(200,40){$\mathcal{C}^{ss}$} \put(300,46){$\mathcal{E}_2$}
\put(240,215){$\hat{f}^j$}
\end{picture}
\caption{The figure shows the covering properties of the blenders
of $\hat{f}$ in the product manifold and the tangent space as seen
in the central direction. The superposition domains are
$\hat{\mathcal{B}}_i= \mathcal{B}_i \times \mathcal{E}_i$.
Condition~\eqref{(T2)} implies that the covering property occurs
simultaneously for the regions $\mathcal{B}_i$ on the manifold,
and $\mathcal{E}_i$ which are inside the respective cones.}
\label{fig6}
\end{center}
\end{figure}

Let $f$ be as above, having a
$cs$-blender $\Gamma_1$ and a $cu$-blender $\Gamma_2$ with
superposition domains $\mathcal{B}_1$ and $\mathcal{B}_2$
respectively.  Observe that $f$ induces a map $\hat{f}$ on the set of $\ell$-dimensional subspaces of the tangent bundle given by 
$$
  \hat f(x,V)= (f(x),Df(x)V) \quad \text{where $x\in \mathcal{M}$ and $V \subset T_x \mathcal{M}$ with $\dim V=\ell$.}
$$
To create robust tangencies the key idea is to construct a robust
cycle as before, but for the map $\hat{f}$. To accomplish this we
will need a cone field  $\mathcal{C}^{uu}$ of dimension $\ell$ on
$\mathcal{B}_1$ which is $n$-th forward invariant and expanded by
$n$-th forward iteration: there are  $C>0$ and $0<\lambda<1$ so
that
 for every $n\in \mathbb{N}$ and $x\in
\mathcal{B} \cap f^{-n}(\mathcal{B})$,
\begin{equation*}
Df^n(x)\mathcal{C}^{uu}_x \subset
\mathrm{int}(\mathcal{C}^{uu}_{f^n(x)}) \quad \text{and} \quad \|Df^n(x)v\|\geq C \lambda^{-n}\|v\| \ \ \text{for all $v\in
\mathcal{C}_x^{uu}$}.
\end{equation*}
Also, assume there is a cone field $\mathcal{C}^{ss}$ on $\mathcal{B}_2$ of dimension $\ell$ which is $n$-th backward invariant and expanded
by $n$-th backward iteration.
To construct the cycle for $\hat{f}$, we require a similar criterion as~\eqref{(RC1)}, that is, two blenders for $\hat{f}$ and
a connection between them along
the strong stable manifold, see Figure~\ref{fig6} and compare with Figures~\ref{fig4} and \ref{fig5}.

\begin{enumerate}[label=(T\arabic*), ref=T\arabic*, itemsep=0.1cm]
\item \label{(T2)} There are a $cs$-blender $\hat\Gamma_1$ and
a $cu$-blender $\hat\Gamma_2$ for $\hat f$ with superposition
domains $\hat{\mathcal{B}}_1$ and $\hat{\mathcal{B}}_2$
respectively so that $\hat{\mathcal{B}}_1 \subset \mathcal{B}_1
\times \mathcal{C}^{uu}$ and $\hat{\mathcal{B}}_2 \subset
\mathcal{B}_2 \times \mathcal{C}^{ss}$, where by
$B\times\mathcal{C}$ we denote the set of points of the form
$(x,v)$ with $x\in \mathcal{B}$ and $v \in \mathcal{C}_x$.
\item \label{(T3)} There is
$X\in \Lambda^s(\hat{\mathcal{B}}_2)$ and $j\in \mathbb{N}$ so
that $W^{ss}_{loc}(\hat f^{-j}(X))$ crosses $\hat{\mathcal{B}}_1$.
\end{enumerate}

This cycle provides a (robust) tangency  as is outlined in the next argument.
Condition~\eqref{(T3)} implies that $X$ belongs to $W^s(\hat{\Gamma}_2)$ and $\hat{f}^{-j}(X)$ is in the superposition domain $\hat{\mathcal{B}}_1$ of the blender $\hat{\Gamma}_1$.
Using the property~\eqref{(B)} of the blender, there is a point $Y=(y,V)$ in the intersection between  $\Lambda^u(\hat{\mathcal{B}}_1)$
and $W^{ss}_{loc}(\hat f^{-j}(X))$. Then $y$ will be a point of tangency between $W^u(\Gamma_1)$ and $W^s(\Gamma_2)$
for the original map $f$, with the tangent vector space $V$. Indeed, by~\eqref{(T2)},
we have that the projection of $\Lambda^u(\hat{\mathcal{B}}_1)$
on the manifold is contained in $\Lambda^u(\mathcal{B}_1)$,
and analogously for $\Lambda^s(\hat{\mathcal{B}}_2)$. 
Also, since $f^j(y)$ is in the projection of $ W^{ss}_{loc}(X) \subset \Lambda^s(\hat{\mathcal{B}}_2)$ then
$$
\text{$y\in \Lambda^u(\mathcal{B}_1)\subset W^{u}(\Gamma_1)$, \
$f^j(y)\in \Lambda^s(\mathcal{B}_2) \subset W^s(\Gamma_2)$ and
thus   $y \in W^u(\Gamma_1)\cap W^s(\Gamma_2)$.}
$$
Now the backward iterates of $Y$ by $\hat{f}$ go to $\hat{\Gamma}_1$, and
in fact are in $\hat{\mathcal{B}}_1$ for some sequence of iterates $\hat{f}^{-m_i}$, $m_i\to\infty$. Since $\hat{\mathcal{B}}_1 \subset \mathcal{B}_1 \times
\mathcal{C}^{uu}$ and $\hat{f}^{-m_i}(Y)=(f^{-m_i}(y), Df^{-m_i}(y)V)$, then
$$f^{-m_i}(y)\in \mathcal{B}_1\cap f^{-m_i}(\mathcal{B}_1)\quad \text{and} \quad  Df^{-m_i}(y)V\in \mathcal{C}^{uu}_{f^{-m_i}(y)}.
$$
Consequently, for a vector $v\in V$, $Df^{-m_i}(y)v$ belongs to
the expanding cone field $\mathcal{C}^{uu}$ and thus
$\|Df^{-m_i}(y)v\|\leq C \lambda^{m_i}\|v\|$, which implies
$\|Df^{-n}(y)v\| \to 0$ as $n\to \infty$. Similarly also
$\|Df^{n}(y)v\| \to 0$ exponentially fast as $n\to\infty$.
Therefore $V$ is a tangent vector space and hence we obtain a
(robust) tangency of dimension $\dim V = \ell$.

\begin{rem} When
$\Gamma_1=\Gamma_2 \equiv \Gamma$, we conclude a (robust) homoclinic tangency of codimension $\ell$ associated with the
$double$-blender $\Gamma$.
\end{rem}

\section{Preliminaries:}
\label{s:preliminares} Let  $\mathscr{A}$ be  a finite set (with
at least two points), that we call an alphabet of symbols, and fix
$0<\nu<1$ and $0<\alpha\leq 1$. Consider the product space
$\Sigma\equiv\Sigma(\mathscr{A},\nu)\eqdef \mathscr{A}^\mathbb{Z}$
of the bi-sequences $\xi=(\xi_i)_{i\in\mathbb{Z}}$ of symbols in
$\mathscr{A}$ endowed with the metric
\begin{equation*}
\label{e:metrica} d_{\Sigma}(\xi,\zeta)\eqdef \nu^{\ell}, \ \
\ell=\min\{i\geq 0: \xi_i\not=\zeta_i\ \text{or} \
\xi_{-i}\not=\zeta_{-i} \}.
\end{equation*}
We remark that results in 
\S\ref{s:preliminares} and \S\ref{s:tangencies} are valid
relaxing the finiteness of the alphabet 
assuming that $\mathscr{A}$ is a compact metric space and $\Sigma$
is endowed with the metric given in~\cite[Example~2.1]{AV10}. Let
$M$ be a topological manifold (not necessarily compact and not
necessarily boundaryless).

\subsection{Symbolic skew-products} \label{s:symbolic}
Given a compact set
$K$ in $M$, we consider the pseudometric in the set $C^0(M)$ of
continuous functions of $M$ given by
\begin{equation}
\label{eq:metrica00}
  d_{C^0}(\phi,\psi)_K \eqdef \max_{x,y\in K} d(\phi(x),\psi(x))
  \quad \text{for any $\phi, \psi \in C^0(M)$.}
\end{equation}
Since $M$ is $\sigma$-compact, there is a sequence of relatively
compact subsets $K_n$ whose union is $M$ and then we can endow
$C^0(M)$ with the weak topology (also called compact-open
topology) induced by the family of
pseudometrics~\eqref{eq:metrica00}. That is,
$$
   d_{C^0}(\phi,\psi)=\sum_{n=1}^\infty 2^{-n} \,
   \frac{d_{C^0}(\phi,\psi)_{K_n}}{1+d_{C^0}(\phi,\psi)_{K_n}},
   \quad
   \text{for any $\phi,\psi \in C^0(M)$}.
$$


\subsubsection{The sets of  skew-products}

We consider skew-product homeomorphisms of the form
\begin{equation}
\label{e:sym-skew}
   \Phi: \Sigma \times M  \to \Sigma \times M, \qquad \Phi(\xi,x)=(\tau(\xi), \phi_\xi(x))
\end{equation}
where the base map $\tau:\Sigma\to\Sigma$ is the lateral shift map
and the fiber maps $\phi_\xi:M\to M$ are homeomorphisms of $M$. In
order to emphasize the role of the fiber maps we write
$\Phi=\tau\ltimes\phi_\xi$ and call it a \emph{symbolic
skew-product}. When no confusion arises we also write
$\mathcal{M}=\Sigma\times M$.

For every $n>0$ and $(\xi,x)\in \mathcal{M}$ set
\begin{equation}
\label{n.seq}
\begin{aligned}
   \phi^n_\xi(x) \eqdef \phi_{\tau^{n-1}(\xi)}\circ\cdots\circ\phi_{\xi}(x) 
    \quad \text{and} \quad
    \phi^{-n}_\xi(x) \eqdef \phi^{-1}_{\tau^{-n}(\xi)}\circ\cdots
    \circ\phi^{-1}_{\tau^{-1}(\xi)}(x)
\end{aligned}
\end{equation}
and hence
$$
\Phi^n(\xi,x)= (\tau^n(\xi),\phi^n_\xi(x)) \quad \text{for all
$n\in\mathbb{Z}$}.
$$

We introduce the set of symbolic skew-products with which we will
work:

\begin{defi}
Denote by
$\mathcal{S}(M)\equiv\mathcal{S}^\alpha_{\mathscr{A},\nu}(M)$ the
set of $\alpha$-H\"older 
symbolic skew-product homeomorphisms of $\mathcal{M}$. 
This is, the set of symbolic skew-products $\Phi=\tau\ltimes
\phi_\xi$ as in~\eqref{e:sym-skew} such that
\begin{enumerate}[leftmargin=0.6cm, label=\textbullet]
\item
$\phi_\xi$ are positive bi-Lipschitz homeomorphisms (uniform in
$\xi$): there are  positive constants $\gamma\equiv\gamma(\Phi)>0$
and $\hat\gamma\equiv \hat\gamma(\Phi)>0$ such that
\begin{equation}
\label{eq:lipschizt}
 \gamma\,d(x,y)< d(\phi_\xi(x),\phi_\xi(y)) <
\hat\gamma^{-1} \, d(x,y), \quad \mbox{ for all $x, y \in M$ and
$\xi \in \Sigma$,}
\end{equation}
\item $\phi_\xi$ depend $\alpha$-H\"older with respect to
$\xi$: there is a non-negative constant  $C_0 \equiv C_0(\Phi)\geq
0$ such that
\begin{equation} \label{eq:Holder} d_{C^0}(\phi^{\pm 1}_\xi,\phi^{\pm
1}_\zeta) \leq C_0 \, d_{\Sigma}(\xi,\zeta)^{\alpha} \quad
\mbox{for all $ \xi, \zeta \in \Sigma$ with $\xi_0=\zeta_0$.}
\end{equation}
\end{enumerate}
\end{defi}

We define in $\mathcal{S}(M)$ the metric
\begin{equation}
\label{eq:metrica0}
  d_{\mathcal{S}}(\Phi,\Psi) \eqdef
  d_0(\Phi,\Psi) + \mathrm{Lip}_0(\Phi,\Psi)+
  \mathrm{Hol}_0(\Phi,\Psi)
\end{equation}
where the symbolic skew-products $\Phi=\tau\ltimes\phi_\xi$ and
$\Psi=\tau\ltimes\psi_\xi$ belong to $\mathcal{S}(M)$ and
\begin{gather*}
\mathrm{Lip}_0(\Phi,\Psi)\eqdef \max_{\xi \in \Sigma}
  \big|\mathrm{Lip}(\phi_\xi^{\pm 1})-\mathrm{Lip}(\psi_\xi^{\pm
  1})\big|, \\
   d_0(\Phi, \Psi)\eqdef \max_{\xi \in \Sigma} \,
 d_{C^0}(\phi^{\pm 1}_\xi,\psi^{\pm 1}_\xi)
\quad  \text{and} \quad \mathrm{Hol}_0(\Phi,\Psi)\eqdef
\big|C_0(\Phi)-C_0(\Psi)\big|
\end{gather*}
with
$$
\mathrm{Lip}(\phi)=\sup_{ x \not = y}
\frac{d(\phi(x),\phi(y))}{d(x,y)} \quad \text{being $\phi$ a
bi-Lipschitz homeomorphism of $M$.}
$$
An important class of $\alpha$-H\"older continuous symbolic
skew-products is the following:

\begin{defi}
\label{def:PHS0} A symbolic skew-product $\Phi=\tau\ltimes
\phi_\xi \in \mathcal{S}(M)$ is \emph{partially hyperbolic}, if
$$\nu^\alpha < \gamma < 1 < \hat{\gamma}^{-1}< \nu^\alpha$$
where $\gamma$ and $\hat\gamma$ are given in~\eqref{eq:lipschizt}.
 We
denote by
$\mathcal{PHS}(M)\equiv\mathcal{PHS}^\alpha_{\mathscr{A},\nu}(M)$
the set of partially hyperbolic symbolic skew-products. Notice
that  $\mathcal{PHS}(M)$ is an open subset of $\mathcal{S}(M)$.

\end{defi}


\subsubsection{Stable and unstable sets for skew-products}

We define the \emph{local stable} and \emph{unstable} set of the
lateral shift map $\tau:\Sigma \to \Sigma$ at $\xi\in \Sigma$
respectively as
\begin{align*}
W^s_{loc}(\xi) \equiv W^s_{loc}(\xi;\tau) \eqdef \{\zeta \in
\Sigma : \zeta_i=\xi_i, \ i\geq 0 \}, \\
\quad W^u_{loc}(\xi) \equiv W^u_{loc}(\xi;\tau) \eqdef \{\zeta \in
\Sigma: \zeta_i=\xi_i, \ i<0\}.
\end{align*}
The \emph{(global) stable set} of the skew-product map
$\Phi:\mathcal{M} \to \mathcal{M}$ at $P\in \mathcal{M}$ is
defined as
$$
W^{s}(P)\equiv W^{s}(P;\Phi)\eqdef \{Q\in \mathcal{M}:
\displaystyle\lim_{n\to\infty} d(\Phi^{n}(Q),\Phi^{n}(P)= 0 \}.
$$
We define the \emph{(global) stable set} of  a compact
$\Phi$-invariant set, i.e.~so that $\Phi(\Gamma)=\Gamma$,  by
$$
W^s(\Gamma) \equiv W^s(\Gamma;\Phi)\eqdef \{P\in \mathcal{M}:
    \lim_{n\to\infty} d(\Phi^{n}(P),\Gamma)= 0\}
$$
or equivalently as the set of the points of $\mathcal{M}$ so that
its $\omega$-limit is contained in $\Gamma$. The set $\Gamma$ is
called \emph{isolated} (or \emph{maximal invariant set}) if there
is a compact neighborhood $\mathcal{U}$ of $\Gamma$, called the
\emph{isolating neighborhood} for $\Gamma$, such that every
invariant subset of $\mathcal{U}$ lies in $\Gamma$.
In such a case, we introduce the \emph{local stable set} of
$\Gamma$ as the forward invariant set of $\Phi$ in the isolating neighborhood
$\mathcal{U}$, that is,
\begin{equation*}
W^s_{loc}(\Gamma) \equiv W^{s}_{loc}(\Gamma;\Phi) \eqdef
  \{P\in\mathcal{M}: \ \Phi^n(P)\in \mathcal{U}  \  \text{for $n\geq 0$}\} = \bigcap_{n\geq 0} \Phi^n(\mathcal{U}).
\end{equation*}
Similarly  $W^u_{loc}(\Gamma) \equiv W^u_{loc}(\Gamma;\Phi)$ and
$W^u(\Gamma) \equiv W^u(\Gamma;\Phi)$ are, respectively, the
\emph{local unstable set} and the \emph{global unstable set} of
$\Gamma$. We have that
$$
   W^s(\Gamma)=\bigcup_{n\geq 0} \Phi^{-n}(W^s_{loc}(\Gamma)) \quad \text{and} \quad W^u(\Gamma)= \bigcup_{n\geq0} \Phi^{n}(W^u_{loc}(\Gamma)).
$$

Finally, given an \emph{$\mathcal{S}$-perturbation} of $\Phi$,
that is a symbolic skew-product $\Psi$ close to $\Phi$ in the
metric given in~\eqref{eq:metrica0}, we denote by $\Gamma_\Psi$
the maximal invariant set in $U$ of $\Psi$. Although isolated sets
vary, a priori, just upper semicontinuously by an abuse of
terminology, we call $\Gamma_\Psi$ the \emph{continuation} of
$\Gamma$ for $\Psi$.

\subsubsection{Strong laminations for partially hyperbolic skew-products}

Under the global assumption of domination introduced in
Definition~\ref{def:PHS0}, the usual graph transform argument
yields a local strong stable $\mathcal{W}^{ss}$ and unstable
$\mathcal{W}^{uu}$ partition for partially hyperbolic symbolic
skew-products:


\begin{prop}[\cite{AV10,ASV12}]
\label{p:inv-s} For every $\Phi\in\mathcal{PHS}(M)$ there exist
unique partitions
\begin{align*}
\mathcal{W}^{ss}=\{W^{ss}_{loc}(\xi,x): \, (\xi,x) \in \mathcal{M}
\} \quad \text{and} \quad \mathcal{W}^{uu}=\{W^{uu}_{loc}(\xi,x):
\, (\xi,x) \in \mathcal{M} \}
\end{align*}
of $\mathcal{M}=\Sigma \times M$ such that it holds
\begin{enumerate}[itemsep=0.1cm]
\item \label{item-kk1} every leaf $W^{ss}_{loc}(\xi,x)$
 is the graph of an $\alpha$-H\"older function
 $\gamma^{s}_{\xi,x} : W^s_{loc}(\xi) \to M$
 with $\alpha$-H\"older constant less or equal than
 $C_0 \cdot (1-\gamma^{-1} \nu^\alpha )^{-1}\geq 0$; 
\item $W^{ss}_{loc}(\xi,x)$
varies continuously with respect to $(\xi,x)$, i.e., the map $
(\xi,\xi',x) \mapsto \gamma^{s}_{\xi,x}(\xi')$ is continuous where
$(\xi, \xi')$ varies in the space of pairs of points in the same
local stable set for $\tau$. \\ Furthermore, it depends
continuously on $\Phi$;
\item \label{item-kk2} $\Phi(W^{ss}_{loc}(\xi,x)) \subset
W^{ss}_{loc}(\Phi(\xi,x))$ for all $(\xi,x) \in \mathcal{M}$;
\item $W^{ss}_{loc}(\xi,x) \subset
W^{s}(\xi,x)$ for all $(\xi,x) \in \mathcal{M}$.
\end{enumerate}
The partition $\mathcal{W}^{uu}$ verifies analogous properties.
\end{prop}

Each leaf of the partition $\mathcal{W}^{ss}$  is called the
\emph{local strong stable set}. We define the \emph{(global)
strong stable set} of $\Phi$ at $P$ as
\begin{equation*}
 W^{ss}(P) \equiv W^{ss} (P;\Phi ) \eqdef
\displaystyle\bigcup_{n\geq 0} \Phi^{-n} (
W^{ss}_{loc}(\Phi^{n}(P))) \subset W^s(P).
\end{equation*}
Let $\Gamma$ be an isolated set. We define the \emph{local strong
stable set} of $\Phi$ at $\Gamma$ as
$$
W^{ss}_{loc}(\Gamma) \equiv
W^{ss}_{loc}(\Gamma;\Phi)\eqdef\bigcup_{P\in\Gamma}W^{ss}_{loc}(P).
$$
In the same manner, the \emph{(global) strong stable set} of
$\Phi$ at $\Gamma$, $W^{ss}(\Gamma)$, is defined. Also the definitions of \emph{local/global
strong unstable sets of $\Phi$ at $\Gamma$} are analogous.

\subsection{Blenders in symbolic skew-products}
\label{s:symbolic-skew-products} Roughly speaking, a blender is a
basic (hyperbolic) set of a dynamical system, which provides that
a non-transversal intersection between stable/unstable manifolds
becomes a robust intersection. In this section, we will first
introduce the notion of hyperbolic set for symbolic skew-products
homeomorphisms. After that we give the formal definition of
blenders and finally we provide a criterion to obtain these local
tools.

\subsubsection{Hyperbolic sets}
\label{ss:hyperbolic} The following definitions and results come
from the literature on hyperbolic
homeomorphisms~\cite{Akin93,Ombach96}.

Fix $\varepsilon>0$ small enough. We introduce the \emph{local
stable set (of size $\varepsilon$)}  of $\Phi$ at $P=(\xi,x)$ as
$$
   W^{s}_{\varepsilon}(P)\equiv W^{s}_{\varepsilon}(P;  \Phi) \eqdef \{Q\in \mathcal{M}: \
   d(\Phi^n(Q),\Phi^n(P))\leq \varepsilon, \ \ n\geq 0
   \} \subset W^s_{loc}(\xi)\times M.
$$
The local unstable set (of size $\varepsilon$), denoted by
$W^u_\varepsilon(P)$, is defined analogously.

\begin{defi}
\label{def:hyp} A compact invariant set $\Gamma\subset
\mathcal{M}$ is \emph{hyperbolic} (for $\Phi$) if there exist
constants $\varepsilon>0$, $K>0$, $0<\theta<1$ such that
\begin{align*}
d(\Phi^n(P),\Phi^n(Q))\leq K \theta^n \ \ \text{for all
$P\in \Gamma$,  $Q\in W^s_{\varepsilon}(P)$ and $n\geq 0$;} \\
d(\Phi^{-n}(P),\Phi^{-n}(Q))\leq K \theta^n \ \ \text{for all
$P\in \Gamma$,  $Q\in W^u_{\varepsilon}(P)$ and $n\geq 0$;}
\end{align*}
and there exists $\delta>0$ such that
$$
\# \, W^s_{\varepsilon}(P) \cap W^u_{\varepsilon}(Q) =1 \ \
\text{for all $P, \, Q \in \Gamma$ with $d(P,Q) \leq \delta$}.
$$
\end{defi}

An isolated set is hyperbolic if and only if it is
\emph{expansive} and has the \emph{shadowing property} (see
definitions in~\cite{Akin93}). Every isolated hyperbolic set
$\Gamma$ for $\Phi$ is topologically stable; i.e., there is an
isolating neighborhood $U$ of $\Gamma$ such that for any
homeomorphism $\Psi$ which is $C^0$ near $\Phi$, the restriction
of $\Psi$ to the maximal invariant set in $U$ (that is, to the
continuation $\Gamma_\Psi$ of $\Gamma$), is semiconjugate to the
restriction of $\Phi$ to $\Gamma$.

We define
$$
[P,Q]=W^s_{\varepsilon}(P)\cap W^u_{\varepsilon}(Q) $$ for every
$P$ and $Q$ in $\Gamma$ with $d(P,Q)\leq \delta$ where
$\varepsilon>0$ and $\delta>0$ follow from the hyperbolicity of
$\Gamma$. The hyperbolic set $\Gamma$ has a \emph{local product
structure} if $[P,Q]\in \Gamma$ or equivalently if $\Gamma$ is an
isolated set of $\Phi$. In such a case, the map $$[ \cdot, \cdot]:
\{(P,Q)\in \Gamma\times \Gamma: d(P,Q)\leq \delta\} \to \Gamma$$
is continuous~(see~\cite{Saki95}). Hence, the map $[\cdot,\cdot]:
M^{cs}_\varepsilon(P)\times M^{cu}_\varepsilon(P) \to M$ is a
homeomorphism onto its image for all $P=(\xi,x)\in \Gamma$ where
$$
M^{cs}_{\varepsilon}(P)=W^{s}_\varepsilon(P)\cap (\{\xi\}\times M)
\quad \text{and} \quad
M^{cu}_{\varepsilon}(P)=W^{u}_\varepsilon(P)\cap (\{\xi\}\times
M).$$ In the sequel we will assume that the topological dimension
(in the sense of the Lebesgue covering dimension) of
$M^{cs}_{\varepsilon}(P)$ and $M^{cu}_{\varepsilon}(P)$ depend
continuously with respect to $P=(\xi,x)\in\Gamma$.

We will now introduce the notion of index of an isolated
transitive hyperbolic set $\Gamma$ in our context.
From the above assumption, the dimensions of
$M^{cs}_{\varepsilon}(P)$ and $M^{cu}_\varepsilon(P)$ are locally
constant.
Hence, being $\Gamma$  transitive, the dimensions 
remain constant for any $P\in\Gamma$. Thus, we may define the
\emph{$cs$-index} (resp.~\emph{$cu$-index}) of $\Gamma$, denoted
by $\mathrm{ind}^{cs}(\Gamma)$ (resp.~$\mathrm{ind}^{cu}(\Gamma)$)
as this dimension. Notice that $\dim M=
\mathrm{ind}^{cs}(\Gamma)+\mathrm{ind}^{cu}(\Gamma)$ and from the
topological stability, the $cs$-index remains constant under small
$\mathcal{S}$-perturbations of $\Phi$.

\subsubsection{Blenders}
In order to introduce the notion of a blender we need first to
define families of $s$-discs and $u$-discs, which provide the
superposition regions. To do this, we will consider a \emph{basic
open} set $\mathcal{B}$ of $\mathcal{M}$, i.e., a set of the form
$\mathsf{V}\times B$ where $\mathsf{V}$ is an open set of $\Sigma$
and $B$ is an open set of $M$.
\begin{defi}[$s$-discs]
\label{discs} A set $\mathcal{D}^s \subset \mathcal{M}$ is called
a \emph{$s$-disc} in $\mathcal{B}$ if there is  $\xi \in
\mathsf{V}$ such that $\mathcal{D}^s$ is a graph of an
$\alpha$-H\"older function from $W^s_{loc}(\xi)\cap \mathsf{V}$ to
$B$.
\end{defi}

We say that two $s$-discs, $\mathcal{D}^s_1, \mathcal{D}^s_2
\subset W^s_{loc}(\xi)\times M$ are close if they are the graphs
of close $\alpha$-H\"older functions. This proximity between discs
allows us to introduce the following:

\begin{defi}[open set of $s$-discs]
We say that a collection of discs $\mathscr{D}^s$ is an \emph{open
set of $s$-discs in $\mathcal{B}$} if given $\mathcal{D}^s_0\in
\mathscr{D}^s$, every $s$-disc $\mathcal{D}^s$ close enough to
$\mathcal{D}^s_0$ is a $s$-disc contained in $\mathcal{B}$ and
belongs~to~$\mathscr{D}^s$.
\end{defi}

Example of $s$-discs are the \emph{almost horizontal discs}
defined as follows: given $\delta>0$ and a point $(\xi, x)\in
\mathcal{B}$, we say that a set $\mathcal{D}^s\equiv
\mathcal{D}^s(\xi, x) \subset \mathcal{M}$ is a
\emph{$\delta$-horizontal disc} in $\mathcal{B}$ if
\begin{enumerate}
 \item[-] $\mathcal{D}^s$ is a graph of a $(\alpha,C)$-H\"older function
 $g : W^s_{loc}(\xi)\cap \mathsf{V} \to B$,
 \item[-] $d(g(\zeta),x)<\delta$ for all
 $\zeta\in W^s_{loc}(\xi)\cap \mathsf{V}$,
 \item[-] $C\nu^\alpha <\delta$.
\end{enumerate}
The set of all $\delta$-horizontal discs in $\mathcal{B}$ is an
open set of $s$-discs in $\mathcal{B}$. Similarly we define
\emph{$u$-discs} in $\mathcal{B}$, \emph{open set of $u$-discs} in
$\mathcal{B}$ and we have that the set of \emph{almost vertical
discs} is an example of an open set of $u$-discs.

Following~\cite{NP12,BKR14}, we introduce symbolic $cs$, $cu$ and
$double$-blenders. 
\vspace{0.2cm}

\begin{defi}[blenders] \label{d:symbolic-blender}
Let $\Phi \in\mathcal{S}(M)$ be a symbolic skew-product. A
transitive hyperbolic maximal invariant set $\Gamma$ in a
relatively compact open set $\mathcal{U}\subset \mathcal{M}$ of
$\Phi$ is called
\begin{enumerate}
 \item \label{eq:cs-blender} $cs$-\emph{blender} if \,$\mathrm{ind}^{cs}(\Gamma)>0$ and there
 exist a basic open set $\mathcal{B}\subset \mathcal{U}$ and
 an open set $\mathscr{D}^s$ of $s$-discs in $\mathcal{B}$
 such that for every small enough
 $\mathcal{S}$-perturbation $\Psi$ of $\Phi$,
 \begin{equation*}
 W^{u}_{loc}(\Gamma_\Psi)\cap \mathcal{D}^s\neq\emptyset \quad \text{for all $\mathcal{D}^s \in \mathscr{D}^s$.}
 \end{equation*}
 \item  \label{eq:cu-blender} $cu$-\emph{blender} if \,$\mathrm{ind}^{cu}(\Gamma)>0$ and there
 exist a basic open set $\mathcal{B}\subset \mathcal{U}$ and
 an open set $\mathscr{D}^u$ of $u$-discs in $\mathcal{B}$
 such that for every small enough
 $\mathcal{S}$-perturbation $\Psi$ of $\Phi$,
 \begin{equation*}
 W^{s}_{loc}(\Gamma_\Psi)\cap \mathcal{D}^u\neq\emptyset \quad
 \text{for all $\mathcal{D}^u\in \mathscr{D}^u$}.
 \end{equation*}
 \item \emph{double-blender} if both~\eqref{eq:cs-blender} and~\eqref{eq:cu-blender} hold (not necessarily for the same $\mathcal{B}$).
\end{enumerate}
The open set $\mathcal{B}$ is called a \emph{superposition domain}
and the open sets of discs $\mathscr{D}^s$ and $\mathscr{D}^u$ are
called the \emph{superposition regions} of the blender. Finally,
the $cs$-blender (resp.~$cu$-blender) with $cs$-index
(resp.~$cu$-index) is equal to $\dim M$ is called a
$contracting$-blender (resp.~$expanding$-blender).
\end{defi}


The above condition in the definition of $cs$-blender
(resp.~$cu$-blender) about the positivity of the  $cs$-index
(resp.~$cu$-index) of $\Gamma$ is imposed in order to avoid
transversal intersections between $s$-discs (resp.~$u$-discs) and
local unstable (resp.~stable) sets of $\Phi$ at points of
$\Gamma$. We understand \emph{transversality} in the sense
of~\cite[Definition~1.9]{DRS98}. In fact, it is more convenient to
introduce the notion of the unstable intersection property in
Euclidian spaces: two compacta $X$ and $Y$ have \emph{unstable
intersection in $\mathbb{R}^c$} if every pair of continuous maps
$f:X\to \mathbb{R}^c$ and $g:Y\to \mathbb{R}^c$ can be
approximated arbitrarily closely by continuous maps $f':X\to
\mathbb{R}^c$ and $g':Y\to \mathbb{R}^c$ with $f'(X)\cap
g'(Y)=\emptyset$. If one of the compacta $X$ and $Y$ is
$0$-dimensional, then $X$ and $Y$ have unstable intersection if
and only if $\mathrm{dim}(X\times Y)<c$~\cite{D00}. In our case,
we assume that $\mathrm{ind}^{cs}(\Gamma)>0$ and consider an
$s$-disc $\mathcal{D}^s \subset W^s_{loc}(\xi)\times M$ and a
local unstable set $W^{u}_\varepsilon(P) \subset
W^u_{loc}(\xi)\times M$ where $P=(\xi,x)\in \Gamma$. Observe that
$\mathcal{D}^s\cap W^{u}_\varepsilon(P) \subset \{\xi\}\times M$.
Set $X$ and $Y$ as $\mathcal{D}^s \cap (\{\xi\}\times M)$ and
$M^{cu}_\varepsilon(P)=W^{u}_\varepsilon(P)\cap (\{\xi\}\times M)$
respectively and take $c=\dim M$. Since $X$ is a singleton, it is
$0$-dimensional and thus $$ \dim (X\times Y) \leq \dim  Y = \dim
M^{cu}_{\varepsilon}(P)=c-\mathrm{ind}^{cs}(\Gamma)<c. $$ This
implies that $X$ and $Y$ have unstable intersection in
$\mathbb{R}^c$ and thus $\mathcal{D}^s$ and $W^{u}_\varepsilon(P)$
have no transversal intersection. Consequently, from the
definition of a $cs$-blender (resp.~$cu$-blender) follows that the
dimension of $W^u_{loc}(\Gamma)$ (resp.~$W^s_{loc}(\Gamma)$) is
robustly larger than the $cu$-index (resp.~$cs$-index) of
$\Gamma$. Namely, the dimension of the projection of this local
unstable (resp. stable) set on $M$ is equal to $\dim M$. Notice
that $double$-blenders have simultaneously these large dimensional
projections.

\subsubsection{Covering property as a criterion to yield blenders}
\label{s:blender} In this section we give a criterion that allows
us to guarantee that a symbolic skew-product has a blender.

Given $i\in \mathscr{A}$  we define the \emph{$i$-horizontal
cylinder} and \emph{$i$-vertical cylinder}, respectively, by
\begin{equation*}
\mathsf{H}_{i} \eqdef \{ \xi \in \Sigma: \, \xi_{0}=i \} \quad
\text{and} \quad \mathsf{V}_{i} \eqdef \{ \xi \in \Sigma:
\,\xi_{-1}=i \}.
\end{equation*}
Let $\mathcal{B}$ be a subset of $\mathcal{M}$ and consider a
symbolic skew-product $\Phi:\mathcal{M}\to\mathcal{M}$. Set
\begin{align*}
   \Lambda^u(\mathcal{B};\Phi) &\eqdef \bigcap_{n\geq 0} \Phi^n(\mathcal{B})=
   \{P\in \mathcal{M} :  \Phi^{-n}(P)\in \mathcal{B} \ \text{for all $n\geq
   0$}\}.
\end{align*}
\begin{thm}[criterion for blenders]
\label{cs-thm} Let $\Phi \in \mathcal{PHS}(M)$ be a partially
hyperbolic symbolic skew-product. Assume that the following
$cs$-covering property holds:

%

There are a finite set $S\subset \mathscr{A}$ and open sets $B
\subset M$ and $B_i \subset M$ for all $i\in S$ such that
\begin{equation}
\label{cs-cover} \mathsf{V}_{i}   \times \overline{B_i} \subset
\Phi(\mathsf{H}_{i}\times B) \quad \text{for all $i\in S$,} \qquad
\overline{B} \subset \bigcup_{i\in S} B_i
\end{equation}
and $C\eqdef C_0 \cdot (1-\gamma^{-1}\nu^\alpha)^{-1} < L$ where
$L$ is the Lebesgue number of~\eqref{cs-cover}. Then for any
$0<\delta<\gamma L/2$ and for every small enough
$\mathcal{S}$-perturbation $\Psi$ of $\Phi$,
\begin{align*}
    \Lambda^u(\mathcal{B};\Psi)  \cap \mathcal{D}^s &\not=\emptyset,
\quad \text{for all $\delta$-horizontal discs $\mathcal{D}^s$ in
$\mathcal{B}\eqdef \Sigma_S^+\times B$ where
$\Sigma_S^+\eqdef\bigcup_{i\in S} \mathsf{H}_{i}$.}
\end{align*}
In addition, assume that $\mathcal{B}$ is contained in a
relatively compact open set $\mathcal{U}$ which is the isolating
neighborhood of a transitive hyperbolic set $\Gamma$ of $\Phi$
with $\mathrm{ind}^{cs}(\Gamma)>0$. Then $\Gamma$ is a
$cs$-blender of $\Phi$ whose superposition region contains the
open set of almost horizontal discs in $\mathcal{B}$.
\end{thm}

\begin{rem}
\label{rem:local-strong-disc} Observe that $C\geq 0$ is a H\"older
constant of the local strong stable partition of $\Phi$. Thus, if
$C\nu^\alpha<\gamma L/2$ then the superposition region of the
$cs$-blender above contains the family of local strong stable sets
in $\mathcal{B}$.
\end{rem}

\begin{rem}
\label{rem:cu-blender} Analogously, we get a $cu$-blender $\Gamma$
whose superposition region contains the open set of almost
vertical discs in
$$
\mathcal{B}\eqdef \Sigma^{-}_S \times B \quad \text{where} \ \
\Sigma_S^{-}\eqdef \bigcup_{i\in S}\mathsf{V}_{i}
$$
In this case $\mathcal{B}$ is contained in the isolating
neighborhood of $\Gamma$, the $cu$-index of $\Gamma$ is positive
and we have the following $cu$-covering property:
\begin{equation}
\label{cu-cover} \mathsf{H}_{i}   \times \overline{B_i} \subset
\Phi^{-1}(\mathsf{V}_{i}\times B) \quad \text{for all $i\in S$,}
\qquad \overline{B} \subset \bigcup_{i\in S} B_i
\end{equation}
with $\hat{C}\eqdef C_0 \cdot (1-\hat\gamma^{-1}\nu^\alpha)^{-1} <
L$. Moreover, if $\hat{C}\nu^\alpha < \hat\gamma L /2$ then the
superposition region contains the family of local strong unstable
sets in $\mathcal{B}$.
\end{rem}

\begin{rem}
A transitive hyperbolic set satisfying the covering
properties~\eqref{cs-cover} and \eqref{cu-cover} (not necessarily
for the same open set $B$) is a $double$-blender.
\end{rem}

A note on the proof of Theorem~\ref{cs-thm}. First notice that the
second part of the theorem is a consequence of the first part.
Indeed, since  $\Gamma$ is the maximal invariant set in
$\mathcal{U}$ and $\mathcal{B} \subset \mathcal{U}$ then
$\Lambda^u(\mathcal{B};\Phi)\subset W^u_{loc}(\Gamma)$.
This inclusion and the first part of the theorem conclude that
$\Gamma$ is a $cs$-blender of $\Phi$ whose superposition region
contains the open set of almost horizontal discs in $\mathcal{B}$.
We give the details of the first part of the proof in Appendix~\ref{AppendixB}.

\subsubsection{On the definition and criterion for blender}

Blenders are actually a power tool in partially hyperbolic
dynamics when the superposition region contains the local strong
stable/unstable set in the superposition domain. For this reason,
without loss of generality, we will assume the following blender
properties:

\begin{scho} Consider $\Phi \in
\mathcal{PHS}(M)$ and let $\Gamma$ be a $cs$-blender with
superposition domain $\mathcal{B}= \mathsf{V}\times B$ and
superposition region $\mathscr{D}^s$. There is a
$\mathcal{S}$-neighborhood $\mathscr{U}$ of $\Phi$ such that for
any $\Psi\in \mathscr{U}$,
\begin{enumerate}[label=(B\arabic*),ref=B\arabic*]
\item \label{cond:B1}
\emph{the open set of discs $\mathscr{D}^s$ contains the family of
local strong stable sets of $\Psi$ in $\mathcal{B}$:
$$
\text{if \ $W^{ss}_{loc}(P;\Phi)\cap (\mathsf{V}\times M) \subset
\mathcal{B}$ \ \ then \ \ $W^{ss}_{loc}(P;\Phi) \cap
(\mathsf{V}\times M) \in \mathscr{D}^s$.}
$$}
\item \label{cond:B2}
\emph{if $W^{ss}_{loc}(P;\Phi)\cap (\mathsf{V}\times M) \subset
\mathcal{B}$ then
$$
         W^u_{loc}(\Gamma_\Psi;\Psi) \cap W^{ss}_{loc}(P';\Psi)
         \not=\emptyset \quad \text{for all $P'$ close enough to $P$.}
$$}
\item \label{cond:B3} for any $P \in \Lambda^u(\mathcal{B};\Phi)$, there exists $P_\Psi$
close to $P$ so that
$$
          W^{uu}_{loc}(P_\Psi;\Psi)\cap\mathcal{B}\subset
          \Lambda^u(\mathcal{B};\Psi)\subset  W^u_{loc}(\Gamma_\Psi;\Psi)
$$
where
$$
   \Lambda^u(\mathcal{B};\Psi)=\bigcap_{n\geq 0} \Psi^n(\mathcal{B}).
$$
\end{enumerate}
Similar conditions are also assumed for $cu$-blenders of partially
hyperbolic skew-products.
\end{scho}

We must to show that the properties~\eqref{cond:B2}
and~\eqref{cond:B3} follow from the definition of blender, i.e.,
Definition~\ref{d:symbolic-blender}.

\begin{proof}
 First of all, notice that the
assumption~\eqref{cond:B1} and Definition~\ref{d:symbolic-blender}
imply that
\begin{equation}
\label{eq:3} \text{if $W^{ss}_{loc}(P;\Phi)\cap (\mathsf{V}\times
M) \subset \mathcal{B}$  \ then \
         $W^u_{loc}(\Gamma;\Phi) \cap W^{ss}_{loc}(P;\Phi)
         \not=\emptyset$ \quad $\mathcal{S}$-robustly.
         }
\end{equation}
A priori, the neighborhood of the $\mathcal{S}$-perturbation of
$\Phi$ where~\eqref{eq:3} holds depends on the $s$-disc
$W^{ss}_{loc}(P;\Phi)$. However, this can be taken independent of
the disc assuming that the disc belongs to a superposition
subdomain $\mathcal{B}_0$ of $\mathcal{B}=\mathsf{V}\times B$.
That is if $ W^{ss}_{loc}(P;\Phi)\cap (V\times M) \subset
\mathcal{B}_0$ where $\mathcal{B}_0=V\times B_0$, $B_0$ is an open
set whose closure is contained in $B$. For this reason, without
loss of generality, we can assume that~\eqref{cond:B2} holds.

With respect to~\eqref{cond:B3}, first notice that $\Gamma$ is the
maximal invariant set in $\mathcal{U}$. Since $\mathcal{B} \subset
\mathcal{U}$ then clearly  $\Lambda^u(\mathcal{B};\Phi) \subset
W^u_{loc}(\Gamma;\Phi)$. Let $\mathcal{B}_1=\mathsf{V}\times B_1$
where $B_1$ an open set containing the closure of $B$. If $P\in
\Lambda^u(\mathcal{B};\Phi)$, by the in phase
result~\cite[Prop.~10]{Ombach96}, there is $Q\in \Gamma \cap
\Lambda^u(\mathcal{B}_1;\Phi)$ such that $P\in
W^u_{\varepsilon}(Q;\Phi)$.
Hence, for any $\mathcal{S}$-perturbation $\Psi$ of $\Phi$, there
exists a continuation point $P_\Psi$ of $P$ and $Q_\Psi$ of $Q$ so
that $P_\Psi\in W^u_{\varepsilon}(Q_\Psi;\Psi)$ and $Q_\Psi\in
\Gamma_\Psi\cap \Lambda^u(\mathcal{B}_1;\Psi)$, and so $P_\Psi\in
\Lambda^u(\mathcal{B}_1;\Psi)$. This implies that $
W^{uu}_{loc}(P_\Psi;\Psi)\cap\mathcal{B}_1\subset
\Lambda^u(\mathcal{B}_1;\Psi)$.

Taking as above a superposition subdomain $\mathcal{B}_0$ of
$\mathcal{B}$ we obtain that the $\mathcal{S}$-neighborhood of the
perturbation is independent of the point $P$. Therefore, without
loss of generality, we can also assume that~\eqref{cond:B3} holds.
\end{proof}

Blenders constructed from covering property have extra useful
properties that are not implied directly from the definition. We
collect here these properties which will be used later on.

\begin{scho}
Let  $\Gamma$ be a $cs$-blender with a superposition domain
$\mathcal{B}=\mathsf{V}\times B$ for $\Phi\in\mathcal{PHS}(M)$,
constructed from the criterion of Theorem~\ref{cs-thm}

\noindent There exists a $\mathcal{S}$-neighborhood $\mathscr{U}$
of $\Phi$ such that for any $\Psi\in \mathscr{U}$,
\begin{enumerate}[label=(B\arabic*),ref=B\arabic*,start=4]
\item \label{cond:B4} if $W^{ss}_{loc}(P;\Psi)\cap (\mathsf{V}\times M) \subset \mathcal{B}$
then
$$
         \Lambda^u(\mathcal{B};\Psi) \cap W^{ss}_{loc}(P';\Psi)
         \not=\emptyset \quad \text{for all $P'$ close enough to $P$;}
$$
\item \label{cond:B5} $B\subset \mathscr{P}(\Lambda^u(\mathcal{B};\Psi))$, where $\mathscr{P}:\mathcal{M}\to M$ is the standard
projection on the fiber-space.
\end{enumerate}
Similar conditions also hold for $cu$-blenders with respect to
$$
   \Lambda^s(\mathcal{B};\Psi)=\bigcap_{n\geq 0} \Psi^{-n}(\mathcal{B})=\{P\in \mathcal{M}: \Psi^{n}(P)\in \mathcal{B} \ \
   \text{for $n\geq 0$}\}.
$$

\end{scho}


%




\section{Criteria: robust cycles and robust tangencies}
\label{s:tangencies}

\subsection{Criterion to yield robust cycles}
\label{ss:construction-cycles} We explain how blenders can be used
to yield robust cyclic intersections between stable and unstable
sets.

\begin{defi}[cycles]
\label{def:symb-cycle} Let $\Phi \in \mathcal{S}(M)$ be a symbolic
skew-product with a pair of disjoint isolated invariant sets
$\Gamma^1$ and $\Gamma^2$. We say that $\Phi$ has a \emph{cycle}
associated with $\Gamma^{1}$ and $\Gamma^{2}$ if their stable and
unstable sets meet cyclically, that is, if
\begin{equation*}
\label{e:ciclo} W^s(\Gamma^1)\cap W^u(\Gamma^2)\ne\emptyset \quad
\mbox{and} \quad W^u(\Gamma^1)\cap W^s(\Gamma^2)\ne\emptyset.
\end{equation*}
Assuming that $\Gamma^1$ and $\Gamma^2$ are also transitive
hyperbolic sets, the cycle is said to be
\emph{$\mathcal{S}$-robust} if any small enough
$\mathcal{S}$-perturbation $\Psi$ of ${\Phi}$ has a cycle
associated with the continuations $\Gamma^1_\Psi$ and
$\Gamma^2_\Psi$ of $\Gamma^1$ and $\Gamma^2$ respectively. We
define the \emph{co-index} of the cycle as the integer
$|\mathrm{ind}^{cs}(\Gamma^1)-\mathrm{ind}^{cs}(\Gamma^2)|=c$. If
$c>0$ the cycle is called \emph{heterodimensional} and otherwise
\emph{equidimensional}.
\end{defi}

Firstly, we would like to obtain a robust intersection between the
unstable and the stable sets of two different isolated
hyperbolic sets. To this end, we use the following criterion:

\begin{thm}[criterion for robust cycles]
\label{thm-cycles1} Let $\Gamma^1$ and $\Gamma^2$ be isolated
transitive hyperbolic sets of  $\Phi\in\mathcal{PHS}(M)$. Assume
that $\Gamma^1$ is a $cs$-blender with a superposition domain
$\mathcal{B}=\mathsf{V}\times B$ such that
\begin{enumerate}[label=(RC1),ref=RC1]
\item \label{cond:C2} there is a point $P\in W^s_{loc}(\Gamma^2)$ such
that
$$
W^{ss}_{loc}(\Phi^{-n}(P))\cap (\mathsf{V}\times M)\subset
\mathcal{B} \quad \text{for some integer $n\geq 0$}.
$$
\end{enumerate}
Then the unstable and stable sets of $\Gamma^1$ and $\Gamma^2$
respectively  meet $\mathcal{S}$-robustly. That is, for any small
$\mathcal{S}$-perturbation $\Psi \in \mathcal{PHS}(M)$ of $\Phi$
it holds $ W^u(\Gamma^1_\Psi)\cap W^s(\Gamma^2_\Psi)\neq\emptyset.
$
\end{thm}
\begin{proof}
Since $\Gamma^1$ is a $cs$-blender, then~\eqref{cond:B1}
and~\eqref{cond:C2} imply that $W^{ss}_{loc}(\Phi^{-n}(P))\cap
(\mathsf{V} \times M)\in \mathscr{D}^s$.
As
$\mathscr{D}^s$ is an open set of $s$-discs, then for every small
enough $\mathcal{S}$-perturbation $\Psi$ of $\Phi$,
$W^{ss}_{loc}(\Psi^{-n}(P'))\cap (\mathsf{V} \times M)\in
\mathscr{D}^s $, where $P'$ is the continuation\footnote{The
existence of $P'$ follows from the in phase
result~\cite[Prop.~10]{Ombach96} using that
 $\Gamma^2$ is an isolated hyperbolic set.} of $P$ in the local stable set of $\Gamma^2_\Psi$.
By definition of the $cs$-blender, $W^{ss}_{loc}(\Psi^{-n}(P'))\cap
W^u_{loc}(\Gamma^1_\Psi)\neq\emptyset$ and then
$\mathcal{S}$-robustly the global unstable set of $\Gamma^1$ meets
the global stable set of $\Gamma^2$.
\end{proof}

 Similar conditions guarantee the robustness of the
 other intersection. Therefore it is possible to construct
 robust  cycles associated with two  $cs$ or
 $cu$-blenders, or a $cs$-blender and a $cu$-blender, or even between two
 $double$-blenders.

\begin{rem}[cycles between non-hyperbolic sets]
The usefulness of the $double$-blender in this setting is that it
permits to construct robust cyclic connections with \emph{any} invariant set
(of saddle type) having a non-empty continuation.
For example, the other set may be a non-hyperbolic fixed point,
and then we obtain the cyclic intersections between the
stable and unstable sets, if the \emph{transition condition}
\eqref{cond:C2} is satisfied.
\end{rem}

\subsection{Criterion to yield robust tangencies}
\label{s:criteria-tangencias} We explain how robust cycles can be
used to robust tangencies. In the sequel $M$ denotes a
differentiable manifold of dimension $c\geq 1$. 

\subsubsection{The set of smooth symbolic skew-products}
Since we will need to work with differentiable fiber maps, it will
be useful to extend the sets of H\"older skew-products to this
setting.

\begin{defi}
For an integer $r\geq 1$,
$\mathcal{S}^{r}(M)\equiv\mathcal{S}^{r+\alpha}_{\mathscr{A},\nu}(M)$
denotes
 the set of  skew-products $\Phi=\tau\ltimes\phi_\xi$ on $\mathcal{M}=\Sigma\times M$
 such that there are  $\gamma\equiv\gamma(\Phi)>0$,
 $\hat\gamma\equiv\hat\gamma(\Phi)>0$ and
 $C_r\equiv C_r(\Phi)\geq 0$  satisfying
\begin{enumerate}
\item[-] $d_{C^r}(\phi^{\pm 1}_\xi,\phi^{\pm 1}_{\zeta}) \leq
C_r \, d_{\Sigma}(\xi,\zeta)^\alpha$ \ \ for all $\xi,\zeta \in
\Sigma$ with $\xi_0=\zeta_0$, and
\item[-] $\phi_\xi$ are $C^r$-diffeomorphisms of $M$ with
$D^r\phi^{\pm 1}_\xi$ Lipschitz (with uniform Lipschitz constant)
$$
\gamma< m(D\phi_\xi(x)) < \|D\phi_\xi(x)\| <\hat{\gamma}^{-1} \ \
\text{for all $(\xi,x) \in \Sigma\times M$}.
$$
\end{enumerate}
For $r=0$, $\mathcal{S}^0(M)$, we understand as $\Phi \in
\mathcal{S}(M)$ with fiber maps $C^1$-diffeomorphisms. In
addition,
$$
   \mathcal{PHS}^{r}(M)\equiv\mathcal{PHS}^{r+\alpha}_{\mathscr{A},\nu}(M)
   \eqdef \mathcal{PHS}(M)\cap \mathcal{S}^{r}(M) \qquad \text{for $r\geq 0$}.
$$
That is, the domination conditions  $\nu^{\theta} < \gamma <1 <
\hat\gamma^{-1} <\nu^{-\theta}$ hold where $0<\alpha\leq
\theta\leq 1$ is the exponent of the H\"older continuity in the
$C^0$-metric. Finally, a partially hyperbolic skew-product is said
to be \emph{fiber bunched} if $\nu^\alpha<\gamma\hat{\gamma}$.
\end{defi}

We endow $\mathcal{S}^r(M)$ for $r\geq 0$ with the metric
$$
 d_{\mathcal{S}^r}(\Phi,\Psi) \eqdef
  d_r(\Phi,\Psi) + \mathrm{Lip}_r(\Phi,\Psi)+
  \mathrm{Hol}_r(\Phi,\Psi)
$$
where
\begin{gather*}
\mathrm{Lip}_r(\Phi,\Psi)\eqdef \max_{\xi \in \Sigma}
  \big|\mathrm{Lip}(D^r\phi_\xi^{\pm 1})-\mathrm{Lip}(D^r\psi_\xi^{\pm
  1})\big| \\
   d_r(\Phi, \Psi)\eqdef \max_{\xi \in \Sigma} \,
 d_{C^r}(\phi^{\pm 1}_\xi,\psi^{\pm 1}_\xi)
\quad  \text{and} \quad
\mathrm{Hol}_r(\Phi,\Psi)\eqdef \big|C_r(\Phi)-C_r(\Psi)\big|.
\end{gather*}
Hence $\mathcal{PHS}^r(M)$ is an open set of $\mathcal{S}^{r}(M)$
and
$ \mathcal{S}^{r+1}(M) \subset \mathcal{S}^{r}(M) \subset
\mathcal{S}(M)$ for
    $r\geq  0$.

\subsubsection{Tangencies in symbolic skew-products} To define the notion of tangency for symbolic
skew-products we first need to introduce the notion of a tangent
direction. 


\begin{defi}[tangent direction] Let
$\Phi=\tau \ltimes \phi_\xi\in \mathcal{S}^{r}(M)$ be a symbolic
skew-product with a pair of transitive hyperbolic sets $\Gamma^1$
and $\Gamma^2$  and suppose $(\xi,x)\in W^{u}(\Gamma^1) \cap
W^{s}(\Gamma^1)$.  A unitary vector $v\in T_x M$ is called
\emph{tangent direction} at $(\xi,x)$ if there are $C>0$ and
$0<\lambda<1$~such~that
$$
\|D\phi^n_\xi(x)v\| \leq C \lambda^{|n|} \quad \text{for all
$n\in\mathbb{Z}$.}
$$
The maximum number of independent tangent directions at $(\xi,x)$
is denoted by $d_T\equiv d_T(\xi,x)$.
\end{defi}

Now we are ready to give the definition of a tangency.

\begin{defi}[tangency]
We say that $\Phi\in \mathcal{S}^r(M)$ has a \emph{(bundle)
tangency} of dimension $\ell> 0$ between $W^{u}(\Gamma^1)$ and
$W^s(\Gamma^2)$ if there exists $(\xi,x)\in W^u(\Gamma^1)\cap
W^s(\Gamma^2)$ such that
$$
      \ell=d_T(\xi,x) \quad\text{and} \quad   \mathrm{ind}^{cu}(\Gamma^1)+\mathrm{ind}^{cs}(\Gamma^2)-\ell < c.
$$
If $\Gamma^1=\Gamma^2$, the tangency is called \emph{homoclinic},
and otherwise \emph{heteroclinic}.
The tangency (of dimension $\ell$) is said to be
\emph{$\mathcal{S}^{r}$-robust} if for any small enough
$\mathcal{S}^{r}$-perturbation $\Psi$ of $\Phi$ has a tangency (of
dimension $\ell$) between the unstable set $W^{u}(\Gamma^1_\Psi)$
and the stable set $W^s(\Gamma^2_\Psi)$.

The \emph{codimension} of the tangency is defined as
$c_T=c-[\mathrm{ind}^{cu}(\Gamma^1)+\mathrm{ind}^{cs}(\Gamma^2)-\ell]$.
\end{defi}

Assume $\Phi$ has a ($\mathcal{S}$-robust) equi/heterodimensional
cycle associated with $\Gamma^1$ and $\Gamma^2$ and
$\mathrm{ind}^{cs}(\Gamma^1)\leq \mathrm{ind}^{cs}(\Gamma^2)$. We
say that
$\Phi$ has a ($\mathcal{S}^r$-robust) \emph{equi/heterodimensional
tangency  on the cycle} if $\Phi$ has a ($\mathcal{S}^r$-robust)
tangency between $W^u(\Gamma^1)$ and $W^s(\Gamma^2)$.

In what follows, we want to construct symbolic skew-products with
tangencies, which will be created locally. For this reason, in
order to be clearer, we will work in local coordinates and may
assume that $M=\mathbb{R}^c$.

\subsubsection{Cone fields in symbolic skew-products}

Consider an integer $1\leq \ell\leq c$. A standard $\ell$-cone in
$\mathbb{R}^c$ is a set of the form $ \mathcal{C}=\{(v,w)\in
\mathbb{R}^c: v\in\mathbb{R}^\ell  \ \text{and} \ \|w\| \leq \rho
\|v\| \ \ \text{for some} \ \rho>0\}$. More generally, a
\emph{$\ell$-cone} is the image of a standard $\ell$-cone under an
invertible linear map.

\begin{defi}[unstable cone]
 Let $\Phi=\tau\ltimes\phi_\xi \in
\mathcal{S}^r(\mathbb{R}^c)$ and consider an open set
$\mathcal{B}$ of $\mathcal{M}=\Sigma \times \mathbb{R}^c$.
An $\ell$-cone $\mathcal{C}^{uu}$ in $\mathbb{R}^c$ is said to be
\emph{unstable} for $\Phi$ on
$\mathcal{B}$ if there is  $0<\lambda<1$ such that
$$
D\phi_\xi(x)\mathcal{C}^{uu} \subset
\mathrm{int}(\mathcal{C}^{uu})
\quad \text{and} \quad 
\text{$\|D\phi_\xi(x)v\|\geq \lambda^{-1} \|v\|$
 \ \ for all $v\in
\mathcal{C}^{uu}$, \  $(\xi,x)\in \mathcal{B}\cap
\Phi^{-1}(\mathcal{B})$.}
$$
\end{defi}

Observe that the unstable $\ell$-cone $\mathcal{C}^{uu}$ is
$\mathcal{S}^0$-robust: there is a $S^0$-neighborhood
$\mathscr{U}$ of $\Phi$ so that $\mathcal{C}^{uu}$ is a unstable
$\ell$-cone for $\Psi$ on $\mathcal{B}$ for any $\Psi
\in\mathscr{U}$. Similarly we define the \emph{stable $\ell$-cone}
$\mathcal{C}^{ss}$ for $\Phi$ on  $\mathcal{B}$ and
$\mathcal{C}^{ss}$ is $\mathcal{S}^0$-robust.

\begin{lem}
\label{lem-unstable-cone} Let $\mathcal{C}^{uu}$ be an unstable
$\ell$-cone for $\Phi$ on $\mathcal{B}$. Then for each
$n\in\mathbb{N}$ it holds that
$$\|D\phi_\xi^{-n}(x)v\|\leq \lambda^n \|v\| \quad
\text{for $v \in \mathcal{C}^{uu}$ and $(\xi,x)\in \mathcal{B}
\cap \Phi(\mathcal{B})\cap \dots \cap \Phi^{n}(\mathcal{B})$}.$$
\end{lem}
\begin{proof}
It suffices to note that  $
  \|v\|=\|D\phi^n_\xi(\phi^{-n}_\xi(x))D\phi^{-n}_\xi(x)v\|
  \geq \lambda^{-n} \|D\phi^{-n}_\xi(x)v\|$.
\end{proof}

An $\ell$-dimensional vector subspace of $\mathbb{R}^c$ is called
a $\ell$-plane. The \emph{Grassmannian manifold} $G(\ell,c)$ is
defined as the set of $\ell$-planes in $\mathbb{R}^c$. Given $E$,
$F$ in $G(\ell,c)$, we define the distance
$$
d_G(E,F)=\inf\{\|A-i_E\|: \ \text{$A:E \to \mathbb{R}^{c}$ is a
linear operator such that $A(E)\subset F$}\}
$$
where $i_E:E\to \mathbb{R}^c$ denotes the inclusion and
$\|\cdot\|$ is the Euclidian operator norm. It is not difficult to
prove that the  infimum is attained at the map $A=P_F|E$ (i.e.,
the orthogonal projection onto $F$ restricted to $E$). Trivially,
$$
\|P_F|E-i_E\|=\sup \big\{ \ \inf_{f\in F} \|f-e\| : \  e\in E, \
\|e\|=1\big\}
$$
which provides another characterization of the above metric in
$G(\ell,c)$. It follows that any $\ell$-cone $\mathcal{C}$ in
$\mathbb{R}^c$ induces an open set in $G(\ell,c)$, which we will
continue denoting by $\mathcal{C}$. Moreover, any unstable
$\ell$-cone $\mathcal{C}^{uu}$ for $\Phi$ on $\mathcal{B}$ is
\emph{strictly invariant} by the action of
$\mathbb{A}=\{D\phi_\xi(x): (\xi,x)\in
\mathcal{B}\cap\Phi^{-1}(\mathcal{B})\}$. That is,
$\mathbb{A}\mathcal{C}^{uu}$ is contained in the interior of
$\mathcal{C}^{uu}$. Here $ \mathbb{A}\mathcal{C}^{uu}$ denotes the
image of $\mathcal{C}^{uu}\subset G(\ell,c)$ under the action of
$\mathbb{A}$, i.e., the subset $\{AE : A\in \mathbb{A}, E \in
\mathcal{C}^{uu}\}$ of $G(k, c)$. Similar observations hold for
stable $\ell$-cones.


\subsubsection{Criterion to yield robust tangencies}
\label{ss:get-tangencias}

Given $\Phi=\tau\ltimes\phi_\xi\in \mathcal{S}^0(\mathbb{R}^c)$
and $1\leq \ell\leq c$, we define the \emph{induced $\ell$-th
Grassmannian skew-product} $\hat\Phi$ on
$$
\hat{\mathcal{M}}\eqdef\Sigma \times \hat{M} \quad \text{where} \
\ \ \text{$\hat{M}\eqdef\mathbb{R}^c\times G(\ell,c)$} $$ as
$$
\hat\Phi:\hat{\mathcal{M}} \to \hat{\mathcal{M}}, \qquad
\hat\Phi(\xi,(x,E))= (\tau(\xi), (\phi_\xi(x), D\phi_\xi(x)E)).$$
\begin{thm}[criterion for robust tangencies]
\label{thm:tangencias} Let $\Gamma^1$ and $\Gamma^2$ be a
$cs$-blender and a $cu$-blender of a skew-product $\Phi\in
\mathcal{PHS}^1(\mathbb{R}^c)$ with, respectively, superposition
domains
$$
\text{$\mathcal{B}_1=\mathsf{V}_1\times B_1$ \ and \
$\mathcal{B}_2=\mathsf{V}_2\times B_2$, \ and \ $cs$-indices $0<
i_1<c$ \ and \ $0<i_2< c$.}
$$
Consider
an integer $ \max\{0,i_2-i_1\}< \ell \leq \min\{c-i_1,i_2\}$.
Assume that
\begin{enumerate}[label=(T\arabic*),ref=T\arabic*]
\item \label{T1} the induced $\ell$-th Grassmannian skew-product
$\hat{\Phi}$ on $\hat{\mathcal{M}}$ belongs to
$\mathcal{PHS}(\hat{M})$;
\end{enumerate}
and there exist
\begin{enumerate}[label=(T\arabic*),ref=T\arabic*,start=2]
\item \label{T2} unstable and stable $\ell$-cones $\mathcal{C}^{uu}$
and $\mathcal{C}^{ss}$ for $\Phi$ on $\mathcal{B}_1$ and
$\mathcal{B}_2$ respectively;
\item \label{T3} a $cs$-blender $\hat\Gamma^1$ and
a $cu$-blender $\hat\Gamma^2$ of $\hat\Phi$ with superposition
domains $\hat{\mathcal{B}}_1$ and $\hat{\mathcal{B}}_2$
respectively, such that $\hat{\Gamma}^1$ satisfies~\eqref{cond:B4}
and
$$
  \hat{\mathcal{B}}_1=\mathsf{V}_1\times \hat{B}_1, \quad
  \hat{\mathcal{B}}_2=\mathsf{V}_2\times \hat{B}_2, \quad
  \text{with \ \ $\hat{B}_1 \subset B_1\times \mathcal{C}^{uu}$ and $\hat{B}_2 \subset
B_2\times \mathcal{C}^{ss}$;}
$$
\item \label{T4} a point
$\hat{P}\in \Lambda^s(\hat{\mathcal{B}}_2;{\hat\Phi}) \subset
W^s_{loc}(\hat\Gamma^2)$ and an integer $n\geq 0$ such that
$$
   W^{ss}_{loc}(\hat\Phi^{-n}(\hat{P})) \cap (\mathsf{V}_1\times \hat{M} )
   \subset \hat{\mathcal{B}}_1.
$$
\end{enumerate}
Then
\begin{enumerate}
\item \label{item1} $W^u(\hat\Gamma^1) \cap W^s(\hat\Gamma^2)\not
=\emptyset$  \ \ $\mathcal{S}$-robustly  (in
$\mathcal{PHS}(\hat{M})$),
\item \label{item2} $W^u(\Gamma^1) \cap W^s(\Gamma^2)\not
=\emptyset$  \ \ $\mathcal{S}$-robustly (in
$\mathcal{PHS}(\mathbb{R}^d)$), and
\item \label{item3} $\Phi$ has a $\mathcal{S}^{1}$-robust tangency of dimension $\ell$
between $W^u(\Gamma^1)$ and $W^s(\Gamma^2)$.
\end{enumerate}
\end{thm}

In order to prove the above theorem, we first need a lemma:

\begin{lem}
\label{lem-projection} Consider $1\leq \ell \leq c$,  $\Phi\in
\mathcal{S}^0(\mathbb{R}^c)$ and let $\hat\Phi$ be the induced
$\ell$-th Grassmannian symbolic skew-product on
$\hat{\mathcal{M}}=\Sigma\times \hat{M}$ with $\hat M
=\mathbb{R}^c\times G(\ell,c)$.
\begin{enumerate}
\item If $\Phi\in \mathcal{S}^1(\mathbb{R}^c)$ then $\hat\Phi\in
\mathcal{S}(\hat{M})$.
\item \label{item-lema2} For every $\mathcal{S}$-neighborhood
$\hat{\mathscr{U}}$ of $\hat\Phi$ there exists a
$\mathcal{S}^1$-neighborhood $\mathscr{U}$ of $\Phi$ such that for
$\Psi \in \mathscr{U}$,  $\hat\Psi$ belongs to
$\hat{\mathscr{U}}$. Here $\hat\Psi$ denotes the induced $\ell$-th
Grassmannian skew-product on $\hat{\mathcal{M}}$ by $\Psi$.
\item \label{item-lema3} If $\hat\Phi \in \mathcal{PHS}(\hat{M})$ then
$\Phi \in \mathcal{PHS}(\mathbb{R}^c)$ and
$$
    \pi W^{ss}_{loc}(\hat{P}) =
    W^{ss}_{loc}(\pi \hat{P}), \ \ \ \ \
    \pi W^{uu}_{loc}(\hat{P}) = W^{uu}_{loc}(\pi \hat{P}) \ \
    \text{for all $\hat{P} \in \hat{\mathcal{M}}$}
$$
where $\pi : \hat{\mathcal{M}} \to \Sigma \times \mathbb{R}^c$ is
the standard projection.
\end{enumerate}
\end{lem}
\begin{proof}
Assume that $\Phi=\tau\ltimes\phi_\xi$ belongs to
$\mathcal{S}^1(\mathbb{R}^c)$ and let $\hat
\Phi=\tau\ltimes\hat{\phi}_\xi$. Notice that $\hat\phi_\xi$
depends $\alpha$-H\"older on the base point $\xi$. Thus, to prove
that $\hat\Phi$ belongs to $\mathcal{S}(\hat{M})$, it suffices to
show that $\hat\phi_\xi$ are bi-Lipschitz homeomorphisms (with
uniform constants that only depends on $\Phi$). First we need the
following claim:

\begin{claim}
\label{lema:Jairo} Let $T$ be a linear automorphism of
$\mathbb{R}^c$. Then the induced transformation on the
Grassmannian $G(\ell,c)$ is bi-Lipschitz with constant
$\|T\|\|T^{-1}\|$.
\end{claim}
\begin{proof}
$ d_G(T\mathrm{V},TW)=\inf_{BT\mathrm{V} \subset TW}
\|B-i_{T\mathrm{V}}\|=\inf_{A\mathrm{V}\subset W}
\|T(A-i_\mathrm{V})T^{-1}|_{T\mathrm{V}}\|\leq \|T\|\|T^{-1}\| \,
d(\mathrm{V},W). $
\end{proof}

Now, using the triangular inequality and the above claim
\begin{align*}
d_{\hat{M}}(\hat\phi_\xi(x,E),\hat\phi_\xi(x',E')) &=
d(\phi_\xi(x),\phi_\xi(x')) + d_G(D\phi_\xi(x)E,D\phi_\xi(x')E')
\\
&\leq \hat\gamma^{-1} \, d(x,x') + (\gamma\hat\gamma)^{-1} \,
d(E,E') + d_G(D\phi_\xi(x)E',D\phi_\xi(x')E').
\end{align*}
By definition of the Grassmannian metric and since
$\Phi\in\mathcal{S}^1(\mathbb{R}^c)$ it holds that
\begin{equation}
\label{eq.estimate} d_G(D\phi_\xi(x)E',D\phi_\xi(x')E') \leq
\|D\phi_\xi(x)^{-1}\| \, \|D\phi_\xi(x)-D\phi_\xi(x')\| \leq
\gamma^{-1}L \, d(x,x')
\end{equation}
where $L\geq 0$ is a uniform (only depending on $\Phi$) Lipschitz
constant for $D\phi^{\pm 1}_\xi$. Thus,
\begin{equation}
\label{eq:lemma}
d_{\hat{M}}(\hat\phi_\xi(x,E),\hat\phi_\xi(x',E')) \leq
\max\{\hat\gamma^{-1} +\gamma^{-1}L,\,(\gamma\hat\gamma)^{-1}\} \,
d_{\hat{M}}((x,E),(x,E')).
\end{equation}
Analogously we can prove a similar inequality for
$\hat\phi^{-1}_\xi$. Therefore, we have obtained that
$\hat\phi_\xi$ is  bi-Lipschitz (with uniform constants on $\xi$).

 Let
$\Psi=\tau\ltimes\psi_\xi$ be a $\mathcal{S}^1$-perturbation of
$\Phi$. We will show that $\hat\Psi=\tau\ltimes\hat\psi_\xi$ is a
$\mathcal{S}$-perturbation of $\hat\Phi$. Indeed, the
$\mathcal{S}$-distance between $\hat\Phi$ and $\hat\Psi$ is less
or equal than the maximum of
\begin{align*}
 d_{C^0}(\phi^{\pm 1}_\xi,\psi^{\pm 1}_\xi) +
 \max_{(x,E)\in \hat{M}} d_G(D\phi_\xi^{\pm 1}(x)E,D\psi_\xi^{\pm
 1}(x)E) +
  |\mathrm{Lip}(\hat\phi^{\pm 1}_\xi)-\mathrm{Lip}(\hat\psi^{\pm
 1}_\xi)|+|C_1(\Phi)-C_1(\Psi)|.
\end{align*}
Similar estimates as in~\eqref{eq.estimate} and~\eqref{eq:lemma}
show that $d_G(D\phi_\xi(x)E,D\psi_\xi(x)E)\leq \gamma^{-1}
\|D\phi_\xi(x)-D\psi_\xi(x)\|$ and $
|\mathrm{Lip}(\hat\phi_\xi)-\mathrm{Lip}(\hat\psi_\xi)|\leq
\gamma^{-1} |\mathrm{Lip}(D\phi_\xi)-\mathrm{Lip}(D\psi_\xi)|$.
Analogously for $\phi^{-1}_\xi$. Thus,
\begin{align*}
  d_\mathcal{S}(\hat\Phi,\hat\Psi) &\leq (\gamma\hat\gamma)^{-1}
   \, \max_{\xi \in \Sigma} \, \big(
  \sup_{\xi \in \Sigma} d_{C^1}(\phi_\xi,\psi_\xi)
   + |\mathrm{Lip}(D\phi^{\pm 1}_\xi)-\mathrm{Lip}(D\psi^{\pm 1}_\xi)| \big)\\
   &+|C_1(\Phi)-C_1(\Psi)|
   \leq (\gamma\hat\gamma)^{-1} d_{\mathcal{S}^1}(\Phi,\Psi).
\end{align*}
This proves~\eqref{item-lema2}.

Finally, it is clear that if $\hat\Phi \in \mathcal{PHS}(\hat{M})$
then $\Phi \in \mathcal{PHS}(\mathbb{R}^c)$. To prove that the
strong laminations of $\hat\Phi$ project on the strong laminations
of $\Phi$, remember that the strong laminations are H\"older
graphs over the stable sets of the symbolic shift. Then it
suffices to see that the projection of the H\"older graphs in
$\Sigma\times \hat{M}$ is as well an invariant lamination of
$\Sigma\times \mathbb{R}^c$ (i.e., the leaves are invariant
H\"older graphs over the strong set of the shift). But this is
clear by the invariance of the strong lamination on
$\hat{\mathcal{M}}$.
\end{proof}

\begin{proof}[Proof of Theorem~\ref{thm:tangencias}]
First of all, notice that~\eqref{T1},~\eqref{T3} and~\eqref{T4}
are the assumptions of Theorem~\ref{thm-cycles1} and
hence~\eqref{item1} follows immediately from this result when
applied to $\hat\Phi$. To obtain~\eqref{item2}, observe that
\eqref{T3} implies
$$
\pi\Lambda^s(\hat{\mathcal{B}}_2;\hat\Phi)\subset\Lambda^s(\mathcal{B}_2;\Phi)\subset
W^s_{loc}(\Gamma^1).$$
Thus, by
Lemma~\ref{lem-projection}(\ref{item-lema3}) and~\eqref{T4} all
the assumptions of Theorem~\ref{thm-cycles1} are satisfied for the
symbolic skew-product $\Phi$, and this concludes~\eqref{item2}.

Now we will show~\eqref{item3}. To do this we basically need to
reprove~\eqref{item1} but extracting more information. Since
$\hat{P}$ belongs to $\Lambda^s(\hat{\mathcal{B}}_2;\hat\Phi)$
where $\hat{\mathcal{B}}_2$ is a superposition domain of the
$cu$-blender $\hat{\Gamma}_2$ it follows from~\eqref{cond:B3} that
\begin{equation}
\label{eq:pp0}
    W^{ss}_{loc}(\hat{P}) \cap \hat{\mathcal{B}}_2
    \subset \Lambda^s(\hat{\mathcal{B}}_2;\hat\Phi)
    \quad \text{$\mathcal{S}$-robustly.}
\end{equation}
On the other hand, according to~\eqref{T4} and since
$\hat{\mathcal{B}}_1$ is a superposition domain of the
$cs$-blender  $\hat\Gamma_1$ which satisfies~\eqref{cond:B4} then
\begin{equation}
\label{eq:pp1}
     \Lambda^u(\hat{\mathcal{B}}_1;\hat\Phi) \cap
     W^{ss}_{loc}(\hat{\Phi}^{-n}(\hat{P}))\not=\emptyset \quad
     \text{$\mathcal{S}$-robustly.}
\end{equation}
Since $\hat\Phi^m(\hat{P})\in
\Lambda^s(\hat{\mathcal{B}}_2;{\hat\Phi})$ for all $m\geq 0$ then
iterating $\hat{P}$ if necessary we can suppose that $n$ is large
enough so that
$$
\hat\Phi^n(W^{ss}_{loc}(\hat{\Phi}^{-n}(\hat{P}))\subset
W^{ss}_{loc}(\hat{P})\cap \hat{\mathcal{B}}_2.
$$

Notice that~\eqref{eq:pp0} and~\eqref{eq:pp1} implies 
that $\mathcal{S}$-robustly $W^s_{loc}(\hat\Gamma^2)$ meets
$W^u_{loc}(\hat\Gamma^1)$. Similar expressions can be obtained for
$\Phi$ applying Lemma~\ref{lem-projection}(\ref{item-lema3}).
However, a priori, the projection on $\mathcal{S}(\mathbb{R}^c)$
of the $\mathcal{S}$-neighborhood $\hat{\mathscr{U}}$ of
$\hat\Phi$ where \eqref{eq:pp0} and~\eqref{eq:pp1} hold, is not a
$\mathcal{S}$-neighborhood of $\Phi$ and thus~\eqref{item2} cannot
follow from this. Anyway, by
Lemma~\ref{lem-projection}(\ref{item-lema2}), there is a
$\mathcal{S}^1$-neighborhood $\mathscr{U}$ of $\Phi$ so that
$\hat\Psi\in \hat{\mathscr{U}}$ for all $\Phi\in \mathscr{U}$.
Consider $\Psi=\tau\times\psi_\xi$ in $\mathscr{U}$ and the
continuation of $\hat{P}$ for $\hat\Psi$ (continue to call
$\hat{P}$) satisfying~\eqref{eq:pp0}. Let $\hat{Q}=(\xi,(x,E))$ be
an intersection point of~\eqref{eq:pp1} for $\hat\Psi$.  Then $
\hat\Psi^{-m}(\hat{Q})\in \hat{\mathcal{B}}_1$ and
$\hat\Psi^{n+m}(\hat{Q}) \in \hat{\mathcal{B}}_2$ for all $m\geq
0$. In particular, by~\eqref{T3}, it follows that
\begin{equation*}
\label{eq:Euu} \text{$\Psi^{-m}(\xi,x)\in \mathcal{B}_1$ and
$D\psi^{-m}_\xi(x)E \in \mathcal{C}^{uu}$ for all $m\geq 0$.}
\end{equation*}
Since $\Gamma^1$ is a $cs$-blender with superposition domain
$\mathcal{B}_1$ then $\Lambda^u(\mathcal{B}_1;\Psi)\subset
W^u_{loc}(\Gamma^1_\Psi)$ and hence $(\xi,x)\in
W^u_{loc}(\Gamma^1_\Psi)$. A similar argument shows that
$\Psi^n(\xi,x)\in W^s_{loc}(\Gamma^2_\Psi)$, and proves that
\begin{equation*}
\label{eq:us} (\xi,x)\in W^u(\Gamma^1_\Psi)\cap
W^s(\Gamma^2_\Psi).
\end{equation*}

On the other hand, by the $\mathcal{S}^0$-robustness of unstable
cones, $\mathcal{C}^{uu}$ is also an unstable $\ell$-cone for
$\Psi$ on $\mathcal{B}_1$ and hence, by
Lemma~\ref{lem-unstable-cone}, since $(\xi,x)\in
\Lambda^u(\mathcal{B}_1;\Psi)$,
$$
\|D\psi_\xi^{-m}(x)v\| \leq \lambda^m\|v\| \quad \text{for all
$m\in \mathbb{N}$ and $v\in E$.}
$$
Analogous argument proves that $\|D\psi_\zeta^m(y)w\|\leq
\lambda^m \|w\|$ for all $m\in\mathbb{N}$ and $w\in
F=D\psi^n_\xi(x)E$ where $(\zeta,y)=\Psi^n(\hat{Q})$. Then, there
is a constant $C>0$ such that $\|D\psi^{m+n}_\xi(x)v\| \leq
C\lambda^m \|v\|$ for all $m\in\mathbb{N}$ and $v\in E$. Adjusting
the constants if necessary we get $C>0$ and $0<\lambda<1$ such
that
$$
   \|D\psi^m_\xi(x)v\| \leq C \lambda^{|m|} \quad  \text{for all
   $m\in \mathbb{Z}$, \,$v\in E$ with $\|v\|=1$.}
$$
Finally, since the dimension of $E$ is equal to $\ell$ and
$\ell>i_2-i_1$ then
$$
 \ell=d_T(\xi,x) \quad \text{and} \quad
 \mathrm{ind^{cu}}(\Gamma^1) +  \mathrm{ind^{cu}}(\Gamma^2) - \ell =
 c-i_1 + i_2 - \ell  < c.
$$
Thus $\Phi$ has a tangency of dimension $\ell$ between the
unstable set of $\Gamma^1_\Psi$ and the stable set of
$\Gamma^2_\Psi$. This proves~\eqref{item3} concluding the theorem.
\end{proof}

If $\Gamma^1$ is equals to $\Gamma^2$ in
Theorem~\ref{thm:tangencias} then $\Gamma\equiv \Gamma^1=\Gamma^2$
is a $double$-blender. In this case the conclusion is that
$\Gamma$ has a $\mathcal{S}^1$-robust homoclinic tangency.
Similarly, the previous theorem allows us to construct a
$\mathcal{S}^1$-robust tangency on an equi/heterodimensional
cycle:

\begin{cor}
\label{cor:tangencias} Let $\Phi\in
\mathcal{PHS}^1(\mathbb{R}^c)$.  Then,
\begin{enumerate}[leftmargin=0.6cm,itemsep=0.1cm]
\item if $\Gamma$ is a $double$-blender satisfying~\eqref{T1},
\eqref{T2},~\eqref{T3} and~\eqref{T4} in
Theorem~\ref{thm:tangencias} for $\Phi$ and $\Phi^{-1}$ then
$\Phi$ has a $\mathcal{S}^1$-robust homoclinic tangency associated
with $\Gamma$.
\item if $\Gamma^1$ and $\Gamma^2$ are a
$cs$-blender and a $cu$-blender respectively with
$\mathrm{ind}^{cs}(\Gamma^1)\leq \mathrm{ind}^{cu}(\Gamma^2)$,
satisfying~\eqref{T1}, \eqref{T2},~\eqref{T3} and~\eqref{T4} in
Theorem~\ref{thm:tangencias} and the transition
property~\eqref{cond:C2} for $\Phi^{-1}$ then $\Phi$ has a
$\mathcal{S}^1$-robust tangency on an equi/heterodimensional cycle
associated with $\Gamma^1$ and $\Gamma^2$.
\item if $\Gamma^1$ and $\Gamma^2$ are a pair of $double$-blenders satisfying~\eqref{T1},
\eqref{T2},~\eqref{T3} and~\eqref{T4} in
Theorem~\ref{thm:tangencias} for $\Phi$ and $\Phi^{-1}$ then
$\Phi$ has a $\mathcal{S}^1$-robust tangency between both cyclic
intersections of the stable and unstable sets of $\Gamma^1$ and
$\Gamma^2$.
\end{enumerate}
\end{cor}

\section{Constructions: cycles and tangencies from one-step maps}
From now on the alphabet is $\mathscr{A}=\{1,\dots,d\}$ with
$d\geq 2$ and the fiber space $M$ is a differentiable manifold of
dimension $c\geq 1$.
We will give the prototypical construction of blenders, robust
cycles and robust tangencies from a particular class of
skew-products, the so-called, one-step maps:
\begin{defi}[one-step skew-product maps]
A symbolic skew-product $\Phi=\tau\ltimes\phi_\xi$ is called
\emph{one-step} if the fiber maps $\phi_\xi$ only depend on the
coordinate $\xi_0$ of the bi-sequences
$\xi=(\xi_i)_{i\in\mathbb{Z}} \in \Sigma$. In this case we have
$\phi_\xi=\phi_i$ if $\xi_0=i$ and write $\Phi=\tau\ltimes
(\phi_1,\dots,\phi_d)$.
\end{defi}
For a one-step skew-product
$\Phi=\tau\ltimes(\phi_1,\dots,\phi_d)$, one has the underline
dynamics given by the semigroup action generated by
$\phi_1,\dots,\phi_d$ (often referred to as the \emph{iterated
function system} or simply $\IFS$). In what follows, $\langle
\phi_1,\dots,\phi_d \rangle^+$ denotes the semigroup generated by
these bi-Lipschitz homeomorphisms.

\subsection{Blenders from one-step maps}  The notion of a blending region was introduced
simultaneously in~\cite{HN13} and~\cite{BKR14} as a local open set
with robust minimal dynamics for an IFS. 
Here we will introduce an extension of this notion.

\begin{defi}[blending region]
\label{def:blending-region1} Let $\Phi=\tau\ltimes
(\phi_1,\dots,\phi_d)\in \mathcal{PHS}(M)$. Consider  bounded open
sets $B$ and $D$ of $M$ with $\overline{B}\subset D$, a subset
$S\subset \mathscr{A}$, and a
hyperbolic transitive set
$$
\Gamma \eqdef \bigcap_{n\in\mathbb{Z}} \Phi^n(S^\mathbb{Z}\times
\overline{D})=\bigcap_{n\in\mathbb{Z}} \Phi^n(\mathsf{V}\times
\overline{D}),
$$
where $\mathsf{V}$ denotes any isolating neighborhood of
$\Sigma^+_S\eqdef\{\xi\in \Sigma: \xi_0 \in S\}$ and
 $\Sigma^-_S\eqdef\{\xi\in \Sigma: \xi_{-1} \in S\}$.
We say that $B$ is a
\emph{$cs/cu/double$-blending region} with respect to $\{ \phi_i:
i\in S\}$ on $D$  if there exists respectively a
\begin{enumerate}
\item \label{cs:cover} \emph{$cs$-cover}: $\{\phi_i(B):i \in S \}$ is an open cover of $\overline{B}$ and $\mathrm{ind}^{cs}(\Gamma)>0$;
\item \label{cu:cover} \emph{$cu$-cover}: $\{\phi_i^{-1}(B): i \in S\}$ is an open cover of $\overline{B}$ and $\mathrm{ind}^{cu}(\Gamma)>0$;
\item \emph{$double$-cover}: both~\eqref{cs:cover} and~\eqref{cu:cover} are true.
\end{enumerate}
We call \emph{$cs$-index} (resp.~\emph{$cu$-index}) of the
blending region $B$ the $cs$-index (resp.~$cu$-index) of $\Gamma$.
As in the case of the blender, if its $cs$-index
(resp.~$cu$-index) is equal to dimension of $X$ the blending
region is called \emph{contracting} (resp.~\emph{expanding}).
\end{defi}


With the above terminology we get the following corollary of
Theorem~\ref{cs-thm}:
\begin{cor}[blenders from one-step maps]
\label{cor:blender} Let $\Phi=\tau\ltimes(\phi_1,\dots,\phi_d) \in
\mathcal{PHS}(M)$. Assume that there are a finite set $S\subset
\mathscr{A}$ and a
\begin{itemize}
 \item[-]  $cs/cu/double$-blending region $B$ with respect to $\{ \phi_i: i\in S\}$ on $D$.
\end{itemize}
Then the maximal invariant set $\Gamma$ in $S^\mathbb{Z}\times
\overline{D}$ is a $cs/cu/double$-blender  of $\Phi$ whose
super\-po\-si\-tion region contains the family of almost
horizontal discs in $\Sigma^+_S\times B$ or/and almost vertical
discs in $\Sigma^-_S\times B$. Moreover, it also contains the
family of local strong stable/unstable sets, i.e.,~\eqref{cond:B1}
holds.
\end{cor}
\begin{proof}
We assume that $B$ is a $cs$-blending region. The remainder cases are
proved analogously.

From the $cs$-cover~\eqref{cs:cover}, for every $i\in S$ we find an
open set $B_i \subset M$ such that
\begin{equation}
\label{eq:cover0}
 \overline{B} \subset \bigcup_{i\in S} B_i \quad
 \text{and} \quad  \overline{B_i}\subset \phi_i(B).
\end{equation}
This implies the covering property~\eqref{cs-cover} of
Theorem~\ref{cs-thm}. As $\Phi$ is a one-step map then the
H\"older constant $C_0=0$. Thus, it holds that $C=C_0\cdot
(1-\gamma^{-1}\nu^\alpha)^{-1}<L$ and $C\nu^\alpha<\gamma L/2$
where $L>0$ is the Lebesgue number of~\eqref{eq:cover0}. Therefore
according to Theorem~\ref{cs-thm} and
Remark~\ref{rem:local-strong-disc}, $\Gamma$ is a symbolic
$cs$-blender whose superposition region contains the family of
almost horizontal discs and local strong sets in $\Sigma_S^+\times
B$.
This concludes the proof.
\end{proof}

\subsubsection{Hyperbolicity and the Conley-Moser conditions}
To construct symbolic blenders from one-step maps one needs to
know when the invariant set is hyperbolic. We will describe some
sufficient conditions on the fiber maps in order to guarantee this
requirement.

 A simple example of a one-step map with a hyperbolic set where
 the fiber maps have both expanding and contracting directions
 is to take a direct product in the following manner.

\begin{exap}
\label{ex:1}
 Let us to
 consider a finite set  $S\subset \mathscr{A}$, the fiber space
 $M=\mathbb{R}^c$ and a subset
 $D=D_{cs}\times D_{cu}$ of $M$ where $D_{cs}$ and $D_{cu}$ are
 bounded open sets of $\mathbb{R}^{cs}$ and
 $\mathbb{R}^{cu}$ respectively.
 We take two one-step symbolic skew-products
$$
\Phi_{cs}=\tau\ltimes (f_1,\dots,f_d) \in
   \mathcal{S}(\mathbb{R}_{cs}) \quad \text{and} \quad
   \Phi_{cu}=\tau\ltimes (h_1,\dots,h_d) \in
   \mathcal{S}(\mathbb{R}_{cu})
$$
such that  $f_i$ and $h_i^{-1}$ are contracting maps of
$\overline{D_{cs}}$ and $\overline{D_{cu}}$ respectively for all
$i\in S$. Let $\Gamma_{cs}$, $\Gamma_{cu}$ be the maximal
invariant sets in the closure of $\Sigma_S^+ \times D_{cs}$ and
$\Sigma_S^+ \times D_{cu}$ of $\Phi_{cs}$ and $\Phi_{cu}$. These sets are attracting and
repelling continuous invariant sections~\cite{HPS77, BKR14}.
Namely,
they are graphs of continuous maps $g_{cs}: S^\mathbb{Z}
\rightarrow D_{cs}$ and $g_{cu}:S^\mathbb{Z}\rightarrow D_{cu}$
such that
$$
g_{cs} \circ \tau(\xi) = f_{\xi_0} \circ g_{cs}(\xi) \quad
\text{and} \quad g_{cu} \circ \tau(\xi) = h_{\xi_0} \circ
g_{cu}(\xi).
$$
Hence, 
$\Phi_{cs}|_{\Gamma_{cs}}$ and $\Phi_{cu}|_{\Gamma_{cu}}$ are
conjugated to the  shift map $\tau:S^\mathbb{Z}\to S^\mathbb{Z}$.
Take
the direct product
$$
\Phi=\tau\ltimes (\phi_1,\dots,\phi_d) \in
   \mathcal{S}(M), \quad \phi_i(x_{cs},x_{cu})=(f_i(x_{cs}),h_i(x_{cu}))
   \quad \text{for $i=1,\dots,d$}.
$$
The maximal invariant set $\Gamma$ of $\Phi$ in $\Sigma_S^+\times
\overline{D}$ is the graph of the continuous function
$$
   g: S^\mathbb{Z} \to D, \quad g(\xi)=(g_{cs}(\xi),g_{cu}(\xi))
   \ \ \text{with} \ \ g\circ \tau(\xi) = \phi_{\xi} \circ g(\xi).
$$
Thus, $\Phi|_{\Gamma}$ is conjugated to the shift map
$\tau:S^\mathbb{Z}\to S^\mathbb{Z}$ and hence $\Gamma$ is an
isolated hyperbolic transitive set of $\Phi$ with $cs$-index
equal to the dimension of $D_{cs}$.
\end{exap}

When $D$ is a bounded open subset of $\mathbb{R}^c$, $c\geq 2$, the topological criteria for a continuous map $f:\overline{D} \rightarrow \mathbb{R}^c$
to have a hyperbolic set in $D$ are known as the \emph{Conley-Moser conditions}
which are explained in detail in~\cite{Wi88}.
The conditions are given with respect to the so-called horizontal
and vertical (contracting and expanding directions) slabs, which
are fattened up horizontal and vertical Lipschitz graphs in $D$.

A modification of the above example using the Conley-Moser conditions is the following.

\begin{exap}
\label{exa:2} Let $\Phi=\tau\ltimes(\phi_1,\dots,\phi_d)\in
\mathcal{S}(\mathbb{R}^c)$, $S\subset\mathscr{A}$ and
$D=D_{cs}\times D_{cu} \subset \mathbb{R}^{c}$.  Assume that $
\phi^{-1}_i(\overline{D})\cap D$ and $\phi_i^{}(\overline{D})\cap
D $ are both, respectively, a unique horizontal and vertical slabs
$H_i$ and  $V_i$ in $D$ for all $i\in S$. Moreover,
\begin{align*}
&\pi_{cs} \circ \phi_i^{}(\cdot,x_{cu}): \overline{D_{cs}} \to D_{cs} \  \text{is a contracting map for all $x_{cu}\in \overline{D_{cu}}$}, \\
&\pi_{cu} \circ \phi_i^{-1}(x_{cs},\cdot): \overline{D_{cu}} \to D_{cu} \ \text{is a contracting map for all $x_{cs}\in \overline{D_{cs}}$}.
\end{align*}
where $\pi_{*}$ is the projection on the $*$-coordinates, $*\in\{cs,cu\}$. 
Hence $ \mathcal{H}_{i}= \mathsf{H}_i \times H_i$ for $i\in S$
where $\mathsf{H}_i=\{\xi \in \Sigma: \xi_0=i\}$ are the
horizontal slabs for $\Phi$ satisfying the Conley-Moser conditions
and therefore the maximal invariant set in the closure of
$\Sigma_S^+\times D$ for $\Phi$ is conjugated to the full shift
$\tau : S^\mathbb{Z}\to S^\mathbb{Z}$.
\end{exap}

Observe that in general, unlike in the above examples, the hyperbolic invariant set
given by the Conley-Moser conditions does not have to be a graph over the base.
It suffices to modify the above example so that  $\phi_i^{-1}(\overline D)\cap D$ and $\phi_i(\overline D)\cap D$ are
two disjoints horizontal and vertical slabs respectively.

\subsubsection{Construction of blending regions}
We will describe now the method to construct a blending region.
First we work in local coordinates with $C^r$-diffeomorphisms with
$r\geq 0$. A $C^0$-diffeomorphism is understood as a bi-Lipschitz
homeomorphism.

\begin{prop}
\label{prop:creation-blender} Consider a $C^r$-diffeomorphism
$\phi$ of $\mathbb{R}^c$ with a hyperbolic
attrac\-ting/repe\-lling/saddle fixed point $x$. Then, there exist
an integer $k\equiv k(\phi,c)\geq 2$, arcs of
$C^r$-diffeomorphisms of $\mathbb{R}^c$,
$\phi_1\equiv\phi_1(\varepsilon),\dots,\phi_k\equiv\phi_k(\varepsilon)$
and bounded open sets $D\equiv
D(\varepsilon)$,  $\varepsilon\geq 0$, 
such that
\begin{enumerate}[leftmargin=0.6cm]
\item[-] $\phi_i(0)=\phi$ for $i=1,\dots,k$;
\item[-] $\phi_i=T_i \circ \phi$ where $T_i\equiv T_i(\varepsilon)$ is a translation (moreover, one can take $\phi_1=\phi$);
\item[-] $B\equiv B_{\delta}(x) \subset D \subset B_{2\varepsilon}(x)$ for some $\delta\equiv \delta(\varepsilon)>0$; 
\item[-] $B$ is $cs/cu/double$-blending region with respect to $\{\phi_1,\dots,\phi_k\}$ on $D$ for all $\varepsilon>0$.
\end{enumerate}
Moreover, the $cs$-index of the blending region is equal to the $s$-index of the hyperbolic fixed point $x$.
\end{prop}

\begin{proof}
Fix $\varepsilon>0$. Assume first that $x$ is an attracting/saddle
fixed point of $\phi$. We will construct a $cs$-blending region.
Repeating the arguments for $\phi^{-1}$ one constructs
$cu$-blending regions and combining both one concludes the
proposition.

Let $D=D_{cs}\times D_{cu}$ be a small neighborhood of $x$ with
$D_*$ being an open ball of radius $\varepsilon>0$ centered at
$\pi_*(x)$ where $*\in \{cs,cu\}$. In the case that $x$ is an
attracting fixed point we allow that $D_{cu}$ could be empty,
i.e., $D=D_{cs}$. Moreover, assume that
\begin{itemize}
 \item[-] $\pi_{cs} \circ \phi(\cdot,x_{cu}): \overline{D_{cs}} \to D_{cs}$
 is a contracting map for all $x_{cu}\in \overline{D_{cu}}$,
\item[-] $\pi_{cu} \circ \phi^{-1}(x_{cs},\cdot):
\overline{D_{cu}} \to D_{cu}$
 is a contracting map for all $x_{cs}\in \overline{D_{cs}}$.
 \end{itemize}
Consider an open ball $B \subset D$ of radius
$0<\delta<\varepsilon$ centered at $x$. Let us denote by $T_v$ the
translation map by the vector $v$. Applying
\cite[Lemma~5.6]{BKR14} (see also~\cite[Proposition~2.3]{NP12})
for the map $\pi_{cs}\circ\phi(\cdot,x_{cu})$ and since
$\pi_{cu}\circ\phi(x_{cs},\cdot)$ is expanding on $D_{cu}$, there
exist $k\in\mathbb{N}$, vectors $u_i$ in the ball of radius 1 so
that
$$
    \overline{B} \subset \bigcup_{i=1}^k  T_{v_i} \circ \phi(B) \quad \text{where $v_i=\delta u_i$ for $i=1,\dots,k$.}
$$
In fact, the number $k$ of translation and the direction of
translation $u_1,\dots,u_k$ only depend on the (uniform)
contraction lower bound of the maps $\pi_{cs}\circ
\phi(\cdot,x_{cu})$ and the dimension $c$. Having in mind
Example~\ref{exa:2} and taking $\phi_i=T_{v_i}\circ \phi$ one
concludes that $B$ is a $cs$-blending region with respect to
$\{\phi_1,\dots,\phi_k\}$ on $D$.

To conclude the proposition is
enough to observe that $v_i(\varepsilon)$ tends
continuously to zero as
$\varepsilon\to 0$. Moreover, we can assume that $\phi_1(\varepsilon)=\phi$ for all
$\varepsilon\geq0$. Consequently
$\phi_1(\varepsilon),\dots,\phi_k(\varepsilon)$, $\varepsilon\geq 0$
are arcs of
$C^r$-diffeomorphisms satisfying the required conditions.
\end{proof}

For every $x\in M$, by means of a arbitrarily small perturbation
of the identity map~\cite{H02,HN13} we can create a map $\phi$
for which $x$ is a hyperbolic fixed point (or periodic of arbitrary large period).
Hence, the above result implies the following:
\enlargethispage{1cm}
\begin{cor}[blending regions homotopic to the identity]
\label{rem:creation-blender}
For every $x\in M$, there exist arcs of $C^r$-diffeomorphisms
$\phi_1\equiv\phi_1(\varepsilon),\dots,\phi_k\equiv\phi_k(\varepsilon)$
of $M$ (where $k$ only depends on the dimension of $M$), bounded
open sets $D\equiv D(\varepsilon)\subset B_{2\varepsilon}(x)$ and
$\delta\equiv \delta(\varepsilon)>0$, for $\varepsilon\geq 0$,
such that
\begin{itemize}
 \item[-] $\phi_i(0)=\mathrm{id}$ for $i=1,\dots,k$ and
 \item[-] $B\equiv B_{\delta}(x)\subset D$ is a blending region with respect to $\{\phi_1,\dots,\phi_k\}$ on $D$ for all $\varepsilon>0$.
\end{itemize}
Moreover, the blending region can be constructed having any $cs$-index between
0 and dimension of~$M$.
\end{cor}


\subsection{Robust cycles from one-step maps}
Theorem~\ref{thm-cycles1} provides a criterion to yield robust
cycles using blenders. Now we translate this criterion to the
dynamics of a one-step skew-product and show how to
construct arcs of one-step maps satisfying this criterion.

\subsubsection{Criterion to yield robust cycles from one-step maps} Let
$\Phi=\tau\ltimes(\phi_1,\dots,\phi_d) \in \mathcal{PHS}(M)$ be a
one-step map. We need a $cs$-blender for $\Phi$. To accomplish
this, from Corollary~\ref{cor:blender}, we can assume that there
is a $cs$-blending region $B$ with respect to
$\{\phi_1,\dots,\phi_d\}$ on $D$. Hence the maximal invariant set
$\Gamma^1$ in the closure of $\Sigma\times D$ for $\Phi$ is a
$cs$-blender with superposition domain $\mathcal{B}=\Sigma \times
B$. Now, we need to translate the transition
property~\eqref{cond:C2}. To do this, let $\Gamma^2$ be another
isolated transitive hyperbolic set of $\Phi$ satisfying the
following condition:

\begin{enumerate}[label=(RC2), ref=RC2]
 \item \label{cond:C2''}  there are $x\in \mathscr{P}(W^s(\Gamma^2))$
and  $T \in \langle\phi_1,\dots,\phi_d\rangle^+$ such that $T^{-1}(x) \in B$.
\end{enumerate}
Take any $\xi \in \Sigma$ such that $T^{-1}= \phi^{-n}_{\tau^{-1}(\xi)}=\phi^{-1}_{\xi_{-n}}\circ \dots \circ \phi_{\xi_{-1}}^{-1}$ for some $n\geq 1$. Then,
$$
   W^{ss}_{loc}(\Phi^{-n}(\xi,x))= W^{s}_{loc}(\tau^{-n}(\xi))\times \{T^{-1}(x)\} \subset \Sigma \times B.
$$
Since $\Gamma^1$ is a $cs$-blender satisfying~\eqref{cond:B1}, the
above inclusion implies $W^{ss}(\xi,x) \cap W^u(\Gamma^1)
\not=\emptyset$. Since $x$ belongs to the projection on the fiber
space of $W^s(\Gamma^2)$, one can choose $\xi$  such that
$W^{ss}(\xi,x) \subset W^s(\Gamma^2)$. Therefore~\eqref{cond:C2''}
is equivalent to~\eqref{cond:C2} and by Theorem~\ref{thm-cycles1}
we get that $\mathcal{S}$-robustly the global stable set
$W^s(\Gamma^2)$ meets the global unstable set $W^u(\Gamma^1)$.

A similar condition to the above guarantees the robustness
of the other intersection. 

\begin{defi}[transition]
\label{def:transition} Let $A_1$ and $A_2$ be two subsets of $M$.
We say that $\langle \phi_1,\dots,\phi_d\rangle^+$ has
\emph{transition from $A_1$ to $A_2$} if there exist $x\in A_1$
and $T\in \langle \phi_1,\dots,\phi_d\rangle^+$ such that $T(x)\in
A_2$.
\end{defi}

With this terminology we have obtained the following consequence
of Theorem~\ref{thm-cycles1}:

\begin{cor}[robust cycles form one-step maps]
\label{cor:ciclos} Let $\Phi=\tau\ltimes(\phi_1,\dots,\phi_d) \in
\mathcal{PHS}(M)$. Suppose that there are bounded open sets $B_1
\subset D_1$ and $B_2\subset D_2$  with $D_1\cap D_2 =\emptyset$
such that
\begin{itemize}
\item[-] $B_1$ is a $cs$-blending region with respect to
$\{\phi_1,\dots,\phi_d \}$ on $D_1$;
\item[-] $B_2$ is a $cu$-blending region with respect to
$\{\phi_1,\dots,\phi_d \}$ on $D_2$;
\item[-] $\langle\phi_1,\dots,\phi_d \rangle^+$ has transition from
$B_1$ to $B_2$ and vice-versa.
\end{itemize}
Then $\Phi$ has a robust cycle associated with a symbolic
$cs$-blender (the maximal invariant set in $\Sigma\times
\overline{D_1}$) and a symbolic $cu$-blender (the maximal
invariant set in $\Sigma\times \overline{D_2}$).
\end{cor}

A point $x\in X$ is called a \emph{periodic point} of $\langle\phi_1,\dots,\phi_d\rangle^+$ if there is $f \in \langle\phi_1,\dots,\phi_d\rangle^+$
such that $f(x)=x$. Observe that the periodic points of $\langle\phi_1,\dots,\phi_d\rangle^+$
are the projection on the fiber space of the periodic points of $\Phi=\tau\ltimes(\phi_1,\dots,\phi_d)$.

\begin{exap}[Robust cycles with non-hyperbolic periodic points]
Suppose $B$ is a $double$-blending region with respect to
$\{\phi_1,\dots,\phi_d\}$ on $D$ and $x \in M\setminus D$ is a
periodic point of $\langle\phi_1,\dots,\phi_d\rangle^+$ such that
there are $T,S\in \langle\phi_d,\dots,\phi_d\rangle^+$ so that
$T^{-1}(x), S(x) \in B$. Then, one gets a symbolic cycle
associated with a $double$-blender and a periodic point. Notice
that this criteria allows us to construct a robust cycle between a
$double$-blender and \emph{any} periodic point that admits
continuation  under small perturbations.
\end{exap}

%

\subsubsection{Construction of blending regions with transition}
We will now explain
how to construct examples of blending regions
with transition on the manifold $M$.  To do this,
we will work in local coordinates providing
a set of $C^r$-diffeomorphisms $\{\phi_1,\dots,\phi_d\}$ from $\mathbb{R}^c$
to itself, as well as disjoint bounded open sets $D_1$ and $D_2$ of $\mathbb{R}^c$
and open sets $B_1 \subset D_1$, $B_2\subset D_2$  such that
the conditions in Corollary~\ref{cor:ciclos} hold.

Let $h$ be a $C^r$-diffeomorphism of $\mathbb{R}^c$ with two
different hyperbolic fixed points $p$ and $q$ so that $ W^s(p;h)
\cap W^u(q;h) \not = \emptyset$. By means of
Proposition~\ref{prop:creation-blender} we get arcs of
$C^r$-diffeomorphisms, $\phi_{2}\equiv\phi_{2}(\varepsilon),\dots,
\phi_{k}\equiv\phi_{k}(\varepsilon)$ and disjoint bounded open
sets $D_1\equiv D_1(\varepsilon)$, $D_2\equiv D_2(\varepsilon)$
so that $\phi_i(0)=h$ and $B_1\equiv
B_{\delta}(p)$, $B_2\equiv B_{\delta}(q)$ are, respectively, a
$cs$ and $cu$-blending regions for
$\{h,\phi_2,\dots,\phi_k\}$ on $D_1$ and $D_2$. Since the stable
set of $p$ and unstable set of $q$ for $h$ have non-empty
intersection, there is a point $y$ arbitrarily close to $q$ such
that $h^n(y)$ converges to $p$. Hence for any $\varepsilon>0$
small enough one can assume that $y\in B_1$ and $h^n(y)\in B_2$.
This provides a transition from $B_2$ to $B_1$. To give the other
transition we consider the map $T(x)=x+q-p$. Since this additional
map is a translation so that $T(D_1)\cap D_1=\emptyset$ one easily
gets that
\begin{enumerate}
 \item[-]  $B_1$ is a $cs$-blending region with respect to $\{h,\phi_2,\dots,\phi_k,T\}$ on $D_1$;
 \item[-]  $B_2$ is a $cu$-blending region with respect to $\{h,\phi_2,\dots,\phi_k,T\}$ on $D_2$;
 \item[-]  $\langle h,\phi_2,\dots,\phi_{k}, T \rangle^+$ has transition from $B_1$ to $B_2$ and vice-versa.
\end{enumerate}

Notice that the above construction can be done homotopic to the
identity taking an arc of $C^r$-diffeomorphisms
$h_\epsilon:\mathbb{R}^c\to \mathbb{R}^c$ as above with the
additional properties that $h_\epsilon$ goes to the identity map
and $\|p(\epsilon)-q(\epsilon)\|\to 0$ as $\epsilon \to 0$.
Therefore we have obtained the following.

\begin{prop}[blending regions with transition homotopic to the identity]
\label{prop:blending-transition} For every $x\in M$, there exist
arcs of $C^r$-diffeomorphisms $\phi_1\equiv
\phi_1(\varepsilon),\dots,\phi_{k+1}\equiv
\phi_{k+1}(\varepsilon),\ \varepsilon\geq 0$ of $M$ where $k$ only
depends on the dimension of $M$ and disjoints bounded open sets
$D_1\equiv D_1(\varepsilon)$ and $D_2\equiv D_2(\varepsilon)$ in a
small neighborhood of $x$ such that
\begin{itemize}
\item[-]  $\phi_i(0)=\mathrm{id}$ for all $i=1,\dots,k+1$,
\item[-]  there is a $cs$-blending region $B_1$ with respect to $\{\phi_1,\phi_2,\dots,\phi_{k+1}\}$ on $D_1$;
\item[-]  there is a $cu$-blending region $B_2$ with respect to $\{\phi_1,\phi_2,\dots,\phi_{k+1}\}$ on $D_2$;
\item[-]  $\langle \phi_1,\dots,\phi_{k+1} \rangle^+$ has transition from $B_1$ to $B_2$ and vice-versa.
\end{itemize}
Moreover, both blending regions can be constructed having any co-index
(difference between the $cs$-indices) between zero and the dimension of $M$.
\end{prop}


\subsubsection{Arcs of one-step maps with robust cycles}
Previous constructions allow us to construct arcs of symbolic
skew-product with robust cycles. Moreover, we obtain these cycles
for \emph{robust  topologically mixing symbolic skew-products}.
That is, for symbolic skew-products $\Phi$ such that for any small
enough $\mathcal{S}^0$-perturbation $\Psi$ and for every pair of
open sets $\mathcal{U}$, $\mathcal{V}$ of $\mathcal{M}$, there is
$n_0> 0$ such that $\Psi^n(\mathcal{U}) \cap \mathcal{V}
\not=\emptyset$ for all $n \geq n_0$. In particular topologically
mixing implies transitivity.

\begin{thm}
\label{thm:mainA-symbolic-setting}
  Let $M$ be a differentiable
  manifold of dimension $c\geq 1$.
  Then, there are an integer
  $d\equiv d(c)\geq 3$ and an arc of one-step maps
  $\Phi_\varepsilon=\tau\ltimes(\phi,\dots,\phi_d) \in \mathcal{PHS}^{r}(M)$,
  $\varepsilon\geq
  0$, such that $\Phi_0=\tau\times \mathrm{id}$ and for every
  $\varepsilon>0$,
  \begin{itemize}
   \item[-] $\Phi_\varepsilon$ is $\mathcal{S}^0$-robustly topologically mixing;
   \item[-] $\Phi_\varepsilon$ has a $\mathcal{S}$-robust (heterodimensional or equidimensional)  cycle.
  \end{itemize}
Moreover, the arc can be taken in such a way that
$\Phi_\varepsilon$ has robust cycles of all intermediate co-index.
\end{thm}

\begin{proof}
Consider a point $x\in M$. By
Proposition~\ref{prop:blending-transition} we get arcs of
$C^r$-diffeomorphisms $\phi_1\equiv
\phi_1(\varepsilon),\dots,\phi_{k+1}\equiv
\phi_{k+1}(\varepsilon)$ homotopic to the identity as $\varepsilon
\to 0^+$, where $k\geq 2$ and only depends on $c$. We obtain as
well the $cs$ and $cu$ blending regions $B_1$ and $B_2$
respectively  in $B_\varepsilon(x)$ with transitions. Without loss
of generality, assume that $B_1$ is a $contracting$-blending
region. We can add to these maps, if necessary, a pair of
Morse-Smale diffeomorphisms (homotopic to the identity) without
periodic points in common and having respectively an
attracting/repelling fixed point in $B_1$ with a dense
stable/unstable manifold. The existence of such a map can be
constructed by perturbations of time-one maps of gradient-like
vector fields. Take $d=k+1+2 \geq 3$, which only depends on $c$,
and set $\Phi_\varepsilon=\tau \ltimes (\phi_1,\dots,\phi_{d}) \in
\mathcal{PHS}^r(M)$. Hence $\Phi_0=\tau\times\mathrm{id}$ and
according to~\cite[Theorem~5.7]{BKR14} and
Corollary~\ref{cor:ciclos}, for any $\varepsilon>0$,
$\Phi_\varepsilon$ is $\mathcal{S}^0$-robustly topologically
mixing and has a $\mathcal{S}$-robust symbolic cycle for any
$\varepsilon>0$. Notice that the co-index of the cycle may be
chosen between zero and~$c$. In fact, repeating the above
procedure and adding more maps if necessary, the arcs may be taken
in such a way so that $\Phi_\varepsilon$ has $\mathcal{S}$-robust
cycles of all intermediate co-indices between zero and $c$. This
completes the proof of the theorem. \end{proof}

\subsection{Robust tangencies from one-step maps}
\label{ss:one-step-tangencias} In Theorem~\ref{thm:tangencias} we
gave conditions~\eqref{T1}, \eqref{T2},~\eqref{T3} and~\eqref{T4}
to construct robust tangencies. We will translate these conditions
to one-step  maps (Corollary~\ref{cor:tangencias-IFS}) and then we
will construct arcs of $\IFS$s unfolding from the identity
satisfying the criterion that imply robust tangencies
(Proposition~\ref{prop:blending-tangencias}).

\subsubsection{Criterion to yield robust tangencies from one-step
maps} Let $\Phi=\tau\ltimes(\phi_1,\dots,\phi_d) \in
\mathcal{PHS}^1(\mathbb{R}^c)$ be a fiber bunched one-step map.
That is,
$$
   \nu^\alpha < \gamma<
   m(D\phi_i)<\|D\phi_i\|<\hat\gamma^{-1}<\nu^{-\alpha} \quad \text{with}
   \quad \gamma\hat\gamma > \nu^\alpha \quad \text{for all $i=1,\dots,d$. }
$$
We take integers $0<i_1,i_2<c$ and $\max\{0,i_2-i_1\}<\ell\leq
\min \{c-i_1,i_2\}$.

In order to translate~\eqref{T1} for one-step maps, we need the
following lemma:

\begin{lem}
If $\Phi=\tau\ltimes\phi_\xi\in \mathcal{PHS}^1(\mathbb{R}^c)$ is
fiber bunched with the uniform Lipschitz constant of $D\phi^{\pm
1}_\xi$ small enough then the induced $\ell$-th Grassmannian
symbolic skew-product $\hat\Phi$ on $\mathcal{M}=\Sigma\times
\hat{M}$ belongs to $\mathcal{PHS}(\hat{M})$, where
$\hat{M}=\mathbb{R}^c\times G(\ell,c)$, with $1\leq \ell\leq c$.
\end{lem}
\begin{proof}
Clearly $\hat\Phi$ is a $\alpha$-H\"older symbolic skew-product.
 According to~\eqref{eq:lemma}, since $\Phi$ is fiber
bunched ($\nu^{\alpha}<\gamma\hat\gamma$) and taking the Lipschitz
constant $L$ of $D\phi_\xi$ less than
$\gamma(\nu^{-\alpha}-\hat\gamma^{-1})$, it follows that the
Lipschitz constant of $\hat\phi_\xi$ satisfies $\max\{
\hat\gamma^{-1}+\gamma^{-1}L, (\gamma\hat\gamma)^{-1}\} <
\nu^{-\alpha}$. A similar inequality holds for
$\hat\phi^{-1}_\xi$, implying that
$\hat\Phi=\tau\ltimes\hat\phi_\xi \in \mathcal{PHS}(\hat{M})$.
\end{proof}

In view of the above lemma, we will assume that the Lipschitz
constant of $D\Phi$, i.e., of $D\phi^{\pm 1}_i$ is small enough
for all $i=1,\dots,d$. Then
$\hat\Phi=\tau\ltimes(\hat\phi_1,\dots,\hat\phi_d)$
satisfies~\eqref{T1} where
$$
\hat\phi_i(x,E)=(\phi_i(x),D\phi_i(x)E), \ \ \ (x,E)\in
\hat{M}\eqdef \mathbb{R}^c\times G(\ell,c), \quad \text{for
$i=1,\dots,d$}.
$$

Now, we will suppose that there are  $cs$ and $cu$ blending
regions $B_1$ and $B_2$
for $\{\phi_1,\dots,\phi_d\}$ with $cs$-indices $i_1$ and $i_2$
respectively. According to Corollary~\ref{cor:blender}, these
conditions provide a $cs$-blender $\Gamma^1$ and a $cu$-blender
$\Gamma^2$ with superposition domains $\mathcal{B}_1=\Sigma \times
B_1$ and $\mathcal{B}_2=\Sigma\times B_2$ respectively. We will
also assume that there are an unstable $\ell$-cone
$\mathcal{C}^{uu}$ and a stable $\ell$-cone $\mathcal{C}^{ss}$ for
$\Phi$ on $\mathcal{B}_1$ and $\mathcal{B}_2$ respectively. Under
these assumptions we have that~\eqref{T2} holds.

To translate property~\eqref{T3} for the one-step case, it
suffices to ask that there are a $cs$-blending region $\hat B_1$
and a $cu$-blending region $\hat B_2$ with respect to
$\{\hat\phi_{1},\dots,\hat\phi_{d}\}$ such that $\hat B_1 \subset
B_1\times \mathcal{C}^{uu}$ and $\hat B_2 \subset B_2\times
\mathcal{C}^{ss}$. Notice that property~\eqref{cond:B4} is
satisfied since we are working with blenders constructed from the
covering property.

Finally, having into account properties~\eqref{cond:B5},
\eqref{cond:C2''} and Definition~\ref{def:transition}, \eqref{T4}
is equivalent in this case to the existence of a transition for
$\langle \hat\phi_1,\dots,\hat\phi_d\rangle^+$ from $\hat B_1$ to
$\hat B_2$.

\begin{defi}[blending region with tangency]
Let $B$, $D$ be bounded open sets of $\mathbb{R}^c$ and
$0<\ell<c$. We say that a $cs$-blending region, $B$, with respect
to $\{\phi_1,\dots,\phi_d\}$ on $D$ has a \emph{tangency} of
dimension $\ell$  if there exist bounded open sets $\hat{B}$ and
$\hat{D}$ in $\hat{M}=\mathbb{R}^c \times G(\ell,c)$ such that
$$
\text{$\hat{B}$ is a $cs$-blending region with respect to
$\{\hat\phi_1,\dots,\hat\phi_d\}$ on $\hat D$ so that $\hat B
\subset B \times \mathcal{C}^{uu}$}.
$$
Here $\mathcal{C}^{uu}$ is an unstable $\ell$-cone  for
$\Phi=\tau\ltimes(\phi_1,\dots,\phi_d)$ on $\Sigma\times B$.  We
say that the set $\hat{B}$ is a \emph{$\ell$-tangency} (of the
blending region $B$). Similarly, we define a $cu$-blending region
with a tangency of dimension $\ell$.
\end{defi}

Using this terminology, we summarize Theorem~\ref{thm:tangencias}
in the one-step setting:

\begin{cor}[robust tangencies from one-step maps]
\label{cor:tangencias-IFS} Let
$\Phi=\tau\ltimes(\phi_1,\dots,\phi_d) \in
\mathcal{PHS}^1(\mathbb{R}^c)$ be fiber bunched one-step map with
small enough Lipschitz constant of $D\phi^{\pm 1}_i$. Consider
integers
$$
0<i_1,i_2<c \quad \text{and} \quad \max\{0,i_2-i_1\}<\ell\leq \min
\{c-i_1,i_2\}.
$$
Suppose that there are bounded open sets $D_1, D_2, B_1, B_2$ of
$\mathbb{R}^c$ such that with respect to
$\{\phi_1,\dots,\phi_d\}$,
\begin{enumerate}[leftmargin=0.6cm]
\item[-] $B_1$ is a $cs$-blending region on
$D_1$,  with $cs$-index $i_1$ and a $\ell$-tangency $\hat{B}_1$;
\item[-] $B_2$ is a $cu$-blending region
on $D_2$, with $cs$-index $i_2$ and a $\ell$-tangency $\hat{B}_2$;
\item[-]$\langle\hat\phi_1,\dots,\hat\phi_d \rangle^+$ has a transition from
$\hat{B}_1$ to $\hat{B}_2$.
\end{enumerate}
Then $\Phi$ has a $\mathcal{S}^1$-robust tangency of dimension
$\ell$ between $W^u(\Gamma^1;\Phi)$ and $W^s(\Gamma^2;\Phi)$ where
$\Gamma^1$ and $\Gamma^2$ are the maximal invariant sets in
$\Sigma\times \overline{D_1}$ and $\Sigma\times \overline{D_2}$.
%
%
%
%
\end{cor}


\subsubsection{Construction of blending regions with tangencies}
\label{ss:IFS-tangencias} We will now explain how to construct
examples of blending regions with tangencies.


\begin{prop}[blending region with tangency]
\label{lema:tang} Let $p$ be a hyperbolic fixed point of a
$C^r$-diffeomorphism $\phi$ of $\mathbb{R}^c$ with $r\geq 2$ and
having a dominated splitting $E^s\oplus E^{cu}\oplus E^{uu}$ where
the unstable direction is $E^u=E^{cu}\oplus E^{uu}$ and
$0<\ell=\dim E^{uu}<c$. Then, there are
\begin{enumerate}
\item an unstable $\ell$-cone $\mathcal{C}^{uu}$ around $E^{uu}$,
\item neighborhoods $B$ of $p$ in $\mathbb{R}^c$ and $G\subset \mathcal{C}^{uu}$ of $E^{uu}$ in
$G(\ell,c)$, and
\item a finite set of $C^0$-diffeomorphisms
$\{\hat\phi_1,\dots,\hat\phi_d\}$ on $\hat{M}=\mathbb{R}^c\times
G(\ell,c)$ induced by $\phi_i=A_i\cdot \phi+c_i$, where $A_i$ is
an orthogonal matrix and $c_i\in \mathbb{R}^c$,
\end{enumerate}
such that
$$
\overline{\hat B}\subset\bigcup_{i=1}^d \hat\phi_i(\hat B) \quad
\text{where} \quad \hat B = B \times G.
$$
That is, $B$ is a $cs$-blending region with respect to
$\{\phi_1,\dots,\phi_d\}$ of $cs$-index $\dim E^s$ and a
$\ell$-tangency. Moreover, the maps $\phi_1,\dots,\phi_d$ can be
constructed as an arc of $C^r$-diffeomorphisms homotopic to
$\phi$.
 \end{prop}
\begin{proof}
 As in Proposition~\ref{prop:creation-blender}, we can take
 a small neighborhood $B$ of $p$ in $\mathbb{R}^c$
 and choose arcs of translations $f_i=\phi+c_i$ homotopic to $\phi$ with
 $c_i\in\mathbb{R}^c$ in order to have that $B$ is
 a $cs$-blending region with respect to $\{f_1,\dots,f_{k_1}\}$. In
 particular, $\{ f_i(B): i=1,\dots,k_1\}$ is an open cover of
 $\overline{B}$. Notice that there exists an $\epsilon>0$,
 such that the above covering property holds for any
 $\epsilon$-perturbation $\phi_i$ of $f_i$.

 On the other hand, $D\phi(p)$ induces a map
 $A$ on $G(\ell,c)$ given by $A(E)=D\phi(p)E$.
 Since $\phi$ is $C^2$, this map
 has a hyperbolic attracting fixed point $E^{uu}$.
 Then, as in Proposition~\ref{prop:creation-blender},
 there are arcs $T_1\equiv T_1(\varepsilon),\dots,T_{k_2}\equiv T_{k_2}(\varepsilon)$
 of translations on $G(\ell,c)$ homotopic to the identity and
 a neighborhood $G\equiv G(\varepsilon)$ of
 $E^{uu}$ in $G({k,c})$, such that
 $
 \{F_j(G): j=1,\dots,k_2\}$
 is
 an open cover of $\overline{G}$ where $F_j= T_j \circ A$.
 Each translation map $F_j$ in the Grassmannian corresponds to a
map of the form $A_j\cdot D\phi(p)$ in the tangent bundle where
$A_j$ is an orthogonal matrix and $\|A_j- \mathrm{id}\|\to 0$ as
$\varepsilon\to 0$ for all $j=1,\dots,k_2$. Moreover, we can take
$G$ small enough so that $G\subset \mathcal{C}^{uu}$ where
$\mathcal{C}^{uu}$ is a unstable $\ell$-cone around $E^{uu}$ for
$\phi$. Taking
$$\|A_j-\mathrm{id}\|< \epsilon / \|\phi\| \quad \text{and}\quad
\phi_{ij}=A_j \cdot \phi + c_i$$ then $\|\phi_{ij}-f_i\|<\epsilon$
and so
$$
\overline{B} \subset \bigcup_{i=1}^{k_1}
\bigcup_{j=1}^{k_2}\phi_{ij}(B).
$$
Finally, by construction, the maps $\phi_{ij}$ induce a set of
maps $\hat{\phi}_{ij}$ on $\hat{M}$ so that
$\{\hat{\phi}_{ij}(\hat{B})\}$ is an open cover of the closure of
$\hat{B}=B\times G$. Moreover, since $A_j$ is an orthogonal matrix
then the hyperbolicity of the blending regions still holds. That
is, $B$ is also a $cs$-blending region with respect to
$\{\phi_{ij}: i=1,\dots,k_1, \ j=1,\dots,k_2\}$ with $cs$-index
equal to $\dim E^s$ and as well is a $\ell$-tangency (the
$cs$-blending region $\hat{B}$ on $\hat{M}$). This completes the
proof.
\end{proof}

Notice that the above construction can be done homotopic to the
identity by means of Corollary~\ref{rem:creation-blender}. Also we
can construct a transition homotopic to the identity as was done
in Proposition~\ref{prop:blending-transition}. Therefore we get
the following:

\begin{prop}[blending regions with tangencies and transition homotopic to the identity]
\label{prop:blending-tangencias} Consider integers $0<i_1,i_2<c$
and $\max\{0,i_2-i_1\}<\ell\leq \min \{c-i_1,i_2\}$. For every
$x\in \mathbb{R}^c$, there are arcs of $C^r$-diffeomorphisms
$\phi_1\equiv
\phi_1(\varepsilon),\dots,\phi_{d}\equiv\phi_{d}(\varepsilon),\
\varepsilon\geq 0$ where $r\geq 2$ and $d\equiv d(c)\geq 2$, and
bounded open sets $B_1\equiv B(\varepsilon)$, $B_2\equiv
B_2(\varepsilon)$, $D_1\equiv D_1(\varepsilon)$ and $D_2\equiv
D_2(\varepsilon)$ in a small neighborhood of $x$ such that with
respect to $\{\phi_1,\dots,\phi_d\}$,
\begin{enumerate}[leftmargin=0.6cm,itemsep=0.05cm]
\item[-] $B_1$ is a $cs$-blending region
on $D_1$,  with $cs$-index $i_1$ and a $\ell$-tangency
$\hat{B}_1$;
\item[-] $B_2$ is a $cu$-blending region
on $D_2$, with $cs$-index $i_2$ and a $\ell$-tangency $\hat{B}_2$;
\item[-]$\langle\hat\phi_1,\dots,\hat\phi_d \rangle^+$ has transition from
$\hat{B}_1$ to $\hat{B}_2$.
\end{enumerate}
\end{prop}

\subsubsection{Arcs of one-step maps with robust tangencies}
We prove Theorem~\ref{thmF} in the~symbolic~setting:

\begin{thm}
\label{thm:mainC-symbolic-setting}
  Let $M$ be a  differentiable
  manifold of dimension $c>1$.
  Consider integers
  $$
  0<i_1,i_2 < c \quad \text{and} \quad
         \max\{0,i_2-i_1\}< d_T \leq \min\{c-i_1,i_2\}.
  $$
  Then, there are an integer
  $d\equiv d(c)\geq 3$ and arcs of one-step maps
  $\Phi_\varepsilon=\tau\ltimes(\phi_1,\dots,\phi_d) \in
  \mathcal{PHS}^{r}(M)$, $\varepsilon\geq
0$, such that $\Phi_0=\tau\times \mathrm{id}$ and
  for every $\varepsilon>0$,
  \begin{itemize}
   \item[-]  there are a $cs$-blender $\Gamma^1_\varepsilon$ and a $cu$-blender
   $\Gamma^2_\varepsilon$ of $\Phi_\varepsilon$ of indices $i_1$ and $i_2$ respectively,
   \item[-] $\Phi_\varepsilon$ has a $\mathcal{S}^1$-robust
   tangency of dimension $d_T$ between $W^u(\Gamma^1_\varepsilon)$ and
   $W^s(\Gamma^2_\varepsilon)$.
   \end{itemize}
Moreover, these arcs can be taken in such a way so that
$\Phi_\varepsilon$ also has a $\mathcal{S}^1$-robust tangency
between $W^s(\Gamma^1_\varepsilon)$ and
$W^u(\Gamma^2_\varepsilon)$. In this case $\Gamma^1_\varepsilon$
and $\Gamma^2_\varepsilon$ are $double$-blenders for
$\Phi_\varepsilon$. Namely, we can take
\begin{enumerate}[leftmargin=0.6cm]
\item[-] $\Phi_\varepsilon$ having a $\mathcal{S}^1$-robust homoclinic tangency associated with $\Gamma^1_\varepsilon=\Gamma^2_\varepsilon$
(hence $i_1=i_2$);
\item[-] $\Phi_\varepsilon$ having a $\mathcal{S}^1$-robust
equi/heterodimensional tangency on a cycle associated with
$\Gamma^1_\varepsilon$ and $\Gamma^2_\varepsilon$.
\end{enumerate}
\end{thm}

\begin{proof} We will work in local coordinates,
and hence it suffices to construct arcs of one-step maps
$$ \Phi=\tau\ltimes(\phi_1,\dots,\phi_d)\in
\mathcal{PHS}^r(\mathbb{R}^c), \ \ \text{with} \ \ c>1 \ \
\text{and} \ \ r\geq 2$$ satisfying the statement of the theorem.
Set
$$
0<i_1,i_2<c \quad \text{and} \quad \max\{0,i_2-i_1\}< \ell \leq
\min\{c-i_1,i_2\}.
$$

By Proposition~\ref{prop:blending-tangencias} we obtain arcs of
$C^r$-diffeomorphisms
$\phi_1\equiv\phi_1(\varepsilon),\dots,\phi_{d}\equiv
\phi_{d}(\varepsilon)$ homotopic to the identity as $\varepsilon
\to 0^+$, where $d\equiv d(c)\geq 2$ only depends on $c$, and
having $cs$ and $cu$-blending regions $B_1$ and $B_2$ of
$cs$-indices $i_1$ and $i_2$ with $\ell$-tangencies $\hat{B}_1$
and $\hat{B}_2$ respectively. Moreover, the IFS generated by the
induced maps $\hat{\phi}_1,\dots,\hat\phi_d$ on
$\hat{M}=\mathbb{R}^c\times G(\ell,c)$ has transition between
$\hat{B}_1$ and~$\hat{B}_2$.

Set $\Phi_\varepsilon=\tau \ltimes (\phi_1,\dots,\phi_{d}) \in
\mathcal{PHS}^r(M)$. Hence $\Phi_0=\tau\times\mathrm{id}$ and
according to Corollary~\ref{cor:tangencias-IFS}, for any
$\varepsilon>0$, $\Phi_\varepsilon$ has a $\mathcal{S}^1$-robust
tangency of dimension $\ell$ between $W^u(\Gamma^1_\varepsilon)$
and $W^s(\Gamma^2_\varepsilon)$. Here $\Gamma^1_\varepsilon$ and
$\Gamma^2_\varepsilon$ are $cs$ and $cu$-blenders of  $cs$-indices
$i_1$ and $i_2$ associated with the blending regions $B_1$ and
$B_2$ respectively. This completes the first part of the theorem.

The second part follows similarly by asking that the above IFS is
constructed (doing the corresponding modifications in
Proposition~\ref{prop:blending-tangencias} and
Corollary~\ref{cor:tangencias-IFS}) in such a way that the
associated one-step skew-product satisfies the required properties
of Corollary~\ref{cor:tangencias}.
\end{proof}

\section{Realizations: Proofs of Theorems~\ref{thm:cor1}, \ref{thmC} and~\ref{thmF}}
\label{s:realization} First of all, we will construct
diffeomorphisms satisfying the hypothesis of Theorem~\ref{thm:mainA-symbolic-setting}. This
smooth realization of robust cycles will done in a similar manner
as in \cite{BKR14}.

\begin{thm}
\label{thmA} For any integer $0<k\leq \dim M$, there is an arc
$\{f_{\varepsilon}\}_{\epsilon \geq \epsilon}$ of
$C^r$-diffeomorphisms of $N \times M$ such that $f_{0} = F \times
\mathrm{id}$ and for every $\varepsilon>0$, any small enough
$C^1$-perturbation $g$ of $f_{\varepsilon}$ has
\begin{enumerate}[leftmargin=0.6cm]
\item[-] a transitive partially hyperbolic  set $\Delta_g
\subset N \times M$ homeomorphic to $\Lambda \times M$ and
\item[-] a $C^1$-robust
heterodimensional  cycle (in $\Delta_g$) of co-index $k$.
\end{enumerate}
\end{thm}

\begin{proof} 
By hypothesis, $F:N\to N$ has a horseshoe $\Lambda\subset N$,
which is conjugated to a full shift of $d$ symbols. Here $d\geq 2$
is a sufficiently large integer coming from the number of
generators of the IFSs required in the constructions of
Theorem~\ref{thm:mainA-symbolic-setting}, which only depends on
the dimension of $M$. Take $R_1,\ldots, R_{d}$ to be the
rectangles in the manifold $N$ such that $\{R_1\cap\Lambda,
\ldots, R_{d}\cap \Lambda\}$ is a Markov partition for
$F|_\Lambda$.

From Theorem~\ref{thm:mainA-symbolic-setting} we have arcs,
$\Phi_{\varepsilon}=\tau\ltimes(\phi_1,\ldots,\phi_d)\in
\mathcal{PHS}^r(M)$, $\varepsilon\geq 0$, of
$\mathcal{S}^0$-robustly topologically mixing symbolic
skew-products with a $\mathcal{S}$-robust cycle of co-index $0\leq
k\leq c$ where $c=\dim M$. We can modify $f_0=F\times\mathrm{id}$
in $R_i\times M$ to get a one-parameter family of diffeomorphisms
$f_\varepsilon$ satisfying
\begin{align*}
f_{\varepsilon} |_{R_i\times M}&=F \times \phi_i  \quad \text{for $i=1,\dots,d$}
\end{align*}
Note that the locally constant  skew-product diffeomorphism
$f_{\varepsilon}$ restricted to the set $\Lambda\times M$ is
conjugated to the one-step symbolic skew-product
$\Phi_\varepsilon$. In fact, according to \cite[Prop.~A.2]{BKR14}
(see also,~\cite{Go06,IN10,PSW12}), for every small enough
$C^1$-perturbation $g$ of $f_\varepsilon$, $g$ has an invariant
set $\Delta_g$ homeomorphic to $\Lambda\times M$ such that
$g|_{\Delta_g}$ is conjugated with a $\mathcal{S}^0$-perturbation
$\Psi_{g}$ of $\Phi_{\varepsilon}$. The conjugation will associate
with the two hyperbolic sets and their cyclic connections
constructed for $\Psi_{g}$, a cycle for the map $g$ restricted to
$\Delta_g$ of the same co-index. Thus, $f_{\varepsilon}$
restricted to $\Delta=\Lambda\times M$ is $C^1$-robustly
transitive partially hyperbolic and has a $C^1$-robust cycle of
co-index $0\leq k \leq c$. This completes the proof.
\end{proof}

Notice that Theorem~\ref{thm:cor1} and~\ref{thmC} are an immediate
combination of the constructions in Theorems~\ref{thmA}
and~\ref{thmF}. Thus, we only need to show Theorem~\ref{thmF}.

\subsection{Proof of Theorem~\ref{thmF}}
From Theorem~\ref{thm:mainA-symbolic-setting} we have an arc
$\Phi_{\varepsilon}=\tau\ltimes(\phi_1,\ldots,\phi_d)\in
\mathcal{PHS}(M)$, $\varepsilon\geq 0$, of symbolic skew-products
having a $\mathcal{S}^1$-robust tangency between the stable and
unstable sets of a $cs$-blender and a $cu$-blender respectively.
As in the proof of Theorem~\ref{thmA}, take
$f_0=F\times\mathrm{id}$ in $R_i\times M$ to be a one-parameter
family of diffeomorphisms $f_\varepsilon$ satisfying
\begin{align*}
f_{\varepsilon} |_{R_i\times M}&=F \times \phi_i  \quad \text{for $i=1,\dots,d$.}
\end{align*}
For a small $C^2$-perturbation $g$ of $f_\varepsilon$, $g$ has an
invariant set $\Delta_g$ homeomorphic to  $\Lambda\times M$ such
that $g|_{\Delta_g}$ is conjugated with an
$\mathcal{S}^{1}$-perturbation of $\Psi_{g}$ of
$\Psi_{\varepsilon}$ (see~\cite{Go06}). The conjugacy will
associate the intersections between the unstable and stable sets
of the blenders for $\Psi_g$ with intersections between the
unstable and stable manifolds of the corresponding blenders for
$g$. We now would like to show that the tangent directions
associated with the tangency for $\Psi_g$ will correspond to
tangent directions for $g$.

Let $h$ be the homeomorphism $h: \Sigma\times M\rightarrow
\Delta_g\subset N\times M$, that conjugates $\Psi_{g}$ and $g$
restricted to $\Delta_g$. Since we are working in the $C^2$
topology, by \cite{Go06}, restricted to the fibers $h$ is $C^2$.
That is for a fixed sequence $\xi$, the map $h_\xi: M \to
\Delta_g$, given by $h_\xi(x)=h(\xi,x)$, is $C^2$. Moreover by
\cite{IN10} (see~\cite[Equations A.5--A.6]{BKR14}), the map
$h_\xi$ is $C^1$-close to the identity of order $\epsilon$, where
$\epsilon$ depends on the size of the initial arc $f_\varepsilon$,
but does not depend on the sequence $\xi$.

Going back to $\Psi_g=\tau\ltimes \psi_\xi$, let $(\zeta,y)\in
\Sigma\times M$ be the tangency point and consider an unitary
tangencial direction $v\in T_yM$. Hence, there are $C>0$ and
$0<\lambda<1$ such that
$$
\|D\psi_\zeta^n(y)v\|\leq C\lambda^{|n|} \quad \text{for all
$n\in\mathbb{Z}$.}
$$
On the other hand, $Y=h_\zeta(y)=h(\zeta,y)$ is a point in the
intersection between the stable and unstable manifolds of the
blenders for $g$. Consider the vector,
$$
w\eqdef Dh_\zeta(y) v\in T_{Y}(M\times N).
$$
We would like to show that $\|Dg^n(Y)w\|$ goes exponentially fast
to zero as $|n|\to \infty$. This would imply that $w$ is a
tangencial direction at $Y$ for $g$.

Using the conjugation on the fibers, $g^n\circ
h_\zeta=h_{\tau^n(\zeta)}\circ \psi_\zeta^n$, and so
\begin{align*}
\| Dg^n(Y)w)\|&=\|Dg^n(h_\zeta(y))Dh_\zeta(y)v\| \\
            &= \|D(g^n\circ h_\xi(y)) v\| =\|D(h_{\tau^n(\zeta)}\circ \psi_\zeta^n(y)) v\| \\
&\leq \|Dh_{\tau^n(\zeta)}(\psi^n_\zeta(y))\| \cdot
\|D\psi_\zeta^n(y) v\| \leq (1+\epsilon) \|D\psi_\zeta^n(y)v\|\leq
(1+\epsilon)C\lambda^n
\end{align*}
where we have used the fact that $Dh_{\tau^n(\zeta)}$ is
$\epsilon$-close to the identity. This proves that that the
tangency at $Y$ for $g$ has the same number of independent
tangencial directions as the tangency $(\zeta,y)$ for $\Psi_g$,
completing the proof of Theorem~\ref{thmF}.

\begin{rem}
Theorem~\ref{thmF} allows the construction of $C^2$-robust
heterodimensional tangencies on a $C^1$-robust cycle. It is
natural to ask, if there exist open sets of diffeomorphisms which
have heterodimensional tangencies on heterodimensional cycles of a
given co-index $k$. This question was posed in~\cite{KS12}, and
was solved in the case of $k=d-2$ where $d$ is the dimension of
the manifold. Hence, by Theorem~\ref{thmF} and the result
of~\cite{KS12}, every $C^2$-manifold of dimension $d\geq 3$ admits
a $C^2$-robust heterodimensional tangency on a cycle of any
co-index $1\leq k\leq d-4$ and $k=d-2$. The case of co-index
$k=d-3$ remains open.
\end{rem}


\appendix
\addtocontents{toc}{\protect\setcounter{tocdepth}{-1}}
\renewcommand{\thesection}{A}
\renewcommand{\theequation}{A.\arabic{equation}}
\setcounter{equation}{0} \setcounter{thm}{0}
\section{Proof of the criterion for blenders}
\label{AppendixB} In this section we prove Theorem~\ref{cs-thm}.
Since blenders are local dynamical tools we can consider only
local perturbations. For this reason, we need to extend the set of
symbolic skew-products.

  Let $M$ be a separable locally compact metric space
with finite Lebesgue covering dimension and consider $D\subset M$
and $0<\lambda<\beta$. A map $\phi\colon M\to M$ is
\emph{$(\lambda, \beta)$-Lipschitz} on $D$ if
$$
\lambda\,d(x,y)< d(\phi(x),\phi(y)) < \beta \, d(x,y), \quad
\mbox{ for all $x, y \in \overline D$.}
$$

For a  subset $S\subset \mathscr{A}$, set $ \Sigma_S^+ =
\{\xi=(\xi)_{i\in\mathbb{Z}}\in \Sigma: \ \xi_0 \in S
   \}$. We say that $\Phi=\tau\ltimes\phi_\xi$ is \emph{locally
$\alpha$-H\"older continuous} on $D$ at $S$
if there is  $C\ge 0$ such that
\begin{equation}
\label{eq:Holder2} d(\phi^{\pm 1}_\xi(x),\phi^{\pm 1}_\zeta(x))
\leq C \, d_{\Sigma}(\xi,\zeta)^{\alpha} \quad \mbox{for all $x\in
\overline D$ and $ \xi, \zeta \in \Sigma_S^+$ with $\xi_0=\zeta_0$.}
\end{equation}
We denote by $C_0(\Phi,D,S)$ the smallest non-negative constant
satisfying~\eqref{eq:Holder2}.

\begin{defi}
\label{d.thesetSS} Consider $D\subset M$, $0<\lambda<\beta$ and
$S\subset \mathscr{A}$. We define
$\mathcal{S}(D,S,\lambda,\beta)\equiv\mathcal{S}_{\mathscr{A},\nu}^\alpha(D,S,\lambda,\beta)$
as the set of symbolic skew-products $\Phi =\tau \ltimes \phi_\xi$
as in~\eqref{e:sym-skew} such that
\begin{itemize}
\item
$\phi_\xi$ is $(\lambda,\beta)$-Lipschitz on $D$ for all
$\xi\in\Sigma_S^+$,
\item $\phi_\xi$ depends locally $\alpha$-H\"older continuously
on $D$  with respect to $\xi\in\Sigma_S^+$.
\end{itemize}
\end{defi}

Assuming $D$ a relatively compact set, the set
$\mathcal{S}(D,S,\lambda,\beta)$ is endowed with the pseudometric
\begin{equation*}
  d_{\mathcal{S}}(\Phi,\Psi)_{D,S} \eqdef
  d_0(\Phi,\Psi)_{D,S}+\mathrm{Lip}_0(\Phi,\Psi)_{D,S}+
  \mathrm{Hol}_0(\Phi,\Psi)_{D,S},
\end{equation*}
where for $\Phi=\tau\ltimes\phi_\xi$ and $\Psi=\tau\ltimes\psi_\xi$
in $\mathcal{S}(D,S,\lambda,\beta)$
\begin{align*}
   d_0(\Phi, \Psi)_{D,S}\eqdef \sup_{\xi \in \Sigma_S^+} \,
 &d_{C^0(\overline{D})}(\phi^{\pm 1}_\xi,\psi^{\pm 1}_\xi),
 \quad
   \mathrm{Lip}_0(\Phi,\Psi)_{D,S}\eqdef
   \sup_{\xi \in \Sigma_S^+}
   \big|L^{\pm1}(\phi_\xi,D)-L^{\pm 1}(\psi_\xi,D)\big| \\
   &\mathrm{Hol}_0(\Phi,\Psi)_{D,S}\eqdef |C_0(\Phi,D,S)-C_0(\Psi,D,S)|
\end{align*}
with
$$
L(\phi,D)=\max_{x,y\in \overline{D},  x \not = y}
\frac{d(\phi(x),\phi(y))}{d(x,y)} \quad \text{and} \quad
L^{-1}(\phi,D)=\min_{x,y\in \phi(\overline{D}),  x \not = y}
\frac{d(\phi^{-1}(x),\phi^{-1}(y))}{d(x,y)}.
$$

Given $\bar{\omega}=\omega_{-\ell}\ldots
\omega_{-1};\omega_0\,\omega_1\ldots \omega_{j-1}$, where
$\ell,j\geq 0$ and $\omega_i \in \{ 1, \dots, d \}$, we define the
\emph{bi-lateral cylinder} by
\begin{equation*}
\label{n.seqs} \mathsf{C}_{\bar{\omega}} \eqdef \{ \xi \in
\Sigma\colon \xi_{i}= \omega_{i}, \ -\ell \leq i \leq j-1 \}.
\end{equation*}
We denote these bi-lateral cylinders by $\mathsf{V}_{\bar\omega}$
and $\mathsf{H}_{\bar{\omega}}$ in the particular cases  that the
finite word $\bar\omega$ has $\ell>0, j=0$ and $\ell=0, j>0$
respectively. Finally, given a finite set of integers
$\CMcal{I}=\{n_1,\dots,n_k\} \subset \mathbb{N}$ and a set
$\mathcal{B}\subset \mathcal{M}=\Sigma\times M$, we define
$$
  \Lambda^u_{\CMcal{I}}(\mathcal{B};\Phi)
  = \{P\in \mathcal{M}:
  \ \ \text{there is $m_i\to \infty$ so that  $m_{i+1}-m_{i}\in \CMcal{I}$ and $\Phi^{-m_i}(P)\in
  \mathcal{B}$}
  \}.
$$

\begin{thm}
\label{cs-thm1} Let $\Phi \in \mathcal{S}(D,S,\lambda,\beta)$ be a
symbolic skew-product with $\nu^\alpha<\lambda$. Assume that the
following $cs$-covering property holds: there
 exist an integer $k\geq 2$, open sets $B, B_i \subset D$,
integers $n_i \in \mathbb{N}$ and words $\alpha_i \in S^{n_i}$ for
all $i=1,\dots,k$ such that
\begin{equation}
\label{cs-cover1} \mathsf{V}_{\alpha_i}   \times \overline{B_i}
\subset \Phi^{n_i}(\mathsf{H}_{\alpha_i}\times B) \quad \text{for
$i=1,\dots,k$ \ \ and} \quad \overline{B} \subset \bigcup_{i=1}^k
B_i
\end{equation}
$$
  \Omega_B\eqdef \bigcup_{i=1}^k \bigcup_{j=0}^{n_i-1}
\Phi^{j}(\mathsf{H}_{\alpha_i}\times B) \subset \Sigma_S^+\times D
\quad \text{and \quad $C(\Phi,D,S) \cdot
(1-\nu^\alpha\lambda^{-1})^{-1} < L$}
$$
where $L$ is the Lebesgue number of~\eqref{cs-cover1}. Then for
any $0<\delta<\lambda L/2$ and for every small enough
$\mathcal{S}$-perturbation $\Psi$ of $\Phi$,
\begin{align*}
    \Lambda^u_\CMcal{I}(\mathcal{B};\Psi)  \cap \mathcal{D}^s &\not=\emptyset,
\quad \text{for all $\delta$-horizontal discs $\mathcal{D}^s$ in
$\mathcal{B}=\mathsf{V}\times B$}
\end{align*}
where
$$
\CMcal{I}=\{n_1,\dots,n_k\},\quad \text{and} \quad
\mathsf{V}=\mathsf{H}_{\alpha_1}\cup \dots\cup
\mathsf{H}_{\alpha_k}.
$$
In addition, assuming that the maximal invariant set $\Gamma$ in
$\Sigma_S^+\times \overline{D}$ is a transitive hyperbolic set of
$\Phi$ with $\mathrm{ind}^{cs}(\Gamma)>0$, then  $\Gamma$ is a
$cs$-blender of $\Phi$ whose superposition region contains the
open set of almost horizontal disks in $\mathsf{V}\times B$.
\end{thm}

\begin{addendum}
\label{cs-Add1}
Under the assumption of Theorem~\ref{cs-thm1},
if $n_i=\ell\geq 1$ for all $i=1,\dots,k$
then for every small enough $\mathcal{S}$-perturbation
$\Psi$ of $\Phi$,
$$
    \Lambda^u(\mathcal{B};\Psi^\ell)  \cap \mathcal{D}^s \not=\emptyset
\ \ \text{for all $\delta$-horizontal discs $\mathcal{D}^s$ in
$\mathcal{B}$ where \ \ $\Lambda^u(\mathcal{B};\Psi^\ell) \eqdef
\bigcap_{n\geq 0} \Psi^{\ell n}(\mathcal{B})$.}
$$
In addition, if the maximal invariant set $\Gamma$
 in $\Sigma_S^+\times \overline{D}$ is a transitive hyperbolic set of
$\Phi^\ell$ with $\mathrm{ind}^{cs}(\Gamma)>0$, then  $\Gamma$ is
a $cs$-blender of $\Phi^\ell$.  
\end{addendum}
Analogously conditions yield $cu$-blenders and $double$-blenders
(see Remark~\ref{rem:cu-blender}).

\begin{rem}
In another approach, Moreira and Silva~\cite{MS12} used the criterion of the recurrent compact set to obtain blenders. Their techniques
have certain similarities to the proof below, and in fact the covering property~\eqref{cs-cover1} implies that the blender-horseshoe constructed satisfies the recurrent
compact criterion. Having in mind the recurrent compact set criterion, a new dynamical definition of a blender was given in~\cite{BBD15}.
The different relationships between these concepts is an intersting problem.
\end{rem}

First of all notice that Theorem~\ref{cs-thm} follows from the
above result. In fact, similar as it was argued
in~\S\ref{s:blender}, since $\Lambda^u_\CMcal{I}(\mathcal{B};\Psi)
\subset W^u_{loc}(\Gamma)$ the second part of
Theorem~\ref{cs-thm1} is an immediate consequence of the first
part. We will now split the proof into three steps.

\subsection{Main Lemma} Given  $\zeta \in \Sigma$ and
$\bar{\omega} =\omega_{-\ell}\ldots \omega_{-1}$, where $\ell\ge
1$ and $\omega_i \in \mathscr{A}$, we define the {\emph{relative
cylinder}} by
$$
\mathsf{C}_{\bar\omega}(\zeta) \eqdef
\mathcal{V}_{\bar{\omega}}\cap W^s_{loc}(\zeta)
$$
and recall the notation~\eqref{n.seq},
$$
\psi^n_\zeta(x) \eqdef \psi^{}_{\tau^{n-1}(\zeta)}\circ\cdots\circ\psi^{}_{\zeta}(x) 
    \quad \text{and} \quad
    \psi^{-n}_\zeta(x) \eqdef \psi^{-1}_{\tau^{-n}(\zeta)}\circ\cdots
    \circ\psi^{-1}_{\tau^{-1}(\zeta)}(x).
$$
\begin{lem}[{\cite{BKR14}}]
\label{l.noholdas} Consider  $\Psi=\tau\ltimes \psi_\xi \in
\mathcal{S}(D,S,\lambda,\beta)$,
a word $\bar{\omega}=\omega_{-n}\ldots \omega_{-1};\omega_0 \ldots
\omega_{n}$ with $\omega_{-i}\in S$ for all $i>0$ and a point
$x\in \overline{D}$ such that
$$
\psi^{-j}_{\zeta}(x) \in \overline{D} \quad \text{for all
$\zeta\in \mathsf{C}_{\bar{\omega}}$  and  $1\leq j\leq n$.}
$$
Then for every $1\leq i\leq n$ and $\xi,\zeta \in
\mathsf{C}_{\bar{\omega}}$,
\begin{align*}
d(\psi^{-i}_{\xi}(x),\, \psi^{-i}_{\zeta}(x)) \leq
 C(\Psi,D,S) \, \nu^{-\alpha i}
\, \sum_{j=0}^{i-1} (\lambda^{-1} \nu^{\alpha})^{j} \,
d_{\Sigma}(\xi,\zeta)^\alpha.
\end{align*}
\end{lem}

\subsection{Choice of the $\mathcal{S}$-neighborhood.} Observe
that~\eqref{cs-cover} is robust under $\mathcal{S}$-perturbations.
%
Thus, we take  a $\mathcal{S}$-neighborhood $\mathscr{U}$ of
$\Phi$ such
that if $\Psi \in \mathscr{U}$, then~\eqref{cs-cover} holds for $\Psi$. 
In particular,
\begin{align*}
\label{rem-per} \psi^{-j}_{\xi}(\overline{B_i}) &\subset D   \quad
\, \text{for all $j=0,\dots,n_i-1$ and \ \ }
\psi^{-n_i}_{\xi}(\overline{B_i}) \subset B \quad \text{for all \
}
   \xi \in \mathsf{V}_{\alpha_i}.
\end{align*}
By hypothesis, shrinking the neighborhood $\mathscr{U}$ if
necessary, we can assume that
$$C(\Psi,D,S)\cdot (1-\nu^\alpha\lambda^{-1})^{-1}<L \ \text{for all} \ \Psi \in \mathscr{U}.$$

\subsection{Existence of an intersection
point} The following proposition provides the final step.

\begin{prop}\label{p.jairo}
Consider $0<\delta <\lambda L/2$ and let
$\mathcal{D}^s=\mathcal{D}^s(\zeta,z)$ be a $\delta$-horizontal
disc in $\mathcal{B}=\mathsf{V}\times B$ with $\alpha$-H\"older
constant $C\geq 0$. Then for every $\Psi=\tau\ltimes \psi_\xi \in
\mathscr{U}$ there are an infinite word
$\bar{\omega}=\ldots\bar\omega_{-j}\ldots\bar\omega_{-1}$ with
$\bar\omega_{-j}= \alpha_{i_j}$ and a sequence of nested compacts
subsets
$\{V_n\}$ of $M$  such that 
\begin{enumerate}[itemsep=0.1cm]
\item $V_n \subset \mathscr{P}
(\mathcal{D}^s\cap(\mathsf{C}_{\bar{\omega}^n}(\zeta)\times B))$,
\item $\Psi^{-j}(\mathsf{C}_{\bar{\omega}^n}(\zeta)\times V_n)
\subset \Omega_B \quad \text{for $j=0,\dots,m_n-1$,}$
\item $\Psi^{-m_n}(\mathsf{C}_{\bar{\omega}^n}(\zeta)\times V_n)
\subset \Omega_B\cap\mathcal{B}$,
\item $\mathrm{diam}(
\psi_{\xi}^{-m_n}(V_n)) \leq C(\lambda^{-1}\nu^\alpha)^{m_n}$ for
all $\xi \in \mathsf{C}_{\bar{\omega}^{n}}(\zeta)$,
\end{enumerate}
where $\bar\omega^n=\bar{\omega}_{-n}\ldots\bar{\omega}_{-1}$,
$m_n$ is the length of the word $\bar{\omega}^n$, and
$\mathscr{P}:\Sigma \times M \to M$ is the projection on  the
fiber space.
\end{prop}

This proposition concludes the first part of
Theorem~\ref{cs-thm1}. Indeed, let
$$
\{x\} =\bigcap_{n\in\mathbb{N}} V_n \subset B \quad \text{and}
\quad  \{\xi\}=\bigcap_{n\in\mathbb{N}}
\mathsf{C}_{\bar{\omega}^{n}}(\zeta)\subset W^s_{loc}(\zeta).
$$
Observe that $(\xi,x)\in \mathcal{D}^s$ and $\Psi^{-m_n}(\xi,x)\in
\Omega_B\cap \mathcal{B}$ for all $n\geq 1$ and $m_{n+1}-m_{n}\in
\CMcal{I}\eqdef \{n_1,\dots,n_k\}$. Hence
$$
(\xi,x)\in \Lambda^u_{\CMcal{I}}(\mathcal{B};\Psi)\cap \mathcal{D}^s.
$$ 

\begin{proof}[Proof of Proposition~\ref{p.jairo}]
Fix a skew-product $\Psi=\tau\ltimes \psi_\xi \in\mathscr{U}$ and
consider the $(\alpha,C)$-H\"older map
\mbox{$h:W^s_{loc}(\zeta;\tau) \to B$} associated with the
$\delta$-horizontal disc
$\mathcal{D}^s=\mathcal{D}^s(\zeta,z)\subset \mathsf{V}\times B$.
The construction of the nested sequence of sets $\{V_n\}$ and the
infinite word
$\bar{\omega}=\ldots\bar\omega_{-j}\ldots\bar\omega_{-1}$ with
$\bar\omega_{-j}= \alpha_{i_j}$ is done inductively. Let $V =
\mathscr{P}(\mathcal{D}^s)\subset B$.
Note that $\mathrm{diam}(V)\leq 2\delta<L$.  By the definition of
the Lebesgue number, we have that
 $V\subset B_{i_1}$ for some $i_1\in \{1,\dots,k\}$.
Consider $$
\text{$\bar{\omega}^{1}=\bar{\omega}_{-1}=\alpha_{i_1}$, \quad
$m_1=n_{i_1}$ \quad \text{and} \quad $V_1\eqdef\mathscr{P} (
\mathcal{D}^s\cap(\mathsf{C}_{\bar{\omega}^{1}}(\zeta) \times
V))$}.
$$
By construction, $V_1 \subset V$ and
$$
\psi_{\xi}^{-j}(V_1) \subset D  \quad \text{for all
$j=0,\dots,m_1-1$} \ \ \text{and} \ \ \psi_{\xi}^{-m_1}(V_1)
\subset B \quad \text{for all $\xi\in
\mathsf{C}_{\bar\omega^{1}}(\zeta)$.}
$$
\begin{claim}
\label{cl.diametro} $\mathrm{diam}(V_1)\leq  C\,\nu^{m_1\alpha}$.
\end{claim}

\begin{proof}
Given $x$ and $y$ in $V_1$ there are $\xi$ and $\eta$ in
$\mathsf{C}_{\bar{\omega}^{1}}(\zeta)$ such that $x=h(\xi)$ and
$y=h(\eta)$. Since $h$ is $(\alpha, C)$-H\"older continuous, $
d(x,y) =d(h(\xi),h(\eta)) \leq C d_{\Sigma}(\xi,\eta)^\alpha \leq
C \, \nu^{m_1\alpha}$, proving the claim.
\end{proof}

By Claim~\ref{cl.diametro}, setting $\delta_1\eqdef
C\,\nu^{m_1\alpha}$ and since the fiber-maps $\psi_\xi$ are
$(\lambda,\beta)$-Lipschitz on $\overline{D}$,
$$
\mathrm{diam} ( \psi_{\xi}^{-n_1}(V_1)  ) \leq
\lambda^{-m_1}\delta_1 \quad \mbox{for all $\xi\in
\mathsf{C}_{\bar{\omega}^{1}} (\zeta)$.}
$$
Since $\lambda^{-1}\nu^\alpha <1$ and $\mathcal{D}^s$ is a
$\delta$-horizontal disc ($C \nu^{\alpha}< \delta$), it follows that $
\lambda^{-m_1}\delta_1 = C(\lambda^{-1} \nu^{\alpha} )^{m_1}<
\lambda^{-1} \delta\leq L/2$. Therefore, the diameter of
$\psi_{\xi}^{-1}(V_1)$ is less or equal than $L/2$.

Arguing inductively, suppose that we have constructed a finite
word $\bar{\omega}^{n}\eqdef \bar\omega_{-n} \dots
\bar\omega_{-1}$ (the word $\bar{\omega}^{j}$ is obtained by adding
$\bar\omega_{-j}=\alpha_{i_j}$ to the word
$\bar{\omega}^{j-1}$). As well we have the closed sets $V_n\subset V_{n-1} \subset
\dots \subset V_1$ with diam$(V_n)\leq \delta_n=C\nu^{m_n\alpha}$, where $m_n$ is the length of the word
$\bar\omega^n$, such that
\begin{equation}
\label{e.hipo} \psi_{\xi}^{-j}(V_n) \subset D \quad \text{for all
$j=0,\dots,m_n-1$} \text{,} \ \psi_{\xi}^{-m_n}(V_n)
\subset B
\end{equation}
and
\begin{equation}
\label{e.hipo2} \mathrm{diam} ( \psi_{\xi}^{-m_n}(V_n)) \leq
\lambda^{-m_n}\delta_n \quad \text{for all
$\xi\in\mathsf{C}_{\bar{\omega}_{-n}}(\zeta)$}.
\end{equation}
We now construct the word $\bar{\omega}^{n+1}$ and the closed set
$V_{n+1} \subset V_n$ satisfying analogous inclusions and
inequalities. By \eqref{e.hipo} and~\eqref{e.hipo2}  we have that
$$
A_n\eqdef \bigcup_{\xi \,\in \,
\mathsf{C}_{\bar{\omega}^{n}}(\zeta)   } \psi_{\xi}^{-m_n}(V_n)
\subset B.
$$

\begin{claim}\label{cl.semnome}
 $\mathrm{diam}(A_n)<L$.
\end{claim}

\begin{proof}
Given $\bar{x}$ and $\bar{y}$ in $A_n$ there are $x, y \in V_n$
and $\xi, \eta \in \mathsf{C}_{\bar{\omega}_{-n}} (\zeta)$ such
that $\bar{x}=\psi^{-m_n}_{\xi}(x)$ and
$\bar{y}=\psi^{-m_n}_{\eta}(y)$. Then
\begin{align*}
d(\bar{x},\bar{y}) = d(\psi^{-m_n}_{\xi}(x),\psi^{-m_n}_{\eta}(y))
\leq C(\Psi,D,S) \nu^{-\alpha m_n} \sum_{j=0}^{m_n-1}
(\lambda^{-1}_{cs}\nu^{\alpha})^j \,d_{\Sigma}(\xi,\eta)^\alpha +
\lambda^{-m_n}\delta_n
\end{align*}
where the last inequality follows from Lemma~\ref{l.noholdas}
and the induction hypothesis~\eqref{e.hipo}. Since $\xi$ and
$\eta$ belong to $\mathsf{C}_{\bar{\omega}^{n}} (\zeta)$ then
$d_{\Sigma}(\xi,\eta)^\alpha \leq \nu^{\alpha m_n}$. Hence
$$
  d(\bar{x},\bar{y})\leq C(\Psi,D,S)  \sum_{j=0}^{m_n-1}
(\lambda^{-1}\nu^{\alpha})^j + \lambda^{-m_n} \delta_n \leq L/2
+\lambda^{-m_n} \delta_n
$$
where the last inequality follows from the fact that
$C(\Psi,D,S)\cdot(1-\nu^\alpha\lambda^{-1})<L$. As
$$
\lambda^{-m_n} \, \delta_n = C\, (\lambda^{-1}\nu^{\alpha})^{m_n}
\leq C\, \lambda^{-1}\nu^{\alpha} < \lambda^{-1} \delta < L/2,
$$
therefore $d(\bar x,\bar y)< L$ and thus $\mathrm{diam}(A_n)<L$,
proving the claim.
\end{proof}

Since $L$ is a Lebesgue number of the covering $\{B_{i}\}$, the
claim implies there is $i_{n+1}\in \{1,\dots,k\}$ such that
$\mathscr{P}(A_n)\subset B_{i_{n+1}}$. Set
$$
\bar{\omega}^{n+1}= \alpha_{i_{n+1}}\omega_{-n} \dots \omega_{-1}
\quad \mbox{and} \quad V_{n+1}=\mathscr{P}(
\mathcal{D}^s\cap(\mathsf{C}_{\bar{\omega}^{n+1}}(\zeta) \times
V_n)).
$$
By the construction $V_{n+1}\subset V_n$, and arguing
as in Claim~\ref{cl.diametro}
$$
\mathrm{diam}(V_{n+1})\leq C\,\nu^{m_{n+1}\,\alpha}
\eqdef\delta_{n+1}
$$
where $m_{n+1}$ is the length of $\bar\omega^n$, that is
$m_{n}+\ell$ where $\ell=n_{i_{n+1}}$ is the length of the word
$\alpha_{i_{n+1}}$. Hence,
$$
\mathrm{diam} (\psi_{\xi}^{-m_{n+1}}(V_{n+1}) ) \leq
\lambda^{-m_{n+1}}\delta_{n+1} \quad \text{for all $\xi\in
\mathsf{C}_{\bar{\omega}^{n+1}}(\zeta)$.}
$$

Since $V_{n+1} \subset V_{n}$ and $A_n \subset B_{i_{n+1}}$, it
follows that  for every $\xi\in
\mathsf{C}_{\bar{\omega}^{n+1}}(\zeta)$,
\begin{align*}
\psi_{\xi}^{-m_{n}-j}(V_{n+1}) \subset
\psi_{\tau^{-m_n}(\xi)}^{-j} \circ \psi_{\xi}^{-m_n}(V_{n})
\subset \psi_{\tau^{-m_n}(\xi)}^{-j}(A_n) \subset
\psi_{\tau^{-m_n}(\xi)}^{-j}(B_{i_{n+1}}) \subset D
\end{align*}
for all $j=0,\dots, \ell -1$ and
\begin{align*}
\psi_{\xi}^{-m_{n+1}}(V_{n+1}) \subset
\psi_{\tau^{-m_n}(\xi)}^{-\ell} \circ \psi_{\xi}^{-m_n}(V_{n})
\subset \psi_{\tau^{-m_n}(\xi)}^{-\ell}(A_n) \subset
\psi_{\tau^{-m_n}(\xi)}^{-\ell}(B_{i_{n+1}}) \subset B.
\end{align*}

Similarly as above the diameter of these sets is less or equal
than $\lambda^{-(n+1)}\delta_{n+1}$. Therefore~\eqref{e.hipo}
and~\eqref{e.hipo2} hold for the $(n+1)$th-step and we can continue
arguing inductively. This completes the construction of the
infinite word $\bar\omega$ and the sequence of nested sets
$\{V_n\}$ in the proposition, ending the proof.
\end{proof}
\bibliographystyle{abbrv}
\bibliography{br-bib}

\end{document}